\documentclass[reqno,12pt]{amsart}

\usepackage{amsmath,amsthm,amssymb,amsfonts,amscd}
\usepackage[all]{xy} 
\usepackage[pdftex]{graphicx} 
\numberwithin{equation}{section}

\theoremstyle{plain} 
	\newtheorem{thm}{Theorem}[section]
	\newtheorem*{thm*}{Theorem}
	\newtheorem{cor}[thm]{Corollary}
	\newtheorem{lem}[thm]{Lemma}
	
	\newtheorem{prop}[thm]{Proposition}
	\newtheorem{conj}[thm]{Conjecture}
	\newtheorem*{conj*}{Conjecture}
	
\theoremstyle{definition}
	\newtheorem{defn}[thm]{Definition}

\theoremstyle{remark}
	\newtheorem{rem}[thm]{Remark}
	
	\newtheorem*{pf}{Proof}
\def\EE{{\mathbb E}}
\def\CC{{\mathbb C}}

\def\HH{{\mathbb H}}
\def\LL{{\mathbb L}}

\def\PP{{\mathbb P}}
\def\QQ{{\mathbb Q}}
\def\RR{{\mathbb R}}

\def\ZZ{{\mathbb Z}}

\def\B{{\mathcal B}}

\def\D{{\mathcal D}}
\def\E{{\mathcal E}}

\def\H{{\mathcal H}}
\def\I{{\mathcal I}}

\def\N{{\mathcal N}}
\def\O{{\mathcal O}}

\def\Q{{\mathcal Q}}

\def\S{{\mathcal S}}
\def\T{{\mathcal T}}

\def\X{{\mathcal X}}
\def\Z{{\mathcal Z}}
\def\ee{{\mathbf e}}
\def\p{{\partial }}

\def\ch{{\rm ch}}
\def\id{{\rm id}}

\newcommand{\Jac}{\mathop{\rm Jac}}
\def\ns{{\nabla}\hspace{-1.4mm}\raisebox{0.3mm}{\text{\footnotesize{\bf /}}}}

\def \mf#1#2#3#4{
\xymatrix{
{#1}\  \ar@<0.4ex>[r]^{{#2}} & \ {#4}
\ar@<0.4ex>[l]^{{#3}}
}
}

\def \mfs#1#2#3#4{\!
\xymatrix@C=1,5em{{#1} \! \ar@<0.2ex>[r]^{{#2}} & \! {#4}
\ar@<0.2ex>[l]^{{#3}}
}
\!}

\def \mfl#1#2#3#4{
\xymatrix@C=2.6em{{#1}\  \ar@<0.4ex>[r]^{{#2}} &\  {#4}
\ar@<0.2ex>[l]^{{#3}}
}
}

\def \mfss#1#2#3#4{\!
\xymatrix@C=1.5em{{#1} \ar@<0.3ex>[r]^{{#2}} & {#4}
\ar@<0.3ex>[l]^{{#3}}
}
\!}
\begin{document}
\title{Gamma integral structure for an invertible polynomial of chain type}
\date{\today}
\author{Takumi Otani}
\address{Department of Mathematics, Graduate School of Science, Osaka University, 
Toyonaka Osaka, 560-0043, Japan}
\email{u930458f@ecs.osaka-u.ac.jp}
\author{Atsushi Takahashi}
\address{Department of Mathematics, Graduate School of Science, Osaka University, 
Toyonaka Osaka, 560-0043, Japan}
\email{takahashi@math.sci.osaka-u.ac.jp}

\begin{abstract}
The notion of the Gamma integral structure for the quantum cohomology of an algebraic variety was introduced by Iritani, Katzarkov--Kontsevich--Pantev. 
In this paper, we define the Gamma integral structure for an invertible polynomial of chain type. 
Based on the $\Gamma$-conjecture by Iritani, we prove that the Gamma integral structure is identified with the natural integral structure for the  Berglund--H\"{u}bsch transposed polynomial by the mirror isomorphism.
\end{abstract}
\maketitle
\section{Introduction}
For a holomorphic function $f:\CC^n\longrightarrow\CC$ with at most an isolated critical point at the origin,
K. Saito developed a theory of primitive forms which yields a Frobenius structure on the base space of the universal unfolding \cite{S-K,ST}. 
In order further to study the (exponential) periods of a primitive form, one needs careful analysis of the structure connections of the Frobenius manifold, especially, their integral structures.
A Frobenius manifold can be associated to a smooth projective variety (or orbifold) $X$ by the genus zero Gromov--Witten invariants of $X$.
The first structure connection of the Frobenius manifold $M_X$, often called the quantum connection in this setting, is a flat meromorphic connection on $\PP^1\times M_X$.
If the ring structure on the tangent space $T_p M_X$ over a point $p\in M_X$ is semi-simple, then the meromorphic connection restricted to $\PP^1=\PP^1\times\{p\}$ has an irregular singular point and a regular singular point.
One can define the monodromy data (at the point $p$) of the Frobenius manifold $M_X$, consisting of the monodromy matrix at the regular singular point, the Stokes matrix $S$ at the irregular singular point, and the central connection matrix $C$ between these two singular points. 
Here, it is very important to note that the monodromy data is locally constant \cite[Isomonodromicity Theorem]{D3}. 
In \cite{D2}, Dubrovin formulated a conjecture on a close relationship between the monodromy data of the Frobenius manifold $M_X$ and the structure of the bounded derived category $\D^b(X)$ of coherent sheaves on $X$ when $X$ is Fano. Its refined version is given in \cite{CDG} (see Conjecture \ref{conj : refined Dubrovin} below). 
The notion of the {\em Gamma integral structure} on the space of multi-valued flat sections of the quantum connection was introduced 
by Iritani \cite{I} (see also~\cite{KKP}). 
This structure is defined to be a framing induced by the Grothendieck group of the derived category $\D^b(X)$. 
To be more precise, suppose that $\D^b(X)$ admits a full exceptional collection $(\E_1,\dots,\E_{\widetilde \mu})$.
Then the framing is given by the free $\ZZ$-module spanned by $\widehat\Gamma_X {\rm Ch}(\E_j)$, $j=1,\dots, \widetilde \mu$, 
where $\widehat\Gamma_X$ is the Gamma class of $X$ and ${\rm Ch}(\E_j)$ is the (modified) Chern character of $\E_j$. 
For a weak Fano orbifold $X$ such that $\D^b(X)$ admits a full exceptional collection, Galkin--Golyshev--Iritani \cite{GGI} conjectured that 
the Gamma integral structure describes the central connection matrix $C$ of the Frobenius manifold $M_X$. 
This conjecture is called the $\Gamma$-conjecture, which is equivalent to a part of the refined Dubrovin's conjecture (see~\cite{I,CDG}). 
{\em Mirror symmetry} is an equivalence of seemingly different two objects in mathematics, algebraic one and geometric one,
which was discovered by physicists in their study of string theory. Consider a polynomial $f_n\in S_n:=\CC[z_1,\dots, z_n]$ of the form 
\begin{equation*}
f_n= f_n(z_1,\dots, z_n):= z_1^{a_1}z_2 + z_2^{a_2}z_3 + \cdots + z_{n-1}^{a_{n-1}}z_n + z_n^{a_n}, \quad a_i\ge 2, 
\end{equation*}
which is called an {\em invertible polynomial of chain type}. 
Let 
\begin{equation*}
G_{f_n} := \left\{ (\lambda_1, \ldots , \lambda_n )\in(\CC^\ast)^n \, \left| \, f_n(\lambda_1z_1, \ldots, \lambda_n z_n) = 
f_n(z_1, \ldots, z_n) \right\} \right.,
\end{equation*}
which is called the group of maximal diagonal symmetries of $f$.
It is expected by Berglund--H\"{u}bsch \cite{BH} that a mirror dual object corresponding to the pair $(f_n,G_{f_n})$ is given by the polynomial $\widetilde f_n\in\CC[x_1,\dots, x_n]$, 
called the {\em Berglund--H\"{u}bsch transpose} of the polynomial $f_n$, defined by 
\begin{equation*}
\widetilde f_n= \widetilde f_n (x_1,\dots, x_n):= x_1^{a_1} + x_1x_2^{a_2} + \cdots + x_{n-1}x_n^{a_n}. 
\end{equation*}
For the pair $(f_n,G_{f_n})$, one can associate a triangulated category ${\rm HMF}^{L_{f_n}}_{S_n}({f_n})$ of maximally graded matrix factorizations, which is considered as an analogue of 
$\D^b(X)$ for a smooth projective variety. 
On the other hand, a (counter-clockwise)  distinguished basis of vanishing cycles in the Milnor fiber of $\widetilde f_n$ can be categorified to an $A_\infty$-category ${\rm Fuk}^{\to}(\widetilde f_n)$ called the directed Fukaya category, whose derived category $\D^b{\rm Fuk}^\to(\widetilde f_n)$ is an invariant of the polynomial $\widetilde f_n$ as a triangulated category. 
Then the {\em homological mirror symmetry conjecture} for invertible polynomials predicts an equivalence of the category ${\rm HMF}^{L_{f_n}}_{S_n}(f_n)$ and the derived directed Fukaya category $\D^b{\rm Fuk}^\to(\widetilde f_n)$ (cf.~\cite{ET,T2}). 
There are many evidences of the above conjecture which follow from related results by several authors 
(cf. \cite{AT,FU1,FU2,Hab,KST,KST2,Kra,LP,Sei,T,T2,Ue}). 
It is known that from the pair $(f,G_f)$ of an invertible polynomial $f$ and the group of maximal diagonal symmetries $G_f$ a Frobenius manifold is constructed by the FJRW theory \cite{FJR} which is analogous to the Gromov-Witten theory in many ways. 
On the other side, one can associate a Frobenius manifold to the Berglund--H\"{u}bsch transpose $\widetilde f$ of $f$ by a choice of a primitive form (cf. \cite{S-K, ST}). 
The {\em classical mirror symmetry conjecture} for invertible polynomials predicts an isomorphism of these Frobenius manifolds, which 
is proven for many cases (cf. \cite{FJR,HLSW,KrSh,LLSS,MS}). 
One can define $\QQ$-graded $\CC$-vector spaces $\Omega_{f_n,G_{f_n}}$ and $\Omega_{\widetilde f_n}$ for $(f_n,G_{f_n})$ and $\widetilde f_n$ respectively as an analogue of the total Hodge cohomology group $\oplus_{p,q}H^q(X,\Omega_X^p)$ of a smooth projective variety $X$ where two $\QQ$-gradings are introduced by {\em exponents} (see Definition \ref{defn : 2.5}).
We can associate a $\CC$-bilinear form $J_{f_n,G_{f_n}}$ (resp., $J_{\widetilde f_n}$) on $\Omega_{f_n,G_{f_n}}$ (resp., $\Omega_{\widetilde f_n}$). 
It is known that there exists an isomorphism $\Omega_{\widetilde f_n}\cong\Omega_{f_n,G_{f_n}}$ of $\QQ$-graded $\CC$-vector spaces called the {\em mirror isomorphism} (cf. \cite{K,EGT}). 
We choose a mirror isomorphism ${\bf mir}:\Omega_{\widetilde f_n}\cong\Omega_{f_n,G_{f_n}}$ so that it is compatible with the $\CC$-bilinear forms $J_{f_n,G_{f_n}}$ and $J_{\widetilde f_n}$ (Proposition \ref{prop : mirror isomorphism}).
Therefore, we shall denote later by $\eta^{(n)}$ the matrix representation of the bilinear form $J_{f_n,G_{f_n}}$ with respect to a specified homogeneous basis on $\Omega_{f_n,G_{f_n}}$,
which is equal to the matrix representation  of the bilinear form $J_{\widetilde f_n}$ with respect to a specified homogeneous basis on $\Omega_{\widetilde f_n}$. 
The homogeneous basis defines a diagonal matrix $\widetilde Q^{(n)}$ whose entries are exponents of $\widetilde f_n$ shifted by $-n/2$. 
In \cite{AT}, Aramaki--Takahashi showed the existence of a full exceptional collection $(E_1,\dots,E_{\widetilde\mu_n})$ on ${\rm HMF}^{L_{f_n}}_{S_n}(f_n)$ which gives a Lefschetz decomposition on ${\rm HMF}^{L_{f_n}}_{S_n}(f_n)$ in the sense of \cite{KuS} where $\widetilde \mu_n$ is the Milnor number of $\widetilde f_n$. 
Based on this result, we associate a matrix ${\rm ch}^{(n)}_\Gamma$ to the full exceptional collection whose $j$-th column is an analogue of 
$\widehat\Gamma_X {\rm Ch}(\E_j)$ for a Fano orbifold $X$ (see Definition \ref{defn : ch Gamma}). 
The following theorem proves the necessary conditions for the matrix $(2\pi)^{-\frac{n}{2}}{\rm ch}^{(n)}_\Gamma$ to give a central connection matrix of a Frobenius manifold whose non-degenerate 
bilinear form on the tangent space, the grading matrix and the Stokes matrix are given by $\eta^{(n)}$, $\widetilde Q^{(n)}$ and $\chi^{(n)}$, respectively:
\begin{thm}[Theorem~\ref{thm : main theorem}]\label{thm : intro}
The matrix ${\rm ch}^{(n)}_\Gamma$ satisfies following equalities: 
\begin{equation*}
\left(\dfrac{1}{(2\pi)^{\frac{n}{2}}}\ch_{\Gamma}^{(n)}\right)^{-1}\ee\left[\widetilde Q^{(n)}\right]\left(\dfrac{1}{(2\pi)^{\frac{n}{2}}}\ch_{\Gamma}^{(n)}\right)={\bf S}^{(n)},
\end{equation*}
\begin{equation*}
\left(\dfrac{1}{(2\pi)^{\frac{n}{2}}}\ch_{\Gamma}^{(n)}\right)^T\ee\left[\dfrac{1}{2}\widetilde Q^{(n)}\right]\eta^{(n)}\left(\dfrac{1}{(2\pi)^{\frac{n}{2}}}\ch_{\Gamma}^{(n)}\right)=\chi^{(n)}, 
\end{equation*}
where ${\bf S}^{(n)}$ is the matrix representation of the automorphism with respect to the exceptional basis $\{[E_i]\}_{i=1}^{\widetilde\mu_n}$ induced by the Serre functor on ${\rm HMF}^{L_{f_n}}_{S_n}(f_n)$ and $\chi^{(n)}$ is the Euler matrix with respect to the exceptional basis. 
\end{thm}
In the present paper, we study the Gamma integral structure for $(f,G_f)$.
Theorem~\ref{thm : intro} enables us to define the Gamma integral structure $\Omega_{f_n,G_{f_n};\ZZ}$ of $\Omega_{f_n,G_{f_n}}$ as the framing of a map defined by the matrix 
$(2\pi\sqrt{-1})^{-n}{\rm ch}^{(n)}_\Gamma$ with respect to a certain basis of $\Omega_{f_n,G_{f_n}}$ (see Definition \ref{defn : Gamma integral structure}). 
On the other hand, there is a natural integral structure on $\Omega_{\widetilde f_n}$ given by the relative homology group $H_n(\CC^n,{\rm Re}(\widetilde f_n)\gg0;\ZZ)$ via 
the isomorphism $H^n(\CC^n,{\rm Re}(\widetilde f_n)\gg0;\CC)\cong \HH^n(\Omega^\bullet_{\CC^n}, d-d\widetilde f_n\wedge)\cong\Omega_{\widetilde f_n}$.
It is known that the middle homology group $H_{n-1}(\widetilde f_n^{-1}(1);\ZZ)$ of the Milnor fiber of $\widetilde f_n$ is a free abelian group of rank $\widetilde\mu_n$ 
generated by vanishing cycles and there is an isomorphism $H_n(\CC^n,{\rm Re}(\widetilde f_n)\gg0;\ZZ)\cong H_n(\CC^n,\widetilde f_n^{-1}(1);\ZZ)\cong H_{n-1}(\widetilde f_n^{-1}(1);\ZZ)$. 
Therefore, the integral structure $\Omega_{\widetilde f_n;\ZZ}$ is also considered as the image of the free $\ZZ$-module $H_{n-1}(\widetilde f_n^{-1}(1);\ZZ)$ 
in $\Omega_{\widetilde f_n}$. 
Now, we can state the following theorem as an analogue of $\Gamma$-conjecture: 
\begin{thm}[Theorem~\ref{thm : Integral structure}]\label{thm : intro 1}
The mirror isomorphism ${\bf mir}: \Omega_{\widetilde f_n}\cong \Omega_{f_n,G_{f_n}}$ can be chosen so that it induces an isomorphism of free $\ZZ$-modules $\Omega_{\widetilde f_n;\ZZ}\cong \Omega_{f_n,G_{f_n};\ZZ}$ and it is compatible with the $\ZZ/d_n\ZZ$-action on both hands side. 
Here, $\ZZ/d_n\ZZ$-action on $\Omega_{\widetilde f_n;\ZZ}$ is induced by its action on $\CC^n$ given by $x_i\mapsto \ee[(-1)^{i-1}/d_i]x_i$, $i=1,\dots, n$ where $d_i:=a_1\cdots a_i$, $\ee[-]=e^{2\pi\sqrt{-1}(-)}$ and the action on $\Omega_{f_n,G_{f_n};\ZZ}$ is induced by the grading shift $(\vec{z}_1)$ on ${\rm HMF}^{L_{f_n}}_{S_n}(f_n)$. 
\end{thm}
Milanov--Zha \cite{MZ} worked on the Gamma integral structure for the case of type ADE, which motivates us to study the one for general invertible polynomials of chain type.
Indeed, their result shows that the theorem holds for $A_l$, $D_l$ and $E_7$ cases (and also $E_6$ and $E_8$ cases after generalizing suitably the conjecture for 
Thom--Sebastiani sum of invertible polynomials of chain type). 
It is an important problem to show the existence of Bridgeland stability conditions~\cite{B} on a given triangulated category. 
In \cite{T,KST}, it is shown that there is a stability condition on ${\rm HMF}^{L_{f}}_{S}(f)$ for an invertible polynomial of type ADE and, 
in terms of the derived directed Fukaya category $\D^b{\rm Fuk}^\to(\widetilde f)\cong {\rm HMF}^{L_{f}}_{S}(f)$ (see \cite{Sei} for the equivalence $\D^b{\rm Fuk}^\to(\widetilde f)\cong\D^b(\CC\vec{\Delta})$ where $\vec{\Delta}$ is the Dynkin quiver of the corresponding type), 
it is expected that the oscillatory integral for a primitive form induces a stability condition (cf. \cite{B2,T}). 
The view point of the homological mirror symmetry naturally leads the following based on the above Theorem. 
\begin{conj}[Conjecture~\ref{conj : stability}]
Let $(\omega_1,\dots,\omega_n)$ be positive rational numbers such that $f_n(\ee[\omega_1]z_1,\dots,\ee[\omega_n]z_n)=f_n(z_1,\dots,z_n)$. 
There exists a Gepner type stability condition $\sigma$ on ${\rm HMF}^{L_{f_n}}_{S_n}(f_n)$ with respect to the auto-equivalence $(\vec{z}_1)$ and ${\bf e}[1/d_n]\in\CC$ 
in the sense of Toda~\cite[Definition 2.3]{To} such that $(\vec{z}_1).\sigma=\sigma.{\bf e}\left[\frac{1}{d_n}\right]$ and its stability function $Z_\sigma:K_0({\rm HMF}^{L_{f_n}}_{S_n}(f_n))\longrightarrow\CC$ is given by 
\small
\[
Z_\sigma(E_j):=\left\{
\begin{array}{ll}
\displaystyle \frac{1}{(2\pi\sqrt{-1})^{n}} \ee\left[ -\frac{j-1}{d_n} \right] \prod_{i=1}^{m} \left(1 - \ee \left[ -\omega_{2i-1} \right] \right) \cdot \int_{(\RR_{\ge 0})^n}e^{-\widetilde{f}_n({\bf x})}d{\bf x}, &\text{if}~ n=2m-1,\\
\displaystyle \frac{1}{(2\pi\sqrt{-1})^{n}} \ee\left[ \frac{j-1}{d_n} \right] \prod_{i=1}^{m} \left(1 - \ee \left[ -\omega_{2i} \right] \right) \cdot \int_{(\RR_{\ge 0})^n}e^{-\widetilde{f}_n({\bf x})}d{\bf x}, &\text{if}~ n=2m,
\end{array}
\right.  
\]
\normalsize
where $d{\bf x}$ denotes the standard volume form $dx_1\wedge \cdots \wedge dx_n$. 
\end{conj}
This conjecture holds for the case of $n=1$ and polynomials of ADE type \cite{T,KST} (see Proposition \ref{prop : stability condition} below).
\bigskip
We briefly outline the contents of the paper. 
In Section \ref{sec : preliminaries}, we recall some notations and terminologies of invertible polynomials, and then discuss the mirror isomorphism ${\bf mir}:\Omega_{\widetilde f_n}\cong\Omega_{f_n,G_{f_n}}$ for an invertible polynomial of chain type (Proposition \ref{prop : mirror isomorphism}). 
Section \ref{sec : main results} is the main part of this paper. After recalling some facts on ${\rm HMF}^{L_{f_n}}_{S_n}(f_n)$, we define the matrix ${\rm ch}^{(n)}_\Gamma$ 
and state Theorem \ref{thm : main theorem}. We also formulate a Gamma integral structure on $\Omega_{f_n,G_{f_n}}$ and state Theorem \ref{thm : Integral structure}.
Section \ref{sec : proof of thm} and Section \ref{sec : proof of integral structure} are devoted to proving Theorem \ref{thm : main theorem} and Theorem \ref{thm : Integral structure}, respectively.
In Section \ref{sec : Stokes matrix}, we study the Stokes matrix $S$ 
and the central connection matrix $C$ of the Frobenius manifold associated to $\widetilde f_n$ and a certain primitive form $\zeta$. 
We describe the relation of the monodromy data $(S,C)$ and the data $(\chi^{(n)}, {\rm ch}^{(n)}_\Gamma)$ of the category ${\rm HMF}^{L_{f_n}}_{S_n}(f_n)$. 
\bigskip
\noindent{\bf Acknowledgements.}
We would like to thank Akishi Ikeda, Hiroshi Iritani and Kohei Iwaki for valuable comments and discussions. 
The second named author is supported by JSPS KAKENHI Grant Number JP16H06337.
\section{Preliminaries}\label{sec : preliminaries}
\subsection{Invertible polynomials}
A polynomial $f=f(z_1,\dots, z_n)\in\CC[z_1,\dots, z_n]$ is called a {\em weighted homogeneous} polynomial  
if there are positive integers ${\rm w}_1,\dots ,{\rm w}_n$ and $d$ such that 
$f(\lambda^{{\rm w}_1} z_1, \dots, \lambda^{{\rm w}_n} z_n) = \lambda^d f(z_1,\dots ,z_n)$ for all $\lambda \in \CC^\ast$.
A weighted homogeneous polynomial $f$ is called {\em non-degenerate}
if it has at most an isolated critical point at the origin in $\CC^n$.
\begin{defn}\label{def:invertible}
A non-degenerate weighted homogeneous polynomial $f\in\CC[z_1,\dots, z_n]$ is called {\em invertible} if 
the following conditions are satisfied.
\begin{itemize}
\item 
The number of variables coincides with the number of monomials in $f$: 
\begin{equation}
f(z_1,\dots ,z_n)=\sum_{i=1}^n c_i\prod_{j=1}^n z_j^{{\mathbb E}_{ij}}
\end{equation}
for some coefficients $c_i\in\CC^\ast$ and non-negative integers 
${\mathbb E}_{ij}$ for $i,j=1,\dots, n$.
\item
The matrix $\EE:=(\EE_{ij})$ is invertible over $\QQ$.
\end{itemize}
\end{defn}
Let $f=\sum_{i=1}^n c_i\prod_{j=1}^n z_j^{\EE_{ij}}$ be an invertible polynomial. 
Define positive rational numbers $\omega_1,\dots,\omega_n$ called the {\em rational weights} of $f$ by the unique solution of the equation
\[
(\omega_1,\dots, \omega_n):= (1,\dots, 1) \EE^{-T}
\]
where $\EE^{-T}:=(\EE^{-1})^T$. Note that $f(\lambda^{\omega_1} z_1, \dots, \lambda^{\omega_n} z_n) = \lambda f(z_1,\dots ,z_n)$ for all $\lambda \in \CC^\ast$.
Without loss of generality one may assume that $c_i=1$ for all $i$ by rescaling the variables. 
An invertible polynomial $f$ can be written as a Thom--Sebastiani sum
$f=f_1\oplus\cdots\oplus f_p$ of invertible ones (in groups of different variables) $f_\nu$, $\nu=1,\dots, p$ of the following types (cf. \cite{KS}):
\begin{enumerate}
\item $z_1^{a_1}z_2+z_2^{a_2}z_3+\dots+z_{m-1}^{a_{m-1}}z_m+z_m^{a_m}$ ({\em chain type}, $m \geq 1$);
\item $z_1^{a_1}z_2+z_2^{a_2}z_3+\dots+z_{m-1}^{a_{m-1}}z_m+z_m^{a_m}z_1$ ({\em loop type}, $m \geq 2$).
\end{enumerate}
For an invertible polynomial $f=f(z_1,\dots, z_n)\in\CC[z_1,\dots, z_n]$, define a $\CC$-algebra called the {\em Jacobian algebra} $\Jac(f)$ as 
\begin{equation}
{\rm Jac}(f):=\CC[z_1,\dots, z_n]\left/\left(\frac{\p f}{\p z_1},\dots, \frac{\p f}{\p z_n}\right)\right..
\end{equation}
Since $f$ is non-degenerate, the Jacobian algebra ${\rm Jac}(f)$ is finite-dimensional. 
Set $\mu_f:=\dim_\CC\Jac(f)$ and call it the {\em Milnor number} of $f$. 
We sometimes include the case $n=0$. If $n=0$, set $\Jac(f):=\CC$, 
naturally considered as the $\CC$-algebra of constant functions on a point, and we have $\mu_f=1$.
Let $\Omega^p(\CC^n)$ be the complex vector space of regular $p$-forms on $\CC^n$. 
Consider the complex vector space 
\begin{equation}
\Omega_f := \Omega^n(\CC^n)/df \wedge \Omega^{n-1}(\CC^n).
\end{equation}
If $n=0$, then set $\Omega_f:=\CC$, the complex vector space of constant functions on a point.
Note that $\Omega_{f}$ is naturally a free $\Jac(f)$-module of rank one, namely, 
by choosing a nowhere vanishing $n$-form $d{\bf z}:=dz_1\wedge \dots \wedge dz_n$ we have the following isomorphism
\begin{equation}\label{eqn : isom from Jac to Omega}
\Jac(f)\stackrel{\cong}{\longrightarrow }\Omega_f,\quad 
[\phi({\bf z})]\mapsto [\phi({\bf z})d{\bf z}].
\end{equation}
\begin{prop}[cf. \cite{H}]
Define a $\CC$-bilinear form $J_{f}:\Omega_{f}\times \Omega_{f}\longrightarrow \CC$ by
\begin{equation}
J_{f}\left([\phi_1({\bf z}) d{\bf z}],[\phi_2({\bf z}) d{\bf z}]\right):=
{\rm Res}_{\CC^n}\left[
\begin{gathered}
\phi_1({\bf z}) \phi_2({\bf z}) dz_1\wedge \dots \wedge dz_n\\
\frac{\partial f}{\partial z_1}, \dots, \frac{\partial f}{\partial z_n}
\end{gathered}
\right].
\end{equation}
Then, the bilinear form $J_{f}$ on $\Omega_f$ is non-degenerate. 
\qed
\end{prop}
\begin{defn}
Let $f\in\CC[z_1,\dots, z_n]$ be an invertible polynomial.
The {\em group of maximal diagonal symmetries} $G_f$ of $f$ is defined as 
\begin{equation}
G_f := \left\{ (\lambda_1, \ldots , \lambda_n )\in(\CC^\ast)^n \, \left| \, f(\lambda_1z_1, \ldots, \lambda_n z_n) = 
f(z_1, \ldots, z_n) \right\} \right..
\end{equation}

Each element $g\in G_f$ has a unique expression of the form $g=({\bf e}[\alpha_1], \dots,{\bf e}[\alpha_n])$ with $0 \leq \alpha_i < 1$,
where ${\bf e}[\alpha] := \exp(2 \pi \sqrt{-1} \alpha)$.
The {\em age} of $g$ is defined as the rational number ${\rm age}(g) := \sum_{i=1}^n \alpha_i$. 
For each $g\in G_f$, denote by ${\rm Fix}(g):=\{z\in\CC^n~|~g\cdot z=z\}$ the fixed locus of $g$, by
$n_g:=\dim_\CC{\rm Fix}(g)$ its dimension and by $f^g:=f|_{{\rm Fix}(g)}$ the restriction of $f$
to the fixed locus of $g$. 
\end{defn}
Note that the function $f^g$ is an invertible polynomial. 
Indeed, for an invertible polynomial of chain type, which is of our main interest in this paper,
we have the following
\begin{prop}
Let $f_n = z_1^{a_1}z_2 + \dots + z_{n-1}^{a_{n-1}}z_n + z_n^{a_n}$ be an invertible polynomial of chain type.
For each $g \in G_{f_n}\backslash\{\id\}$, there exists $1\le k \le n$ such that
${\rm Fix}(g)=\{z_i=0\,|\,1\le i\le k\}$.
\qed
\end{prop}
We shall use the fact that, for each $g\in G_f$, 
$\Omega_{f^g}$ admits a natural $G_f$-action by restricting the $G_f$-action on $\CC^n$ to ${\rm Fix}(g)$ since $G_f$ acts diagonally on $\CC^n$.
Since $\Omega_f$ has a structure of $\QQ$-graded complex vector space with respect to the rational weights $\omega_1,\dots, \omega_n$,
we can associate the following $\QQ$-graded $\CC$-vector space:
\begin{defn}\label{defn : 2.5}
Define $\QQ$-graded complex vector space $\Omega_{f,G_f}$ by
\begin{equation}
\Omega_{f,G_f}:=\bigoplus_{g\in G_f}\Omega_{f,g},\quad \Omega_{f,g}:=(\Omega_{f^g})^{G_f}(-{\rm age}(g)),
\end{equation}
where $(\Omega_{f^g})^{G_f}$ denotes the $G_f$-invariant subspace of $\Omega_{f^g}$
and $(\Omega_{f^g})^{G_f}(-{\rm age}(g))$ denotes the $\QQ$-graded complex vector space $(\Omega_{f^g})^{G_f}$ shifted by $-{\rm age}(g)\in\QQ$.
\end{defn}
\begin{defn}
Define a non-degenerate $\CC$-bilinear form $J_{f,G_{f}}:\Omega_{f,G_{f}}\times \Omega_{f,G_{f}}\longrightarrow \CC$ by
\begin{subequations}
\begin{equation}
J_{f,G_{f}}:=\bigoplus_{g\in G_{f}}J_{f,g},
\end{equation}
where $J_{f,g}$ is a perfect $\CC$-bilinear form 
$J_{f,g}:\Omega_{f,g}\times \Omega_{f,{g^{-1}}}\longrightarrow \CC$ defined by 
\begin{equation}
J_{f,g}\left(\xi_1,\xi_2\right):=\dfrac{1}{|G_{f}|}\cdot {\rm Res}_{{\rm Fix}(g)}\left[
\begin{gathered}
\phi_1 \phi_2 dz_{n-n_g+1}\wedge\dots\wedge dz_{n}\\
\frac{\partial f^g}{\partial z_{n-n_g+1}},\dots,\frac{\partial f^g}{\partial z_{n}}
\end{gathered}
\right]
\end{equation}
for $\xi_1=[\phi_1 dz_{n-n_g+1}\wedge\dots\wedge dz_{n}]\in \Omega_{f,g}$ and 
$\xi_2=[\phi_2 dz_{n-n_{g^{-1}}+1}\wedge\dots\wedge dz_{n}]\in \Omega_{f,g^{-1}}$ (note that $n_g=n_{g^{-1}}$).
In particular for each $g\in G_{f}$ with $n_g=0$, we define 
\begin{equation}
J_{f,g}\left({\bf 1}_g,{\bf 1}_{g^{-1}}\right):=\dfrac{1}{|G_{f}|},
\end{equation}
\end{subequations}
where ${\bf 1}_g$ denotes the constant function $1$ on ${\rm Fix}(g)=\{0\}$.
\end{defn} 
We have the {\em mirror isomorphism} between $\Omega_{\widetilde f}$ and $\Omega_{f,G_f}$.
\begin{prop}[cf. \cite{K, EGT}]
Let $f=\sum_{i=1}^n \prod_{j=1}^n z_j^{{\mathbb E}_{ij}}$ be an invertible polynomial. 
There exists an isomorphism of $\QQ$-graded $\CC$-vector spaces
\[
{\bf mir}: \Omega_{\widetilde f}\cong \Omega_{f,G_f},
\]
where $\widetilde f:=\sum_{i=1}^n \prod_{j=1}^n x_j^{{\mathbb E}_{ji}}$ is the Berglund--H\"{u}bsch transpose of $f$.
\qed
\end{prop}
See Proposition~\ref{prop : mirror isomorphism} below for the details of the mirror isomorphism ${\bf mir}$ when $f$ is of chain type.
\subsection{Invertible polynomials of chain type}
From now on, we shall only consider invertible polynomials of chain type.
Set 
\begin{gather*}
f_n := f_n(z_1,\dots, z_n):=z_1^{a_1}z_2 + \dots + z_{n-1}^{a_{n-1}}z_n + z_n^{a_n},\\
\widetilde{f}_n:=\widetilde{f}_n(x_{1},\dots,x_{n}):=x_{1}^{a_{1}}+x_{1}x_{2}^{a_{2}}+\cdots+x_{n-1}x_{n}^{a_{n}}.
\end{gather*}
For simplicity, assume that $a_i\geq 2$ for all $i=1,\dots, n$. 
\begin{defn}\label{defn : monomial basis on Jac}
For each non-negative integer $n$, define sets $B'_{\widetilde{f}_n}, B_{\widetilde{f}_n}$ of monomials in $\CC[x_1,\dots, x_n]$ inductively as follows: 
Let $B'_{\widetilde{f}_0}:=\{1\}$ and 
\[
B'_{\widetilde{f}_n}:=\left\{x_{1}^{k_{1}}x_{2}^{k_{2}}\cdots x_{n}^{k_{n}}\,|\,0 \le k_{i} \le a_{i}-1\ (i=1,\dots, n-1),\ 0 \le k_{n} \le a_{n}-2\right\},\quad n\ge 1.
\]
For $n=0,1$, let $B_{\widetilde{f}_0}:=B_{\widetilde{f}_0}=\left\{1\right\}$ and $B_{\widetilde{f}_1}:=B'_{\widetilde{f}_1}=\left\{x_{1}^{k_{1}}\,|\,0 \le k_{1} \le a_{1}-2\right\}$.

For $n\ge 2$, let
\[
B_{\widetilde{f}_n}:=B'_{\widetilde{f}_n}\cup\left\{\phi^{(n-2)}(x_{1},\dots,x_{n-2})x_{n}^{a_{n}-1}\,|\, \phi^{(n-2)}(x_{1},\dots,x_{n-2})\in B_{\widetilde{f}_{n-2}}\right\}.
\] 
\end{defn}
\begin{prop}[\cite{K}]
The set $B_{\widetilde{f}_n}$ defines a $\CC$-basis of the Jacobian algebra $\Jac(\widetilde{f}_n)$. Namely, we have 
$\Jac(\widetilde{f}_n)=\langle [\phi^{(n)}({\bf x})]\,|\, \phi^{(n)}({\bf x})\in B_{\widetilde{f}_n}\rangle_\CC$.\qed
\end{prop}
Set $d_0:=1$ and $d_i:=a_1\cdots a_i$ for $i=1,\dots ,n$. Define a positive integer $\widetilde\mu_n$ by 
$\widetilde \mu_{n}:=\sum_{i=0}^n(-1)^{n-i}d_i$, which satisfies $\widetilde\mu_n=d_n-d_{n-1}+\widetilde\mu_{n-2}$.
\begin{cor}
The Milnor number $\mu_{\widetilde{f}_n}=\dim_\CC{\rm Jac}(\widetilde {f}_n)$ is given by $\widetilde \mu_{n}$.\qed
\end{cor}
Denote by $\EE$ the invertible matrix associated to $f_n$, which is given by 
\[
\EE=(a_i\delta_{i,j}+\delta_{i+1,j})=
\begin{pmatrix}
a_1 & 1 & 0 & \cdots & 0 \\
0 & a_2 & \ddots & \ddots & \vdots\\
\vdots & \ddots & \ddots & \ddots & 0\\
\vdots & \ddots & \ddots & \ddots & 1 \\
0 & \cdots & \cdots &0 & a_n
\end{pmatrix}. 
\]
\begin{defn}
For each ${\bf k}=(k_1,\dots, k_n)$, define rational numbers ${\bf \omega}^{(n)}_{{\bf k},1},\dots, {\bf \omega}^{(n)}_{{\bf k},n}$ by
\begin{equation}
(\omega^{(n)}_{{\bf k},1},\dots, {\bf \omega}^{(n)}_{{\bf k},n}):=(k_1+1,\dots, k_n+1)\EE^{-T}. 
\end{equation}
In particular, $\omega^{(n)}_{{\bf 0},1},\dots, \omega^{(n)}_{{\bf 0},n}$ are nothing but the rational weights of $f_n$ and 
will be denoted simply by $\omega^{(n)}_{1},\dots, \omega^{(n)}_{n}$.
It is easy to see that 
\begin{equation}
\omega^{(n)}_i=\sum_{l=i}^n(-1)^{l-i}\frac{d_{i-1}}{d_l},\quad i=1,\dots, n.
\end{equation}
\end{defn}
\begin{prop}
For each ${\bf k}=(k_1,\dots, k_n)$ such that
$0 \le k_{i} \le a_{i}-1\ (i=1,\dots, n-1),\ 0 \le k_{n} \le a_{n}-2$, we have $0< \omega^{(n)}_{{\bf k},i}<1$ for all $i=1,\dots, n$.
Moreover, for ${\bf k}^\ast:=(a_1-1,\dots, a_{n-1}-1,a_{n}-2)-{\bf k}$, we have the duality property
\begin{equation}
\omega^{(n)}_{{\bf k}^\ast,i}=1-\omega^{(n)}_{{\bf k},i},\quad i=1,\dots, n.
\end{equation}
\end{prop}
\begin{pf}
Since $\omega_{{\bf k},n}^{(n)}=(k_n+1)/a_n$, we have $0<\omega_{{\bf k},n}^{(n)}<1$. 
For $i=1,\dots,n-1$, we have $\omega_{{\bf k},i}^{(n)}=(k_i+1-\omega_{{\bf k},i+1}^{(n)})/a_{n-1}$, and hence we obtain $0<\omega_{{\bf k},i}^{(n)}<1$ inductively. 

By the definition of ${\bf k}^\ast$, we have 
\begin{equation*}
\left(\omega^{(n)}_{{\bf k},1},\dots, {\bf \omega}^{(n)}_{{\bf k},n}\right) + \left(\omega^{(n)}_{{\bf k}^\ast,1},\dots, {\bf \omega}^{(n)}_{{\bf k}^\ast,n}\right) 
= (a_1-1,\dots,a_{n-1}-1,a_n-2) \EE^{-T}
= (1,\dots,1) .
\end{equation*}
\qed
\end{pf}
On the structure of the group $G_{f_n}$, we have the following
\begin{prop}
The group $G_{f_n}$ is a cyclic group of order $d_n$ generated by the element 
\[
\left(\ee\left[\frac{1}{d_n}\right],\dots, \ee\left[(-1)^{i-1}\frac{d_{i-1}}{d_n}\right],\dots ,\ee\left[(-1)^{n-1}\frac{d_{n-1}}{d_n}\right]\right).
\]
\qed
\end{prop}
For each non-negative integer $n$, define sets $I'_n, I_n$ inductively as follows: 
Let $I'_0:=\{1\}$ and $I'_n:=\left\{\kappa\in\{1,\dots,d_n\}\,|\,a_n\nmid\kappa\right\}$, $n\ge1$.
For $n=0,1$, set $I_0:=I'_0=\{1\}$ and $I_1:=I'_1=\{1,2,\dots,a_1-1\}$. 
For $n\ge2$, set $I_n:=I'_n\sqcup I_{n-2}$. 
Note that $|I'_n|=d_n-d_{n-1}$ and $|I_n|=\widetilde\mu_n$. 

For $\kappa\in I'_n$ and $i=1,\dots, n$, set
\[
\omega^{(n)}_{\kappa,i}:=(-1)^{i-1}\frac{d_{i-1}}{d_n}\cdot \kappa-\left\lfloor(-1)^{i-1}\frac{d_{i-1}}{d_n}\cdot \kappa\right\rfloor.
\]
\begin{prop}\label{prop : duality of omega}
For $\kappa\in I'_n$ and $i=1,\dots, n$, we have the following duality property;
\[
\omega^{(n)}_{d_n-\kappa,i}=1-\omega^{(n)}_{\kappa,i} .
\]
\end{prop}
\begin{pf}
It follows from a direct calculation: 
\begin{align*}
\omega^{(n)}_{d_n-\kappa,i} & = (-1)^{i-1}\frac{d_{i-1}}{d_n}\cdot (d_n-\kappa)-\left\lfloor(-1)^{i-1}\frac{d_{i-1}}{d_n}\cdot (d_n-\kappa)\right\rfloor \\
& = 1-\left( (-1)^{i-1}\frac{d_{i-1}}{d_n}\cdot \kappa-\left\lfloor(-1)^{i-1}\frac{d_{i-1}}{d_n}\cdot \kappa\right\rfloor\right)  = 1-\omega^{(n)}_{\kappa,i} .
\end{align*}
\qed
\end{pf}
The following correspondence $\psi$ plays an important role in this paper. 
\begin{prop}\label{prop: 1-1}
The map 
\[
\psi:\{{\bf k}=(k_1,\dots, k_n)\,|\,0\leq k_i\leq a_i-1\ (i=1,\dots ,n-1),\ 0\leq k_n\leq a_n-2\}\\
\stackrel{\cong}{\longrightarrow}
I'_n
\]
defined by 
\[
\psi({\bf k})=\psi(k_1,\dots, k_n):=\sum_{l=1}^n(-1)^{l-1}\frac{d_n}{d_l}(k_l+1) .
\]
is a bijection of sets of $d_n-d_{n-1}$ elements and satisfies $\omega^{(n)}_{\psi({\bf k}),i}=\omega^{(n)}_{{\bf k},i}$ for $i=1,\dots, n$. 

In particular, we have $\psi({\bf 0})=\psi(0,\dots, 0)=\sum_{l=1}^n (-1)^{l-1}\frac{d_n}{d_l}$.
\end{prop}
\begin{pf}
One can check easily by a direct calculation. 
\qed
\end{pf}
\begin{cor}\label{cor : gamma integral}
For each ${\bf k}=(k_1,\dots, k_n)$ such that ${\bf x}^{\bf k}:=x_1^{k_1}\cdots x_n^{k_n}\in B'_{\widetilde{f}_n}$, we have
\[
\int_{(\RR_{\ge 0})^n}e^{-\widetilde{f}_n({\bf x})}{\bf x}^{\bf k}d{\bf x}=\frac{1}{d_n}\prod_{l=1}^n\Gamma(1-\omega_{d_n-\psi({\bf k}),l}^{(n)}),
\]
where ${\bf x}^{\bf k}d{\bf x}:=x_1^{k_1}\dots x_n^{k_n}dx_1\wedge \dots \wedge dx_n\in\Omega^n(\CC^n)$ and $\Gamma(s)$ denotes the Gamma function
\[
\Gamma(s):=\int_0^\infty y^{s-1}e^{-y}dy.
\]
\end{cor}
\begin{pf}
By using the change of variables $y_i:=\prod_{j=1}^nx_j^{\EE_{ij}}$, $i=1,\dots, n$ and the definition of 
${\bf \omega}^{(n)}_{{\bf k},1},\dots, {\bf \omega}^{(n)}_{{\bf k},n}$, we have 
\[
\int_{(\RR_{\ge 0})^n}e^{-\widetilde{f}_n({\bf x})}{\bf x}^{\bf k}d{\bf x}=\frac{1}{d_n}\prod_{l=1}^n\Gamma(\omega_{{\bf k},l}^{(n)}) .
\]
Hence, we obtain the statement by Proposition \ref{prop : duality of omega} and Proposition \ref{prop: 1-1}. 
\qed 
\end{pf}
\begin{rem}
For any invertible polynomial, one can define rational numbers $\omega_{{\bf k},i}$, $\omega_{\kappa, i}$ and a map $\psi$. 
Hence, we obtain the similar statement of Corollary \ref{cor : gamma integral} in the same way. 
\end{rem}
For each $\kappa\in I'_n$ define an element $g_\kappa \in G_{f_n}$ of order $d_n$ by
\[
g_\kappa:=\left(\ee\left[\frac{1}{d_n}\kappa\right],\dots, \ee\left[(-1)^{i-1}\frac{d_{i-1}}{d_n}\kappa\right],\dots ,\ee\left[(-1)^{n-1}\frac{d_{n-1}}{d_n}\kappa\right]\right) .
\]
We will define a basis $\{\zeta^{(n)}_\kappa\}_{\kappa\in I_n}$ (resp., $\{\xi^{(n)}_\kappa\}_{\kappa\in I_n}$) on $\Omega_{\widetilde{f}_n}$ (resp., $\Omega_{f_n,G_{f_n}}$) as follows: 
\begin{enumerate}
\item Let $\zeta^{(0)}_1:=\phi^{(0)}=1$ for $n=0$ and, for $n=1$,
\[
\zeta^{(1)}_\kappa := [\phi^{(1)}_\kappa(x_1)dx_1]=[x_1^{\kappa-1}dx_1],\quad \kappa\in I_1. 
\]
For $n\ge2$, let
\[
\zeta^{(n)}_\kappa:=
\begin{cases}
[{\bf x}^{\bf k}d{\bf x}]=[x_1^{k_1}\dots x_n^{k_n}dx_1\wedge \dots \wedge dx_n], & \kappa\in I'_n , \\
-[\phi^{(n-2)}_{\kappa}(x_1\,\dots,x_{n-2})x_n^{a_n-1}dx_1\wedge\cdots\wedge dx_n], & \kappa\in I_{n-2}, 
\end{cases}
\]
where ${\bf k}=(k_1,\dots,k_n)=\psi^{-1}(\kappa)$ for $\kappa\in I'_n$. 

\item Let $\xi^{(0)}_1:=1$ for $n=0$ and, for $n=1$,
\[
\xi^{(1)}_\kappa := {\bf 1}_{g_\kappa},\quad \kappa\in I_1, 
\]
where ${\bf 1}_{g_\kappa}$ denotes the generator of $\Omega_{f_1^{g_\kappa}}$, namely, the constant function $1$ on ${\rm Fix}(g_\kappa)=\{0\}$. 

For $n\ge2$, note that the natural inclusion map 
\[
G_{f_{n-2}}\hookrightarrow G_{f_n},\quad ({\bf e}[\alpha_1], \dots,{\bf e}[\alpha_{n-2}])\mapsto ({\bf e}[\alpha_1], \dots,{\bf e}[\alpha_{n-2}], 1,1)
\]
defines an subspace of $\Omega_{f_n,G_{f_n}}$ spanned by 
\[
[\xi\wedge d(z_{n-1}^{a_{n-1}})\wedge dz_n],\quad [\xi]\in (\Omega_{f_{n-2}^g})^{G_{f_{n-2}}},\ g\in G_{f_{n-2}}.
\]
Let
\[
\xi^{(n)}_\kappa:=
\begin{cases}
{\bf 1}_{g_\kappa} & \kappa\in I'_n, \\
[\overline{\xi}^{(n-2)}_{\kappa}\wedge d(z_{n-1}^{a_{n-1}})\wedge dz_n], & \kappa\in I_{n-2}, 
\end{cases}
\]
where ${\bf 1}_{g_\kappa}$ denotes the generator of $\Omega_{f_n^{g_\kappa}}$, namely, the constant function $1$ on ${\rm Fix}(g_\kappa)=\{0\}$ and $\overline{\xi}^{(n-2)}_{\kappa}$ does a differential form representing $\xi^{(n-2)}_{\kappa}$. 
\end{enumerate}
\begin{rem}\label{rem : index}
In this paper, it is important that matrices with respect to the basis $\{\xi^{(n)}_\kappa\}_{\kappa\in I_n}$ (resp., $\{\zeta^{(n)}_\kappa\}_{\kappa\in I_n}$) on $\Omega_{f_n,G_{f_n}}$ (resp., $\Omega_{\widetilde f_n}$) are labelled by the set $I_n$. 
\end{rem}
We have the isomorphism of $\QQ$-graded complex vector spaces
\begin{subequations}\label{eq: Omega_f induction}
\begin{equation}\label{eqn : decomposition in B}
\Omega_{\widetilde{f}_n}\cong \left(\bigoplus_{\kappa\in I'_n}\CC\cdot\zeta^{(n)}_{\kappa}\right)
\bigoplus\Omega_{\widetilde{f}_{n-2}}(-1),\ n\ge 2,
\end{equation}
where $\Omega_{\widetilde{f}_{n-2}}(-1)$ is identified with the subspace of $\Omega_{\widetilde{f}_n}$ spanned by $\zeta^{(n)}_\kappa,~ \kappa\in I_{n-2}$. 
On the other hand, we also have the isomorphism of $\QQ$-graded complex vector spaces
\begin{equation}\label{eqn : decomposition in A}
\Omega_{f_n,G_{f_n}}\cong \left(\bigoplus_{\kappa\in I'_n}\CC\cdot \xi^{(n)}_{\kappa}\right)\bigoplus\Omega_{f_{n-2},G_{f_{n-2}}}(-1),\ n\ge 2.
\end{equation}
\end{subequations}
where $\Omega_{f_{n-2},G_{f_{n-2}}}(-1)$ is identified with the subspace of $\Omega_{f_n,G_n}$ spanned by $\xi^{(n)}_{\kappa},~\kappa\in I_{n-2}$. 
Denote by $\eta^{(n)}$ the matrix representation of $J_{f_n,G_{f_n}}$ with respect to the basis $\{\xi^{(n)}_\kappa\}_{\kappa\in I_n}$ of $\Omega_{f_n,G_{f_n}}$. 
By the following proposition, it turns out that $\eta^{(n)}$ is also the matrix representation of $J_{\widetilde{f}}$ with respect to the basis $\{\zeta^{(n)}_\kappa\}_{\kappa\in I_n}$ of $\Omega_{\widetilde f_n}$. 
\begin{prop}\label{prop : mirror isomorphism}
The bijection $\psi$ in Proposition~\ref{prop: 1-1} induces the following bijection
\[
\Psi:\{x_1^{k_1}\cdots x_n^{k_n}dx_1\wedge\cdots\wedge dx_n\,|\,x_1^{k_1}\cdots x_n^{k_n}\in \widetilde{B}'_n\}
\stackrel{\cong}{\longrightarrow}
\{{\bf 1}_{g_\kappa}\,|\,\kappa\in I'_n\},
\]
which yields, via \eqref{eq: Omega_f induction}, the mirror isomorphism ${\bf mir}$.
In particular, the map $\Psi$ is chosen so that ${\bf mir}$ is an isomorphism of pairs 
of $\QQ$-graded complex vector spaces and non-degenerate bilinear forms:
\begin{equation}\label{eqn : mirror map}
{\bf mir}: (\Omega_{\widetilde f_n},J_{\widetilde{f}_n})\cong (\Omega_{f_n,G_{f_n}}, J_{f_n,G_{f_n}}),\quad
\zeta^{(n)}_\kappa\mapsto\xi^{(n)}_\kappa .
\end{equation}
Moreover, $\eta^{(n)}$ is given by 
\[
\eta^{(0)}=(1) , \quad \eta^{(1)}=\left( \dfrac{1}{a_{1}}\delta_{\kappa+\lambda,a_{1}} \right) , 
\]
and, for $n\ge2$,
\begin{equation}
\eta^{(n)}=
\begin{pmatrix}
\dfrac{1}{d_{n}}\delta_{\kappa+\lambda,d_{n}} & 0 \\
0 & -\dfrac{1}{a_{n}}\eta^{(n-2)}
\end{pmatrix} ,
\end{equation}
where $\kappa$ and $\lambda$ run the set $I'_n$. 
\end{prop}
\begin{pf}
First part is shown by \cite[Proposition 12]{EGT}. 
That is, the bijection $\Psi$ defined by
\[
\Psi(x_1^{k_1}\cdots x_n^{k_n}dx_1\wedge\cdots\wedge dx_n)={\bf 1}_{g_{\psi({\bf k})}},
\]
for ${\bf k}=(k_1,\dots,k_n)$ satisfying $x_1^{k_1}\cdots x_n^{k_n}\in B'_{\widetilde{f}_n}$ induces the isomorphism ${\bf mir}:\Omega_{\widetilde f_n}\cong\Omega_{f_n,G_{f_n}}$ by the decompositions \eqref{eq: Omega_f induction}. 

Here, we shall show by induction that the mirror isomorphism ${\bf mir}$ preserves non-degenerate bilinear forms $J_{\widetilde{f}_n}$ and $J_{f_n,G_{f_n}}$. 
It is clear for the case of $n=0$. If $n=1$, then for $k,k'=0,\dots,a_1-2$ we have 
\[
J_{\widetilde{f}_1}([x_1^kdx_1],[x_1^{k'}dx_1])=
{\rm Res}_{\CC}\left[
\begin{gathered}
x_{1}^{k+k'} dx_{1}\\
\frac{d \widetilde{f}_1}{d x_{1}}
\end{gathered}
\right]=
\begin{cases}
\frac{1}{a_1}, & \text{if }k+k'=a_1-2 \\
0, & \text{otherwise}\\
\end{cases}.
\]
On the other hand, we have a decomposition $G_{f_1}\cong\bigoplus_{\kappa=1}^{a_1-1}\CC\cdot{\bf 1}_{g_\kappa}$. 
Hence, they coincide via the mirror isomorphism. 

We shall show the general case $n\ge2$. 
For ${\bf k}=(k_1,\dots,k_n)$ and ${\bf k}'=(k'_1,\dots,k'_n)$ such that $x_1^{k_1}\cdots x_n^{k_n},~x_1^{k'_1}\cdots x_n^{k'_n}\in B'_{\widetilde{f}_n}$, we have 
\begin{align*}
J_{\widetilde{f}_n}([{\bf x}^{\bf k}d{\bf x}],[{\bf x}^{\bf k'}d{\bf x}])&={\rm Res}_{\CC^n}\left[
\begin{gathered}
x_1^{k_1+k'_1}\cdots x_n^{k_n+k'_n} dx_{1}\wedge\dots\wedge dx_{n}\\
\frac{\partial \widetilde{f}_n}{\partial x_{1}},\dots,\frac{\partial \widetilde{f}_n}{\partial x_{n}}
\end{gathered}
\right] \\
& =
\begin{cases}
\frac{1}{d_n} & ({\bf k}+{\bf k}'=(a_1-1,\dots,a_{n-1}-1,a_n-2)) \\
0 & (\text{otherwise}) \\
\end{cases},
\end{align*}
and since $a_{n-1}x_{n-2}x_{n-1}^{a_{n-1}-1}+x_{n}^{a_{n}}=0$ in ${\rm Jac}(\widetilde f_n)$ for $\phi,~\phi' \in B_{\widetilde{f}_{n-2}}$ we have 
\begin{align*}
 J_{\widetilde{f}_n}([\phi\, x_{n}^{a_{n}-1} d{\bf x}],[\phi'\, x_{n}^{a_{n}-1} d{\bf x}]) 
=& {\rm Res}_{\CC^n}\left[
\begin{gathered}
\phi\,\phi' x_{n}^{2a_{n}-2} dx_{1}\wedge\dots\wedge dx_{n}\\
\frac{\partial \widetilde{f}_n}{\partial x_{1}},\dots,\frac{\partial \widetilde{f}_n}{\partial x_{n}}
\end{gathered}
\right] \\
=& -\frac{1}{a_n}{\rm Res}_{\CC^{n-2}}\left[
\begin{gathered}
\phi\,\phi' dx_{1}\wedge\dots\wedge dx_{n-2}\\
\frac{\partial \widetilde{f}_{n-2}}{\partial x_{1}},\dots,\frac{\partial \widetilde{f}_{n-2}}{\partial x_{n-2}}
\end{gathered}
\right] \\
=& -\frac{1}{a_n}J_{\widetilde{f}_{n-2}}([\phi d{\bf x}],[\phi' d{\bf x}]).
\end{align*}
On the other side, since $|G_{f_n}|=d_n$ we have $J_{f_n,G_{f_n}}({\bf 1}_{g_\kappa},{\bf 1}_{g_\kappa^{-1}})=1/d_n$ for $\kappa\in I'_n$. 
For any $\xi=[\phi\,dz_{n-n_g+1}\wedge\cdots\wedge dz_{n-2}]\in(\Omega_{f^g_{n-2}})^{G_{f_{n-2}}}$, $\xi'=[\phi'dz_{n-n_{g^{-1}}+1}\wedge\cdots\wedge dz_{n-2}]\in(\Omega_{f^{g^{-1}}_{n-2}})^{G_{f_{n-2}}}$ we have 
\begin{align*}
& J_{f_n,G_{f_n}}([\overline{\xi}\wedge d(z_{n-1}^{a_{n-1}})\wedge dz_n],[\overline{\xi'}\wedge d(z_{n-1}^{a_{n-1}})\wedge dz_n]) \\
=& \dfrac{1}{d_n}\cdot {\rm Res}_{{\rm Fix}(g)}\left[
\begin{gathered}
\phi\,\phi'dz_{n-n_g+1}\wedge\dots\wedge dz_{n-2}\wedge d(z_{n-1}^{a_{n-1}})\wedge dz_n \\
\frac{\partial f_n^g}{\partial z_{n-n_g+1}},\dots,\frac{\partial f_n^g}{\partial z_{n}}
\end{gathered}
\right] \\
=& -\dfrac{1}{a_{n}}\cdot\dfrac{1}{d_{n-2}}\cdot {\rm Res}_{{\rm Fix}(g)}\left[
\begin{gathered}
\phi\,\phi'dz_{n-n_g+1}\wedge\dots\wedge dz_{n-2} \\
\frac{\partial f_{n-2}^g}{\partial z_{n-n_g+1}},\dots,\frac{\partial f_{n-2}^g}{\partial z_{n-2}}
\end{gathered}
\right] = -\frac{1}{a_n}J_{f_{n-2},G_{f_{n-2}}}(\xi,\xi'),
\end{align*}
since $z_{n-1}^{a_{n-1}}+a_nz_n^{a_n-1}=0$ in ${\rm Jac}(f_n)$. 

Therefore, we obtain the statement by the decompositions \eqref{eq: Omega_f induction}. 
\qed
\end{pf}
\section{Main results}\label{sec : main results}
First, we recall some notations and terminologies on the category of matrix factorizations. See \cite{AT} for details.
Consider the {\em maximal grading} $L_{f_n}$ of $f_n$, which is the abelian group defined by the quotient 
\begin{equation*}
L_{f_n}:=\left.\left(\bigoplus_{i=1}^n\ZZ\vec{z}_i\oplus \ZZ\vec{f}_n\right)\right/\left(\vec{f}_n-\sum_{j=1}^n\EE_{ij}\vec{z}_j;\ i=1,\dots ,n\right). 
\end{equation*}
Then $S_n:=\CC[z_1,\dots,z_n]$ has the natural $L_{f_n}$-graded $\CC$-algebra structure and 
one can define a triangulated category ${\rm HMF}^{L_{f_n}}_{S_n}(f_n)$ associated to $f_n$, the homotopy category of $L_{f_n}$-graded matrix factorizations.  
Denote by $K_{0}({\rm HMF}^{L_{f_n}}_{S_n}(f_n))$ the Grothendieck group of ${\rm HMF}^{L_{f_n}}_{S_n}(f_n)$. 
For each $\vec{l}\in L_{f_n}$, there is an auto-equivalence $(\vec{l})$ on ${\rm HMF}^{L_{f_n}}_{S_n}(f_n)$ defined by grading shift. 
In particular, we have $T^2=(\vec{f}_n)$, where $T$ is the translation functor on ${\rm HMF}^{L_{f_n}}_{S_n}(f_n)$. 
For each $j\in\ZZ$, set 
\begin{equation}\label{eqn : defn of exc obj}
E_j:=
\begin{cases}
\left(S_n/(z_1, z_3, \dots, z_{2m-1})\right)^{stab}(-(j-1)\vec{z}_1) & \text{if}\quad n=2m-1, \\
\left(S_n/(z_2, z_4, \dots, z_{2m})\right)^{stab}((j-1)\vec{z}_1) & \text{if}\quad n=2m ,
\end{cases}
\end{equation}
where $(-)^{stab}$ is the stabilization defined as follows: 
\begin{itemize}
\item For $n=2m-1$, the Koszul resolution of $S_n/(z_1, z_3, \dots, z_{2m-1})$ yields an $L_{f_n}$-graded matrix factorization $\overline{F}=(\mfs{F_0}{f_0}{f_1}{F_1})$ of $f_n$ such that
\begin{equation}
F_0:=\bigoplus_{k}\left(\bigwedge^{2k} \bigoplus_{i=1}^{m}S_n(-\vec{z}_{2i-1})\right)(k\vec{f}_n),\quad F_1:=\bigoplus_{k}\left(\bigwedge^{2k-1} \bigoplus_{i=1}^{m}S_n(-\vec{z}_{2i-1})\right)(k\vec{f}_n), 
\end{equation}
and denote by $\left(S_n/(z_1, z_3, \dots, z_{2m-1})\right)^{stab}$ the $L_{f_n}$-graded matrix factorization $\overline{F}$. 
\item For $n=2m$, the Koszul resolution of $S_n/(z_2, z_4, \dots, z_{2m})$ yields an $L_{f_n}$-graded matrix factorization $\overline{F}=(\mfs{F_0}{f_0}{f_1}{F_1})$ of $f_n$ such that
\begin{equation}
F_0:=\bigoplus_{k}\left(\bigwedge^{2k} \bigoplus_{i=1}^{m}S_n(-\vec{z}_{2i})\right)(k\vec{f}_n),\quad F_1:=\bigoplus_{k}\left(\bigwedge^{2k-1} \bigoplus_{i=1}^{m}S_n(-\vec{z}_{2i})\right)(k\vec{f}_n), 
\end{equation}
and denote by $\left(S_n/(z_2, z_4, \dots, z_{2m})\right)^{stab}$ the $L_{f_n}$-graded matrix factorization $\overline{F}$. 
\end{itemize}
Motivated by the homological mirror symmetry, it is proven that ${\rm HMF}^{L_{f_n}}_{S_n}(f_n)$ admits a full exceptional collection which is a Lefschetz decomposition in the sense of Kuznetsov--Smirnov \cite{KuS}.
\begin{prop}[{\cite[Theorem 1.3]{AT}}]\label{prop : fec and Euler}
The sequence $(E_1,\dots,E_{\widetilde{\mu}_n})$ is a full exceptional collection such that $\chi(E_i,E_j)=\chi^{(n)}_{ij}$ where 
\begin{equation}
\chi(E_i,E_j):=\sum_{p\in\ZZ}(-1)^p\dim_\CC{\rm HMF}^{L_{f_n}}_{S_n}(f_n)(E_i,T^p E_j) , 
\end{equation}
and $\chi^{(n)}$ is a matrix of size $\widetilde\mu_n$ defined by $\chi^{(n)}=1/\varphi_n(N)$, $N:=(\delta_{i+1,j})$, and 
\[
\varphi_n(t):=\prod_{i=0}^{n}\left(1-t^{d_{i}}\right)^{(-1)^{n-i}}\ (n\ge 1),\quad \varphi_0(t):=1-t.
\]
\qed
\end{prop}
It is known that there exists the Serre functor $\S^{(n)}$ on ${\rm HMF}^{L_{f_n}}_{S_n}(f_n)$, which is given by $T^{n}(-\vec{z}_1-\dots -\vec{z}_n)$. Denote by ${\bf S}^{(n)}$ the matrix representation of the automorphism induced by the Serre functor $\S^{(n)}$ with respect to $\{[E_j]\}_{j=1}^{\widetilde\mu_n}$, which satisfies
\begin{equation}
{\bf S}^{(n)}=(\chi^{(n)})^{-1}(\chi^{(n)})^T .
\end{equation}
Define a diagonal matrix $\widetilde Q^{(n)}$ of size $\widetilde{\mu}_{n}$ inductively as follows: 
\[
\widetilde Q^{(0)}:=(0) , \quad \widetilde Q^{(1)}:=\left( \left( \omega^{(1)}_{\kappa,1} -\frac{1}{2} \right) \delta_{\kappa\lambda} \right), \quad 
\]
and, for $n\ge 2$, 
\[
\widetilde Q^{(n)}:=
\begin{pmatrix}
\widetilde P^{(n)} & 0 \\
0 & \widetilde Q^{(n-2)}
\end{pmatrix},
\]
where $\widetilde  P^{(n)}=(\widetilde P^{(n)}_{\kappa\lambda})$ is a matrix of size $(d_n-d_{n-1})$ given by 
\[
\widetilde P^{(n)}_{\kappa\lambda}:= \left(\sum_{l=1}^{n} \left( \omega^{(n)}_{\kappa,l} -\frac{1}{2} \right) \right) \delta_{\kappa\lambda},\quad 
\kappa,\lambda\in I'_n.
\]
Now we are ready to introduce objects of our main interest.
\begin{defn}\label{defn : ch Gamma}
For $\kappa\in I'_n$, define $c^{(n)}_{\kappa}\in \CC$ by 
\begin{align}
c^{(n)}_{\kappa}:
& =\left\{
\begin{array}{ll}
\displaystyle \prod_{l=1}^{n}\Gamma \left( 1 - \omega^{(n)}_{\kappa,l} \right) \cdot \prod_{i=1}^{m} \left(1 - \ee \left[ \omega^{(2m-1)}_{\kappa,2i-1} \right] \right), &\text{if}\quad  n=2m-1,\\
\displaystyle \prod_{l=1}^{n}\Gamma \left( 1 - \omega^{(n)}_{\kappa,l} \right) \cdot \prod_{i=1}^{m} \left(1 - \ee \left[ \omega^{(2m)}_{\kappa,2i} \right] \right), &\text{if}\quad n=2m.
\end{array}
\right. 
\end{align}

For $j\in\ZZ$, define $\ch^{(n)}_{\Gamma,j}\in\CC^{\widetilde{\mu}_{n}}=\CC^{(d_n-d_{n-1})+\widetilde\mu_{n-2}}$ inductively by 
\begin{subequations}
\begin{gather}
\ch^{(0)}_{\Gamma,j} := (1) , \quad
\left(\ch^{(1)}_{\Gamma,j}\right)_{\kappa} := c^{(1)}_{\kappa}\ee\left[\omega_{\kappa,1}^{(1)}(j-1)\right],\ \kappa\in I_1=I'_1,\\
\left(\ch^{(n)}_{\Gamma,j}\right)_{\kappa} := \begin{cases}
c^{(n)}_\kappa\ee\left[(-1)^{n-1}\omega_{\kappa,1}^{(n)}(j-1)\right], & \ \text{if}\ \kappa\in I'_n,\\
2\pi\sqrt{-1}~\left(\ch^{(n-2)}_{\Gamma,j}\right)_\kappa, & \ \text{if}\ \kappa\in I_{n-2}=I_n\setminus I'_n.
\end{cases}
\end{gather}
\end{subequations}
Define a matrix $\ch^{(n)}_{\Gamma}$ of size $\widetilde \mu_n$ by 
\begin{equation}
\ch^{(n)}_{\Gamma}:=\left( \ch^{(n)}_{\Gamma,1},\cdots,\ch^{(n)}_{\Gamma,\widetilde{\mu}_{n}} \right) .
\end{equation}
\end{defn}
Note that the matrix ${\rm ch}^{(n)}_\Gamma$ is invertible.
\begin{rem}
The first part $\prod_{l=1}^n \Gamma(1-\omega_{\kappa,l}^{(n)})$ of $c_\kappa^{(n)}$ can be considered as the $\Gamma$-class on the $\kappa$-sector $\Omega_{f_n,{g_\kappa}}$ (see \cite[Definition 2.17]{CIR}) and the last part of $c_\kappa^{(n)}$ can be considered as the Chern character of $E_1$ (see \cite{PV,PV2} for Chern characters for $L_{f_n}$-graded matrix factorizations). 
The part $\ee\left[(-1)^{n-1}\omega_{\kappa,1}^{(n)}(j-1)\right]$ comes from the automorphism on the Grothendieck group $K_0({\rm HMF}^{L_{f_n}}_{S_n}(f_n))$ induced by the auto-equivalence $((-1)^nj\vec{z}_1)$ whose matrix representation is given by 
\[
\begin{pmatrix}
\left(\ee\left[(-1)^{n-1}\omega_{\kappa,1}^{(n)}(j-1)\right]\delta_{\kappa\lambda}\right) & 0 \\
0 & (1)
\end{pmatrix},
\]
which acts on the vector $\ch^{(n)}_{\Gamma,1}$ to get $\ch^{(n)}_{\Gamma,j}$.
Therefore, $\ch^{(n)}_{\Gamma,j}$ can be considered as the matrix representation of ``$\widehat{\Gamma}_{f_n,G_{f_n}} {\rm Ch}(E_j)$'', the Chern character of $E_j$ multiplied by the Gamma class of the pair $(f_n,G_{f_n})$, with respect to the basis $\{\xi^{(n)}_{\kappa}\}_{\kappa\in I_n}$. See below for the usual Gamma class on a smooth projective variety.
\end{rem}
The following theorem is the analogue in our setting of results of \cite[Proposition 2.10]{I}, namely, 
the necessary conditions for the matrix $(2\pi)^{-\frac{n}{2}}{\rm ch}^{(n)}_\Gamma$ to give a central connection matrix of a Frobenius manifold whose non-degenerate bilinear form on the tangent space, the grading matrix and the Stokes matrix are given by $\eta^{(n)}$, $\widetilde Q^{(n)}$ and $\chi^{(n)}$, respectively.
\begin{thm}\label{thm : main theorem}
We have the following equalities:
\begin{subequations}
\begin{equation}
\label{eqn : Main theorem 1}
\left(\dfrac{1}{(2\pi)^{\frac{n}{2}}}\ch_{\Gamma}^{(n)}\right)^{-1}\ee\left[\widetilde Q^{(n)}\right]\left(\dfrac{1}{(2\pi)^{\frac{n}{2}}}\ch_{\Gamma}^{(n)}\right)={\bf S}^{(n)},
\end{equation}
\begin{equation}
\label{eqn : Main theorem 2}
\left(\dfrac{1}{(2\pi)^{\frac{n}{2}}}\ch_{\Gamma}^{(n)}\right)^T\ee\left[\dfrac{1}{2}\widetilde Q^{(n)}\right]\eta^{(n)}\left(\dfrac{1}{(2\pi)^{\frac{n}{2}}}\ch_{\Gamma}^{(n)}\right)=\chi^{(n)}.
\end{equation}
\end{subequations}
\end{thm}
We prove Theorem \ref{thm : main theorem} in Section \ref{sec : proof of thm}. 
Here, we recall (refined) Dubrovin's conjecture and $\Gamma$-conjecture II. 
Given an algebraic variety $X$ and $\widetilde{\mu}:=\dim_\CC \bigoplus_{p\in\ZZ}H^p(X;\CC)$, the Gamma class $\widehat{\Gamma}_X$ of $X$ is defined to be 
\begin{equation}
\widehat{\Gamma}_X:=\prod_{i=1}^{\widetilde{\mu}}\Gamma(1+\delta_i),
\end{equation}
where $\delta_1,\dots,\delta_{\widetilde{\mu}}$ are the Chern roots of the tangent bundle $TX$. 

We refer the reader to \cite{D2,D3,CDG} for the monodromy data of quantum cohomologies (see also Section \ref{sec : Stokes matrix}). 
\begin{conj}[\cite{D2,CDG,GGI}]\label{conj : refined Dubrovin}
Let $X$ be an $n$-dimensional smooth Fano variety satisfying $h^{p,q}(X)=0$ for any $p \ne q$. 
\begin{enumerate}
\item The quantum cohomology of $X$ is semi-simple if and only if there exists a full exceptional collection in the derived category $\D^b(X)$.
\item If the quantum cohomology of $X$ is semi-simple, then for any oriented line $\ell$ (with $\phi\in[0,1)$) there is a correspondence between the monodromy data $(S^\phi,C^\phi)$ and full exceptional collections $(\E^\phi_1,\dots,\E^\phi_{\widetilde\mu})$ in $\D^b(X)$. 
\item The monodromy data $(S^\phi,C^\phi)$ are related to the following geometric data of the corresponding full exceptional collection $(\E^\phi_1,\dots,\E^\phi_{\widetilde\mu})$: 
\begin{enumerate}
\renewcommand{\labelenumii}{(\alph{enumii})}
\item the Stokes matrix is equal to the Euler form on $K_0(X)_\CC:=K_0(\D^b(X))\otimes_\ZZ \CC$, computed with respect to the exceptional basis $\{[\E^\phi_j]\}_{j=1}^{\widetilde{\mu}}$
\begin{equation}
S^\phi_{ij}=\chi(\E^\phi_i,\E^\phi_j) .
\end{equation}
\item ($\Gamma$-conjecture II) the $j$-th column of the central connection matrix $C^\phi$ on a homogeneous basis $\{T_j\}_{j=1}^{\widetilde\mu}$ of the $H^*(X;\CC)$ is given by 
\begin{equation}
C^\phi_j=\frac{1}{(2\pi)^{\frac{n}{2}}}\widehat{\Gamma}_X {\rm Ch}(\E^\phi_j) , 
\end{equation}
where ${\rm Ch}(\E)=\sum_{k=1}^{{\rm rank}\,\E}{\bf e}[\delta_k(\E)]$ for the Chern roots $\delta_1(\E),\dots,\delta_{{\rm rank}\,\E}(\E)$ of $\E$. 
\end{enumerate} 
\end{enumerate}
\end{conj}
Let us return to our setting. 
\begin{defn}\label{defn : Gamma integral structure}
We define the {\em $K$-group framing} ${\rm ch}_\Gamma^{(n)} : K_0({\rm HMF}^{L_{f_n}}_{S_n}(f_n))\longrightarrow\Omega_{f_n,G_{f_n}}$ of the $\CC$-vector space $\Omega_{f_n,G_{f_n}}$ by 
\begin{equation}\label{eqn : Gamma integral str}
{\rm ch}_\Gamma^{(n)} : K_0({\rm HMF}^{L_{f_n}}_{S_n}(f_n))\longrightarrow\Omega_{f_n,G_{f_n}},\quad [E_j]\mapsto
\sum_{\kappa\in I_n}{\rm ch}_{\Gamma,\kappa j}^{(n)}\,\xi^{(n)}_\kappa ,
\end{equation}
where $E_j$ is the exceptional object given by \eqref{eqn : defn of exc obj}. 
We call the image 
\begin{equation}
\Omega_{f_n,G_{f_n};\ZZ}:=\frac{1}{(2\pi\sqrt{-1})^n}{\rm ch}_\Gamma^{(n)}\left(K_0({\rm HMF}^{L_{f_n}}_{S_n}(f_n))\right) .
\end{equation}
of the $K$-group framing the {\em Gamma integral structure} of $\Omega_{f_n,G_{f_n}}$. 
\end{defn}
See \cite[Definition 2.9]{I} and \cite[Definition 2.19]{CIR} for the related notion, the Gamma integral structure of the quantum D-modules.
\begin{rem}
The constant factor $(2\pi\sqrt{-1})^{-n}$ is due to the fact that $\Omega_{f_n, G_{f_n}}$ (and $\Omega_{\widetilde f_n}$) has the natural weight $n$ from the view point of the Hodge theory.
\end{rem}
The Gamma integral structure $\Omega_{f_n,G_{f_n};\ZZ}$ admits a $\ZZ/d_n\ZZ$-action. 
This action is induced by the grading shift functor $(\vec{z}_1)$ of ${\rm HMF}^{L_{f_n}}_{S_n}(f_n)$. 
Under the isomorphism $H_n(\CC^n,{\rm Re}(\widetilde f_n)\gg0;\ZZ)\cong H_n(\CC^n,\widetilde f_n^{-1}(1);\ZZ) \cong H_{n-1}(\widetilde f_n^{-1}(1);\ZZ)$, 
define a ``Poincar\'{e} Duality'' map ${\mathbb D}: H_{n-1}(\widetilde{f}_n^{-1}(1);\ZZ)\to\Omega_{\widetilde f_n}$ by 
\[
{\mathbb D}(L):=\frac{1}{(2\pi\sqrt{-1})^{n}}\sum_{\kappa\in I_n}\left(\sum_{\lambda\in I_n}\eta^{\lambda\kappa}\int_\Gamma e^{-\widetilde f_n} \zeta^{(n)}_\lambda\right)\zeta^{(n)}_\kappa,\quad L\in H_{n-1}(\widetilde{f}_n^{-1}(1);\ZZ)
\]
where $\Gamma$ is the element of $H_n(\CC^n,{\rm Re}(\widetilde f_n)\gg0;\ZZ)$ corresponding to $L\in H_{n-1}(\widetilde{f}_n^{-1}(1);\ZZ)$,
$\zeta^{(n)}_\lambda$ in the integrand denotes the class in $\HH^n(\Omega^\bullet_{\CC^n}, d-d\widetilde f_n\wedge)$ 
corresponding to the class $\zeta^{(n)}_\lambda\in \Omega_{\widetilde f_n}$ under the isomorphism
$\HH^n(\Omega^\bullet_{\CC^n}, d-d\widetilde f_n\wedge)\cong\Omega_{\widetilde f_n}$, 
and $(\eta^{\lambda\kappa})$ is the inverse matrix of $\eta^{(n)}$. 
We shall denote by $\Omega_{\widetilde{f}_n;\ZZ}$ the natural integral structure induced by the Milnor homology, that is, 
\begin{equation}
\Omega_{\widetilde{f}_n;\ZZ}:={\mathbb D}\left(H_{n-1}(\widetilde{f}_n^{-1}(1);\ZZ)\right) .
\end{equation} 
The $\ZZ/d_n\ZZ$-action on the polynomial ring $\CC[x_1,\dots,x_n]$ given by 
\[
(x_1,\dots, x_k,\dots, x_n)\mapsto \left(\ee\left[\frac{1}{d_1}\right]x_1,\dots, \ee\left[\frac{(-1)^{k-1}}{d_k}\right]x_k,\dots ,\ee\left[\frac{(-1)^{n-1}}{d_n}\right]x_n\right)
\]
induces the one on the integral structure $\Omega_{\widetilde{f}_n;\ZZ}$. 
We obtain the following theorem as an analogue of the $\Gamma$-conjecture, which is proven for weak Fano toric orbifolds (\cite[Theorem 1.1]{I} and \cite[Theorem 31]{CCIT}). 
\begin{thm}[{See Milanov--Zha~\cite[Theorem 1.3]{MZ} for ADE cases}]\label{thm : Integral structure}
The mirror isomorphism ${\bf mir}: \Omega_{\widetilde f_n}\cong \Omega_{f_n,G_{f_n}}$ induces an isomorphism of free $\ZZ$-modules $\Omega_{\widetilde f_n;\ZZ}\cong \Omega_{f_n,G_{f_n};\ZZ}$ and the isomorphism is $\ZZ/d_n\ZZ$-equivariant.
\end{thm}
We prove the theorem in Section \ref{sec : proof of integral structure}.
Theorem \ref{thm : Integral structure} can naturally be generalized to the Thom--Sebastiani sum of invertible polynomials of chain type. 
We omit and leave it to the readers.
It is expected that the oscillatory integral $\int e^{-\widetilde{f}_n({\bf x})} d{\bf x}$ defines a stability condition as a central charge on $\D^b{\rm Fuk}^\to (\widetilde{f}_n)$ (see \cite{T,B2}). 
From the view point of the homological mirror symmetry conjecture for invertible polynomials, the mirror object dual to $\int_{\Gamma_j} e^{-\widetilde{f}_n({\bf x})} d{\bf x}$ is given by $\sum_{\lambda\in I_n}\eta_{\psi({\bf 0})\lambda}{\rm ch}^{(n)}_{\Gamma,\lambda j}$.  
Based on Corollary \ref{cor : gamma integral} and Theorem \ref{thm : Integral structure}, we expect the following 
\begin{conj}\label{conj : stability}
There exists a Gepner type stability condition $\sigma$ on ${\rm HMF}^{L_{f_n}}_{S_n}(f_n)$ with respect to the auto-equivalence $(\vec{z}_1)$ and ${\bf e}[1/d_n]\in\CC$ 
in the sense of Toda~\cite[Definition 2.3]{To} such that $(\vec{z}_1).\sigma=\sigma.{\bf e}\left[\frac{1}{d_n}\right]$ and its stability function $Z_\sigma:K_0({\rm HMF}_{S_n}^{L_{f_n}}(f_n))\longrightarrow\CC$ is given by
\small
\begin{equation*}
Z_\sigma([E_j]):=\left\{
\begin{array}{ll}
\displaystyle \frac{1}{(2\pi\sqrt{-1})^{n}} \ee\left[-\frac{j-1}{d_n}\right]\prod_{i=1}^{m} \left(1 - \ee \left[ -\omega^{(n)}_{2i-1} \right] \right) \cdot \int_{(\RR_{\ge 0})^n}e^{-\widetilde{f}_n({\bf x})}d{\bf x}, &\text{if}~ n=2m-1,\\
\displaystyle \frac{1}{(2\pi\sqrt{-1})^{n}} \ee\left[\frac{j-1}{d_n}\right]\prod_{i=1}^{m} \left(1 - \ee \left[ -\omega^{(n)}_{2i} \right] \right) \cdot \int_{(\RR_{\ge 0})^n}e^{-\widetilde{f}_n({\bf x})}d{\bf x}, &\text{if}~ n=2m,
\end{array}
\right.  
\end{equation*}
\normalsize
where $E_j$ is the exceptional object given by \eqref{eqn : defn of exc obj}.
\end{conj}
Note that this conjecture can also be generalized to the Thom--Sebastiani sum of invertible polynomials of chain type. 
A stability condition naturally associated to ${\rm HMF}_{S_n}^{L_{f_n}}(f_n)$ is constructed 
in \cite{T} for $n=1$ and in \cite{KST} for the cases when $f_3$  is of ADE type and $L_{f_3}\cong \ZZ$. 
For this stability condition, indecomposable objects of ${\rm HMF}_{S_n}^{L_{f_n}}(f_n)$ are semi-stable and each irreducible morphism between two indecomposable objects has phase $1/h$, where $h$ is the Coxeter number of the root system of the type of $\widetilde f_3$. 
Since the space of stability conditions has the $\CC$-action, we have the following result. 
\begin{prop}\label{prop : stability condition}
Conjecture \ref{conj : stability} holds for $n=1$ and for invertible polynomials of ADE type in two and three variables which is 
the Thom--Sebastiani sum of invertible polynomials of chain type. 
\end{prop}
\begin{pf}
Note that the Berglund--H\"{u}bsch transpose respects the set of  invertible polynomials of ADE type. 
One can easily show by the same argument as in \cite{KST2, AT} that for any invertible polynomial $f$ of ADE type in two or three variables, 
there exists an equivalence of triangulated categories ${\rm HMF}_{S}^{L_{f}}(f)\cong\D^b(\CC\vec{\Delta})$
where $\vec{\Delta}$ is the Dynkin quiver of the type of $\widetilde f$. 
In particular, there are finitely many indecomposable objects up to isomorphism and translation, any indecomposable object in ${\rm HMF}_{S}^{L_{f}}(f)$ can be constructed from one indecomposable object (for example, the one corresponding to the edge vertex of the longest arm of $\vec{\Delta}$) by using Auslander--Reiten triangles and any morphism is a composition of irreducible ones.
By requiring that all indecomposable objects are semi-stable and the auto-equivalences $(\vec{z}_i)$ respect semi-stable objects which increase those phases by $2\omega_{i}$, 
one obtains the expected stability condition (up to $\CC$-action) similarly to the case when $L_{f}\cong \ZZ$ in \cite{KST} since 
the Serre functor increases the phase $\sum_{i=1}^n(1-2\omega_i)=1-2/h$ where $h$ is the Coxeter number of the root system associated to $\vec{\Delta}$. 
\qed
\end{pf}
\section{Proof of Theorem \ref{thm : main theorem}}\label{sec : proof of thm}
Since we have 
\begin{equation*}
\widetilde Q^{(n)}\eta^{(n)}=-\eta^{(n)}\widetilde Q^{(n)},
\end{equation*}
by Proposition~\ref{prop : duality of omega}, it is easy to see that \eqref{eqn : Main theorem 1} follows from \eqref{eqn : Main theorem 2}. 
Therefore, only the proof for \eqref{eqn : Main theorem 2} is necessary, which will be done by induction. 
We prepare some lemmas. Let 
\[
p_n(t):=\frac{1}{\varphi_n(t)}\cdot (1-t^{d_n})=\prod_{i=1}^{n}\left(1-t^{d_{i-1}}\right)^{(-1)^{n-i}}\ (n\ge 1),\quad p_0(t):=1.
\]
Note that $p_n(t)$ is a polynomial in $t$ since $p_1(t)=1-t$ and 
\[
p_{n}(t)= p_{n-2}(t)\cdot\dfrac{1-t^{d_{n-1}}}{1-t^{d_{n-2}}},\quad d_{n-1}=d_{n-2}\cdot a_{n-1}.
\]
\begin{lem}\label{lem : calculation of coefficients of LHS}
We have 
\[
\frac{1}{(2\pi)^{n}}c^{(n)}_{\kappa}\ee\left[\frac{1}{2}\sum_{l=1}^{n}\left(\omega^{(n)}_{\kappa,l}-\frac{1}{2}\right)\right]\frac{1}{d_{n}}c^{(n)}_{d_{n}-\kappa}
=\frac{1}{d_{n}} p_{n}\left( \ee \left[ (-1)^{n-1}\omega_{\kappa,1}^{(n)} \right] \right).
\]
\end{lem}
\begin{pf}
By the Euler's reflection formula and the definition of $\omega^{(n)}_{\kappa,l}$, we have
\[
\Gamma \left( 1 - \omega^{(n)}_{\kappa,l} \right)\Gamma \left( 1 - \omega^{(n)}_{d_n-\kappa,l} \right) \ee\left[\frac{1}{2}\left(\omega^{(n)}_{\kappa,l}-\frac{1}{2}\right)\right]
= \frac{2\pi}{1-\ee\left[-\omega^{(n)}_{\kappa,l}\right]}.
\]
Therefore, for $n=2m-1$ we have
\begin{eqnarray*}
& &\frac{1}{(2\pi)^{2m-1}}c^{(2m-1)}_{\kappa}\ee\left[\frac{1}{2}\sum_{l=1}^{2m-1}\left(\omega^{(2m-1)}_{\kappa,l}-\frac{1}{2}\right)\right]\frac{1}{d_{2m-1}}c^{(2m-1)}_{d_{2m-1}-\kappa}\\
&=&\frac{1}{d_{2m-1}}\prod_{i=1}^{m} \left(1 - \ee \left[ \omega^{(2m-1)}_{\kappa,2i-1} \right] \right) \left(1 - \ee \left[ \omega^{(2m-1)}_{d_{2m-1}-\kappa,2i-1} \right] \right) \\
& &\quad \cdot \frac{1}{(2\pi)^{2m-1}}
\prod_{l=1}^{2m-1} \Gamma \left( 1 - \omega^{(2m-1)}_{\kappa,l} \right)\Gamma \left( 1 - \omega^{(2m-1)}_{d_{2m-1}-\kappa,l} \right) \ee\left[\frac{1}{2}\left(\omega^{(2m-1)}_{\kappa,l}-\frac{1}{2}\right)\right] \\
&=&\frac{1}{d_{2m-1}} \prod_{i=1}^{m} \left(1 - \ee \left[ \omega^{(2m-1)}_{\kappa,2i-1} \right] \right)\cdot \prod_{i=1}^{m}\left(1 - \ee \left[ -\omega^{(2m-1)}_{\kappa,2i} \right] \right)^{-1}\\
&=&\frac{1}{d_{2m-1}} p_{2m-1}\left( \ee \left[ \omega^{(2m-1)}_{\kappa,1} \right] \right) ,
\end{eqnarray*}
and for $n=2m$ we have
\begin{eqnarray*}
& &\frac{1}{(2\pi)^{2m}}c^{(2m)}_{\kappa}\ee\left[\frac{1}{2}\sum_{l=1}^{2m}\left(\omega^{(2m)}_{\kappa,l}-\frac{1}{2}\right)\right]\frac{1}{d_{2m}}c^{(2m)}_{d_{2m}-\kappa}\\
&=&\frac{1}{d_{2m}}\prod_{i=1}^{m} \left(1 - \ee \left[ \omega^{(2m)}_{\kappa,2i} \right] \right) \left(1 - \ee \left[ \omega^{(2m)}_{d_{2m}-\kappa,2i} \right] \right)  \\
& &\quad \cdot \frac{1}{(2\pi)^{2m}}
\prod_{l=1}^{2m} \Gamma \left( 1 - \omega^{(2m)}_{\kappa,l} \right)\Gamma \left( 1 - \omega^{(2m)}_{d_{2m}-\kappa,l} \right) \ee\left[\frac{1}{2}\left(\omega^{(2m)}_{\kappa,l}-\frac{1}{2}\right)\right] \\
&=&\frac{1}{d_{2m}} \prod_{i=1}^{m} \left(1 - \ee \left[ \omega^{(2m)}_{\kappa,2i} \right] \right)\cdot \prod_{i=1}^{m}\left(1 - \ee \left[ -\omega^{(2m)}_{\kappa,2i-1} \right] \right)^{-1}\\
&=&\frac{1}{d_{2m}} p_{2m}\left( \ee \left[ -\omega^{(2m)}_{\kappa,1} \right] \right) .
\end{eqnarray*}
\qed
\end{pf}
\begin{lem}\label{lem : induction for Euler}
The $(i,j)$-entry $\chi^{(n)}_{i,j}$ of the matrix $\chi^{(n)}$ is given by
\begin{equation}
\chi^{(n)}_{i,j}=\dfrac{1}{d_{n}}\sum_{a=1}^{d_{n}}p_{n} \left( \ee \left[ \frac{a}{d_{n}} \right] \right) \ee \left[ \frac{a}{d_{n}}(i-j) \right] .
\end{equation}
\end{lem}
\begin{pf}
Suppose that $n=2m-1,~m\in\ZZ_{\ge 1}$. 
Since the Poincar\'{e} polynomial $p'_n(t)$ of the graded ring $\CC[z_{2} ,z_{4}, \dots , z_{2m-2}]/(z_{2}^{a_{2}},z_{4}^{a_{4}},\dots,z_{2m-2}^{a_{2m-2}})$ with 
respect to the degrees $\deg(z_{2i})=d_{2i-1}$, $i=1,\dots,m$ is given by
\[
p'_n(t)=\prod_{i=1}^{m-1}\frac{1-t^{d_{2i}}}{1-t^{d_{2i-1}}},
\]
we have
\[
\dim_{\CC}\left( \CC[z_{2} ,z_{4}, \dots , z_{2m-2}]/(z_{2}^{a_{2}},z_{4}^{a_{4}},\dots,z_{2m-2}^{a_{2m-2}}) \right)_{l} 
 = \displaystyle \dfrac{1}{d_{n}}\sum_{a=1}^{d_{n}}p'_{n} \left( \ee \left[ \frac{a}{d_{n}} \right] \right) \ee \left[ -\frac{a}{d_{n}}l \right].
\]
Since $p_n(t)=(1-t)p'_n(t)$, we obtain the statement. 

Next, suppose that $n=2m,~m\in\ZZ_{\ge 1}$. 
The Poincar\'{e} polynomial of the graded ring $\CC[z_{1} ,z_{3}, \dots , z_{2m-1}]/(z_{1}^{a_{1}},z_{3}^{a_{3}},\dots,z_{2m-1}^{a_{2m-1}})$ with respect to the degrees 
$\deg(z_{2i-1})=d_{2i-2}$, $i=1,\dots,m$ is given by $p_n(t)$.
We obtain the statement since 
\[
\dim_{\CC}\left( \CC[z_{1} ,z_{3}, \dots , z_{2m-1}]/(z_{1}^{a_{1}},z_{3}^{a_{3}},\dots,z_{2m-1}^{a_{2m-1}}) \right)_{l} 
 = \displaystyle \dfrac{1}{d_{n}}\sum_{a=1}^{d_{n}}p_{n} \left( \ee \left[ \frac{a}{d_{n}} \right] \right) \ee \left[ -\frac{a}{d_{n}}l \right].
\]
\qed
\end{pf}
Now, we give a proof of the equality \eqref{eqn : Main theorem 2}.
For $n=0$, it is clear since ${\rm ch}_\Gamma^{(0)}=1$, $\widetilde Q^{(0)}=0$, $\eta^{(0)}=1$ and $\chi^{(0)}=1$.
\begin{lem}
The equality \eqref{eqn : Main theorem 2} holds for $n=1$. 
\end{lem}
\begin{pf}
The $(i,j)$-entry of the LHS of \eqref{eqn : Main theorem 2} is given by 
\[
\frac{1}{2\pi}\sum_{\kappa=1}^{a_1-1}c^{(1)}_{\kappa}\ee\left[\frac{1}{2}\left(\omega^{(1)}_{\kappa,l}-\frac{1}{2}\right)\right]\frac{1}{d_{1}}c^{(1)}_{d_{1}-\kappa}
 \ee \left[\omega^{(1)}_{\kappa,1}(i-j) \right]=
\frac{1}{d_{1}} \sum_{\kappa=1}^{a_1-1}p_{1}\left( \ee \left[ \omega_{\kappa,1}^{(1)} \right] \right)  \ee \left[ \omega_{\kappa,1}^{(1)}(i-j) \right],
\]
which is equal to $\chi^{(1)}_{i,j}$ by Lemma~\ref{lem : induction for Euler}.
\qed
\end{pf}
For each $n\in\ZZ_{\ge 0}$ and $l\in\ZZ$, define $X^{(n)}_l\in\CC$ by 
\begin{equation}
X^{(n)}_{l}:=\dfrac{1}{d_{n}}\sum_{a=1}^{d_{n}}p_{n} \left( \ee \left[\frac{a}{d_{n}} \right] \right) \ee \left[\frac{a}{d_{n}}l \right]\ (n\ge 1), \quad X^{(0)}_{l}:=0 . 
\end{equation}
\begin{lem}\label{lem : X and chi}
Suppose that $i\le j$ for $i,j=1,\dots,\widetilde{\mu}_n$. 
Then, we have 
\[
X^{(n)}_{l}=\chi^{(n)}_{i,j} ,
\]
where $l\equiv i-j ~({\rm mod}~d_n)$. 
\end{lem}
\begin{pf}
The statement follows from Lemma \ref{lem : induction for Euler} and the definition of $X^{(n)}_{l}$. 
\qed
\end{pf}
\begin{lem}\label{lem : induction for X}
For $n\ge 2$ and $l\in\ZZ$, we have 
\begin{equation}
X^{(n)}_{l} = \dfrac{1}{d_{n}}\sum_{\kappa\in I'_n}p_{n}\left( \ee \left[ \omega_{\kappa,1}^{(n)} \right] \right) 
\ee \left[ \omega_{\kappa,1}^{(n)}\cdot l \right] + \dfrac{1}{a_{n}}X^{(n-2)}_{l} .
\end{equation}
\end{lem}
\begin{pf}
Note that 
\begin{equation}
p_{n}(t)= p_{n-2}(t)\cdot\dfrac{1-t^{d_{n-1}}}{1-t^{d_{n-2}}} = p_{n-2}(t)\cdot\sum_{b=1}^{d_{n-1}}(t^{d_{n-2}})^{b}.
\end{equation}
Therefore, if $a_{n}|a$, $a_{n-1}\nmid a$ we have $p_{n}( \ee [ a/d_n])=0$ and 
\begin{align*}
X^{(n)}_{l} = & \dfrac{1}{d_{n}}\sum_{a=1}^{d_{n}}p_{n}\left( \ee \left[ \dfrac{a}{d_{n}} \right] \right) \ee \left[ -\dfrac{a}{d_{n}}l \right] \\
= & 
\dfrac{1}{d_{n}}\sum_{\kappa\in I'_n}p_{n}\left( \ee \left[ \dfrac{\kappa}{d_{n}} \right] \right) \ee \left[ \dfrac{\kappa}{d_{n}}l \right] 
+ \dfrac{1}{d_{n-2}}\dfrac{1}{a_{n}}\sum_{a'=1}^{d_{n-2}}p_{n-2}\left( \ee \left[ \dfrac{a'}{d_{n-2}} \right] \right) \ee \left[ \dfrac{a'}{d_{n-2}}l \right] \\
= & 
\dfrac{1}{d_{n}}\sum_{\kappa\in I'_n}p_{n}\left( \ee \left[ \omega_{\kappa,1}^{(n)} \right] \right) \ee \left[ \omega_{\kappa,1}^{(n)}\cdot l \right]
 + \dfrac{1}{a_{n}}X^{(n-2)}_{l}.
\end{align*}
\qed
\end{pf}
By Lemma \ref{lem : calculation of coefficients of LHS}, Lemma \ref{lem : X and chi} and Lemma \ref{lem : induction for X}, the $(i,j)$-entry of the LHS of \eqref{eqn : Main theorem 2} is given by 
\begin{align*}
& \sum_{\kappa\in I'_n}\frac{1}{(2\pi)^{n}}c^{(n)}_{\kappa}\ee\left[\frac{1}{2}\sum_{l=1}^{n}\left(\omega^{(n)}_{\kappa,l}-\frac{1}{2}\right)\right]\frac{1}{d_{n}}c^{(n)}_{d_{n}-\kappa}\ee\left[(-1)^{n-1}\omega_{\kappa,1}^{(n)}(i-j)\right] + \frac{1}{a_n}X^{(n-2)}_{i-j}\\
= & \frac{1}{d_{n}}\sum_{\kappa\in I'_n}p_{n}\left( \ee \left[ (-1)^{n-1}\omega_{\kappa,1}^{(n)} \right] \right)\ee\left[(-1)^{n-1}\omega_{\kappa,1}^{(n)}(i-j)\right] + \frac{1}{a_n}X^{(n-2)}_{i-j}\\
= & \frac{1}{d_{n}}\sum_{\kappa\in I'_n}p_{n}\left( \ee \left[ \omega_{\kappa,1}^{(n)} \right] \right)\ee\left[\omega_{\kappa,1}^{(n)}(i-j)\right] + \frac{1}{a_n}X^{(n-2)}_{i-j}
=  \chi^{(n)}_{i,j} .
\end{align*}

Hence, we have finished the proof of Theorem \ref{thm : main theorem}. \qed
\section{Proof of Theorem \ref{thm : Integral structure}}\label{sec : proof of integral structure}
We shall prove Theorem \ref{thm : Integral structure} by induction on $n$. 
\begin{lem}\label{lem : integral structure for n=1}
Theorem \ref{thm : Integral structure} holds for $n=1$. 
\end{lem}
\begin{pf}
For $j=1,\dots,a_1-1$, set
\[
\gamma_j:(-\infty,\infty)\longrightarrow\CC,\quad
\gamma_j(t):=
\begin{cases}
{\bf e}\left[-\frac{j-1}{a_1}\right]\cdot t, & t\ge0, \\
{\bf e}\left[-\frac{j}{a_1}\right]\cdot t, & t\le0 .
\end{cases}
\]
and denote by $\Gamma^{(1)}_j$ the image of the map $\gamma_j:(-\infty,\infty)\longrightarrow\CC$.
Then, $\Gamma^{(1)}_j$ defines a relative homology class of $H_1(\CC^1,{\rm Re}(\widetilde f_1)\gg0;\ZZ)$ and $\{\Gamma^{(1)}_1,\dots,\Gamma^{(1)}_{a_1-1}\}$ generates $H_1(\CC^1,{\rm Re}(\widetilde f_1)\gg0;\ZZ)$.
For $\kappa\in I_1$, namely, for $\kappa=1,\dots,a_1-1$, we have
\begin{align*}
{\rm ch}^{(1)}_{\Gamma,\kappa j}
& = \Gamma\left(1-\omega_{\kappa,1}^{(1)}\right)\left(1-{\bf e}[\omega^{(1)}_{\kappa,1}]\right){\bf e}[\omega^{(1)}_{\kappa,1}(j-1)] \\
& = a_1 \left(1-{\bf e}\left[\frac{\kappa}{a_1}\right]\right){\bf e}\left[\frac{\kappa}{a_1}(j-1)\right] \int_0^\infty e^{-x_1^{a_1}} x_1^{a_1-1-\kappa}dx_1 \\
& = \eta^{a_1-\kappa,\kappa}\int_{\Gamma^{(1)}_j} e^{-\widetilde f_1} \zeta^{(1)}_{a_1-\kappa}
\end{align*}
By the definition of $\Gamma^{(1)}_j$, the action $x_1\mapsto{\bf e}[1/a_1]x_1$ yields $\Gamma^{(1)}_j\mapsto \Gamma^{(1)}_{j-1}$.
On the other side, since $E_{j-1}=E_j(\vec{z}_1)$ on ${\rm HMF}^{L_{f_1}}_{S_1}(f_1)$, the action induced by $(\vec{z}_1)$ is ${\rm ch}^{(1)}_{\Gamma,j}\mapsto{\rm ch}^{(1)}_{\Gamma,j-1}$.
\qed
\end{pf}
\begin{lem}
Theorem \ref{thm : Integral structure} holds for $n=2$. 
\end{lem}
\begin{pf}
Consider the projection
\[
\widetilde f_2^{-1}(1)=\{(x_1,x_2)\,|\, x_1^{a_1}+x_1x_2^{a_2}=1\}\longrightarrow \CC,\quad (x_1,x_2)\mapsto x_1,
\]
which is a branched covering with branching points $\ee[k/a_1]$, $k=0,\dots, a_1-1$.
Let $\epsilon$ be a sufficiently small positive number and put 
\begin{eqnarray*}
l'_{0}&:=&\left\{(x_1,x_2)\in\CC^2\,\middle|\,x_1\in[\epsilon,1],~x_2=\left(\frac{1-x_1^{a_1}}{x_1}\right)^{\frac{1}{a_2}}\right\}, \\
l'_{1}&:=&\left\{(x_1,x_2)\in\CC^2\,\middle|\,x_1\in[\epsilon,1],~x_2={\bf e}\left[-\frac{1}{a_2}\right]\left(\frac{1-x_1^{a_1}}{x_1}\right)^{\frac{1}{a_2}}\right\}.
\end{eqnarray*}
There exists a counter-clockwise loop $C_\epsilon=\{x_1(\theta)\in\CC\,|\,0\le\theta\le1\}$ around $x_1=0$ with starting and end point $x_1=\epsilon$ such that its lift in $\widetilde f_2^{-1}(1)$ is given by 
\[
C'_\epsilon:=\left\{\left(x_1(\theta),{\bf e}\left[-\frac{\theta}{a_2}\right]\left(\frac{1-\epsilon^{a_1}}{\epsilon}\right)^{\frac{1}{a_2}}\right)\in\CC^2\,\middle|\,0\le\theta\le1\right\},
\]
and $\lim_{\epsilon\to0}x_1(\theta)=0$ for all $0\le\theta\le1$.
\begin{figure}[h]
\centering
\includegraphics[pagebox=cropbox]{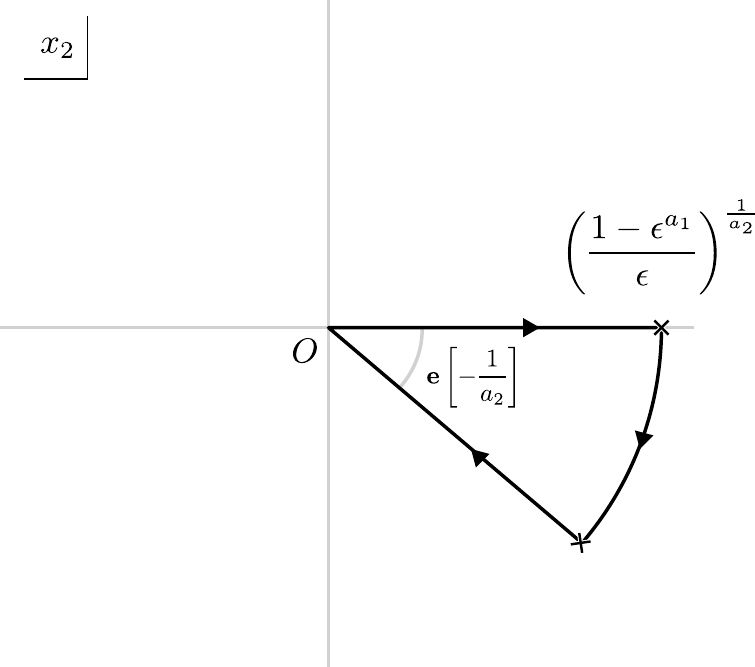}
\caption{}
\label{fig : cycle n=2}
\end{figure}

Then the cycle $S^{1}_1:=l'_1\circ C'_\epsilon \circ (l'_0)^{-1}$ in $\widetilde f_2^{-1}(1)$ defines a homology class $L_1^{(2)}$ in $H_{1}(\widetilde f_2^{-1}(1);\ZZ)$ and for the corresponding relative homology class $\Gamma^{(2)}_1\in H_2(\CC^2,{\rm Re}(\widetilde f_2)\gg0;\ZZ)$ we have
\[
\int_{\Gamma^{(2)}_1} e^{-\widetilde f_2({\bf x})} {\bf x}^{\bf k}d{\bf x} =
\lim_{\epsilon\to 0}\left(\int_\epsilon^\infty\!\!\!\int_0^\infty-\int_\epsilon^\infty\!\!\!\int_0^{\ee[-1/a_2]\infty}
+\int_{D_\epsilon}\right)e^{-\widetilde{f}_2({\bf x})} {\bf x}^{\bf k}d{\bf x},
\]
where $D_\epsilon:=\{(x_1(\theta), r\ee[-\theta/a_2])\in\CC^2\,|\,0\leq \theta \leq 1,\ r\in\RR_{\ge 0}\}$.
By a direct calculation, it follows that
\[
\int_{\Gamma^{(2)}_1} e^{-\widetilde f_2(x_1,x_2)} x_1^{k_1}x_2^{k_2}dx_1dx_2=
(1-{\bf e}[-\omega^{(2)}_{{\bf k},2}])\int_{(\RR_{\ge 0})^2} e^{-\widetilde f_2(x_1,x_2)} x_1^{k_1}x_2^{k_2}dx_1dx_2
\]
if $0\le k_1\le a_1-1$ and $0\le k_2 \le a_2-2$ and
\[
\int_{\Gamma^{(2)}_1} e^{-\widetilde f_2(x_1,x_2)}(-x_2^{a_2-1})dx_1dx_2=-\frac{2\pi\sqrt{-1}}{a_2}.
\]
Therefore, by Corollary \ref{cor : gamma integral}, for $\kappa\in I_2$ we obtain the equality 
\[{\rm ch}^{(2)}_{\Gamma,\kappa 1}=\sum_{\lambda\in I_2}\eta^{\lambda\kappa}\int_{\Gamma^{(2)}_1} e^{-\widetilde f_2} \zeta^{(2)}_\lambda.
\]
Define a homology class $L_j^{(2)}\in H_{1}(\widetilde f_2^{-1}(1);\ZZ)$ as the image of $L_1^{(2)}$ by the action $(x_1,x_2)\mapsto \left(\ee\left[\frac{(j-1)}{a_1}\right]x_1,\ee\left[\frac{-(j-1)}{a_1a_2}\right]x_2\right)$ on $\widetilde f_2^{-1}(1)$ and 
$\Gamma^{(2)}_j\in H_2(\CC^2,{\rm Re}(\widetilde f_2)\gg0;\ZZ)$ as the corresponding class of $L_j^{(2)}$.
Hence, we have the equality 
\[
{\rm ch}^{(2)}_{\Gamma,\kappa j}=\sum_{\lambda\in I_2}\eta^{\lambda\kappa}\int_{\Gamma^{(2)}_j} e^{-\widetilde f_2} \zeta^{(2)}_\lambda
\]
for $j=1,\dots,\widetilde\mu_2$.
The statement follows.
\qed
\end{pf}
Next, we show the general case.
For simplicity, we put ${\bf x}_{n}':=(x_2,x_3,\dots,x_{n})$ and ${\bf x}_{n}'':=(x_3,x_4,\dots,x_{n})$ for ${\bf x}={\bf x}_n=(x_1,\dots,x_{n})$, respectively.
\begin{rem}\label{rem : construction of cycle}
Let $S^{n-1}_1$ be an $(n-1)$-dimensional sphere in $\widetilde f_n^{-1}(1)$ whose homology class defines $L_1^{(n)}\in H_{n-1}(\widetilde f_n^{-1}(1);\ZZ)$.
Since $\widetilde f_n$ is weighted homogeneous, for each $w\in\RR_{>0}$ there exists an $(n-1)$-dimensional sphere $S_w^{n-1}$ in $\widetilde f_{n}^{-1}(w)$ whose homology class defines $L^{(n)}_1(w)\in H_{n-1}(\widetilde f^{-1}_{n}(w);\ZZ)$ which is the image of $L^{(n)}_1$ by the isomorphism $H_{n-1}(\widetilde{f}_{n}^{-1}(1);\ZZ)\cong H_{n-1}(\widetilde{f}_{n}^{-1}(w);\ZZ)$ induced by the parallel transformation.
\end{rem}
We prove the case of $n=2m-1$, $m\in\ZZ_{\ge2}$.
If $x_1\ne0$, by using the change of variables 
\[
y_2^{a_2}=x_1x_2^{a_2},\quad y_{i-1}y_i^{a_i}=x_{i-1}x_i^{a_i}, i=3,4,\dots, 2m-1,
\]
we have $\widetilde f_{2m-1}({\bf x}_{2m-1})=x_1^{a_1}+\widetilde f_{2m-2}({\bf y}'_{2m-1})$. 
\begin{lem}\label{lem : construction of cycle n=2m-1}
Suppose that there exists a $(2m-3)$-dimensional sphere $S_1^{2m-3}$ in $\widetilde f_{2m-2}^{-1}(1)$ such that its homology class defines $L^{(2m-2)}_1\in H_{2m-3}(\widetilde f^{-1}_{2m-2}(1);\ZZ)$ and for the corresponding class $\Gamma^{(2m-2)}_1$ we have 
\[
{\rm ch}^{(2m-2)}_{\Gamma,\kappa 1}=\sum_{\lambda\in I_{2m-2}}\eta^{\lambda\kappa}\int_{\Gamma^{(2m-2)}_1} e^{-\widetilde f_{2m-2}} \zeta^{(2m-2)}_\lambda.
\]
Then there exists an element $L^{(2m-1)}_1\in H_{2m-2}(\widetilde f^{-1}_{2m-1}(1);\ZZ)$ whose corresponding class $\Gamma^{(2m-1)}_1\in H_{2m-1}(\CC^{2m-1},{\rm Re}(\widetilde f_{2m-1})\gg0;\ZZ)$ satisfies that 
\begin{equation}\label{eqn : Lefschetz thimble for n=2m-1}
\begin{split}
& \int_{\Gamma^{(2m-1)}_1}e^{-\widetilde f_{2m-1}({\bf x}_{2m-1})}x_1^{k_1}\phi({\bf x}_{2m-1}')d{\bf x}_{2m-1} \\
= & \int_{\Gamma^{(1)}_1}e^{-x_1^{a_1}}x_1^{k_1+\bar{\phi}}dx_1\cdot
\int_{\Gamma^{(2m-2)}_1}e^{-\widetilde f_{2m-2}({\bf y}_{2m-1}')}\phi({\bf y}_{2m-1}')d{\bf y}_{2m-1}' ,
\end{split}
\end{equation}
for $x_1^{k_1}\cdot\phi({\bf x}_{2m-1}')\in B_{\widetilde f_{2m-1}}$, where $\bar{\phi}$ is a rational number given by $\phi({\bf x}_{2m-1}')=x_1^{\bar{\phi}}\cdot\phi({\bf y}_{2m-1}')$. 
\end{lem}
\begin{pf}
We construct a $(2m-2)$-dimensional sphere so that it gives a homology class $L_1\in H_{2m-2}(\widetilde f^{-1}_{2m-1}(1);\ZZ)$ in the statement.
Consider the projection 
\[
\widetilde f_{2m-1}^{-1}(1)\longrightarrow\CC,\quad (x_1,\dots,x_{2m-1}) \mapsto x_1.
\]
Let $\epsilon$ be a sufficiently small positive number, 
$l_{0}:=\{t\in\CC\,|\,t\in[\epsilon,1]\},~l_{1}:=\{{\bf e}[-1/a_1]t\in\CC\,|\,t\in[\epsilon,1]\}$ line segments, 
and $C_\epsilon:=\{\epsilon{\bf e}[-\theta/a_1]\in\CC\,|\,0\le\theta\le1\}$.
Set $A^{(2m-1)}_1:=l_0 \circ C_\epsilon \circ l_1^{-1}$ (see Figure \ref{fig : n=2m-1}).
\begin{figure}[h]
\centering
\includegraphics[pagebox=cropbox]{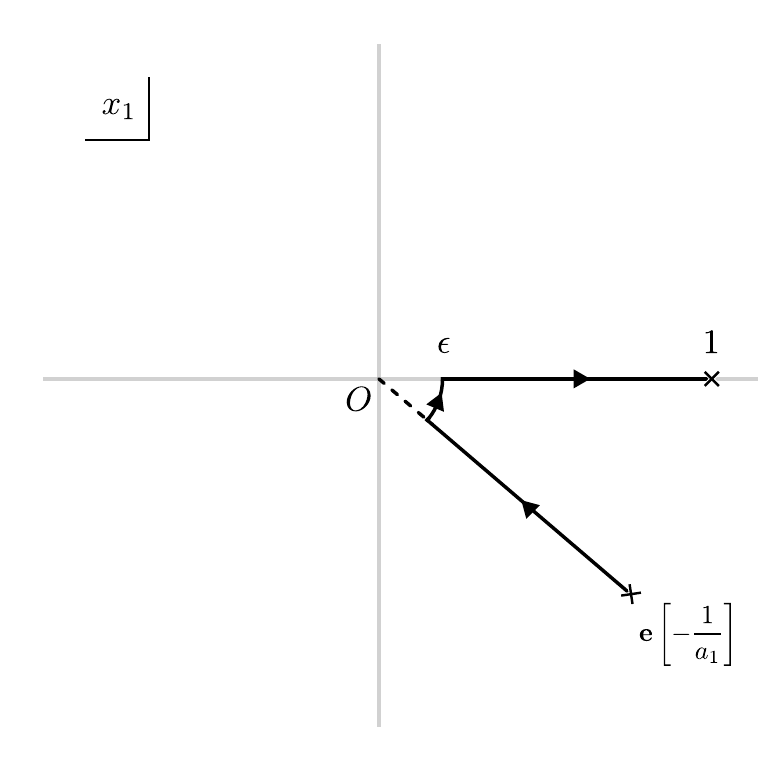}
\caption{Domain $A^{(2m-1)}_1$.}
\label{fig : n=2m-1}
\end{figure}

Since $x_1\in A^{(2m-1)}_1$ implies $x_1\ne0$, we have $\widetilde f_{2m-1}({\bf x}_{2m-1})=x_1^{a_1}+\widetilde f_{2m-2}({\bf y}'_{2m-1})$.
Since $S^{2m-3}_{1-x_1^{a_1}}$ vanishes on the fiber over $x_1=1,{\bf e}[-1/a_1]$, $S^{2m-2}_1:=\bigcup_{x_1\in A^{(2m-1)}_1}S^{2m-3}_{1-x_1^{a_1}}$ is homeomorphic to the $(2m-2)$-dimensional sphere by the assumption and Remark \ref{rem : construction of cycle}.

Define $L_1^{(2m-1)}\in H_{2m-2}(\widetilde f_{2m-1}^{-1}(1);\ZZ)$ by the homology class of $\bigcup_{x_1\in A_1}S^{2m-3}_{1-x_1^{a_1}}$.
For the corresponding class $\Gamma^{(2m-1)}_1$, we have 
\begin{align*}
& \int_{\Gamma^{(2m-1)}_1}e^{-\widetilde f_{2m-1}({\bf x}_{2m-1})}x_1^{k_1}\phi({\bf x}_{2m-1}')d{\bf x}_{2m-1} \\
= & \lim_{\epsilon\to0}
\left(\int_{\epsilon}^{\infty}-\int_\epsilon^{\ee[-1/a_1]\infty}+\int_{C_\epsilon}\right)e^{-x_1^{a_1}}x_1^{k_1+\bar{\phi}}dx_1\cdot
\int_{\Gamma^{(2m-2)}_1}e^{-\widetilde f_{2m-2}({\bf y}_{2m-1}')}\phi({\bf y}_{2m-1}')d{\bf y}_{2m-1}' ,
\end{align*}
It gives the equation \eqref{eqn : Lefschetz thimble for n=2m-1}.
\qed
\end{pf}
By Lemma \ref{lem : construction of cycle n=2m-1}, for ${\bf k}=(k_1,\dots,k_{2m-1})$ satisfying ${\bf x}_{2m-1}^{\bf k}\in B'_{\widetilde f_{2m-1}}$ we obtain 
\begin{eqnarray*}
& & \int_{\Gamma^{(2m-1)}_1}e^{-\widetilde f_{2m-1}({\bf x}_{2m-1})}{\bf x}_{2m-1}^{\bf k}d{\bf x}_{2m-1} \\ 
& = & \int_{\Gamma^{(1)}_1}e^{-x_1^{a_1}}x_1^{a_1\omega^{(2m-1)}_{{\bf k},1}-1}dx_1 \cdot
\int_{\Gamma^{(2m-2)}_1}e^{-\widetilde f_{2m-2}({\bf y}_{2m-1}')}({\bf y}_{2m-1}')^{{\bf k}'}d{\bf y}_{2m-1}' \\
& = & \frac{1}{a_1}\Gamma(\omega^{(2m-1)}_{{\bf k},1})(1-{\bf e}[-\omega^{(2m-1)}_{{\bf k},1}]) 
\int_{\Gamma^{(2m-2)}_1}e^{-\widetilde f_{2m-2}({\bf y}_{2m-1}')}({\bf y}_{2m-1}')^{{\bf k}'}d{\bf y}_{2m-1}' ,
\end{eqnarray*}
where ${\bf k}':=(k_2,\dots,k_{2m-1})$ and $({\bf y}_{2m-1}')^{{\bf k}'}:=y_2^{k_2}y_3^{k_3}\dots y_{2m-1}^{k_{2m-1}}$.
For each $x_1^{k_1}\cdot\phi({\bf x}_{2m-3}')\in B_{\widetilde f_{2m-3}}$, it inductively follows that 
\begin{eqnarray*}
& & \int_{\Gamma^{(2m-1)}_1}e^{-\widetilde f_{2m-1}({\bf x}_{2m-1})}x_1^{k_1}\phi({\bf x}_{2m-3}')x_{2m-1}^{a_{2m-1}-1}d{\bf x}_{2m-1} \\
& = & \int_{\Gamma^{(1)}_1}e^{-x_1^{a_1}}x_1^{k_1+\bar{\phi}}dx_1\cdot
\int_{\Gamma^{(2m-2)}_1}e^{-\widetilde f_{2m-2}({\bf y}_{2m-1}')}\phi({\bf y}'_{2m-3})y_{2m-1}^{a_{2m-1}-1}d{\bf y}'_{2m-1} \\
& = & \int_{\Gamma^{(1)}_1}e^{-x_1^{a_1}}x_1^{k_1+\bar{\phi}}dx_1\cdot
\left(-\frac{2\pi\sqrt{-1}}{a_{2m-1}}\right)\int_{\Gamma^{(2m-4)}_1}e^{-\widetilde f_{2m-4}({\bf y}_{2m-3}')}\phi({\bf y}_{2m-3}')d{\bf y}_{2m-3}' \\
& = & -\frac{2\pi\sqrt{-1}}{a_{2m-1}}\int_{\Gamma^{(2m-3)}_1}e^{-\widetilde f_{2m-3}({\bf x}_{2m-3})}x_1^{k_1}\phi({\bf x}_{2m-3}')d{\bf x}_{2m-3} .
\end{eqnarray*}
Therefore, for $\kappa \in I_{2m-1}$ we obtain 
\[
{\rm ch}^{(2m-1)}_{\Gamma,\kappa 1}=\sum_{\lambda\in I_{2m-1}}\eta^{\lambda\kappa}\int_{\Gamma^{(2m-1)}_1} e^{-\widetilde f_{2m-1}} \zeta^{(2m-1)}_\lambda.
\]
Define a homology class $L_j^{(2m-1)}\in H_{2m-2}(\widetilde f_{2m-1}^{-1}(1);\ZZ)$ as the image of $L_1^{(2m-1)}$ by the action 
$(x_1,\dots,x_k,\dots,x_{2m-1})\mapsto \left(\ee\left[-\frac{(j-1)}{a_1}\right]x_1,\dots,\ee\left[(-1)^{k}\frac{(j-1)}{d_{k}}\right]x_k,\dots,\ee\left[-\frac{(j-1)}{d_{2m-1}}\right]x_{2m-1}\right)$
on $\widetilde f_{2m-1}^{-1}(1)$ and $\Gamma^{(2m-1)}_j\in H_{2m-1}(\CC^{2m-1},{\rm Re}(\widetilde f_{2m-1})\gg0;\ZZ)$ as the corresponding class of $L_j^{(2m-1)}$.
For $j=1,\dots,\widetilde\mu_{2m-1}$, we have the equality 
\[
{\rm ch}^{(2m-1)}_{\Gamma,\kappa j}=\sum_{\lambda\in I_{2m-1}}\eta^{\lambda\kappa}\int_{\Gamma^{(2m-1)}_j} e^{-\widetilde f_{2m-1}} \zeta^{(2m-1)}_\lambda.
\]
Finally, we prove the case of $n=2m$, $m\in\ZZ_{\ge2}$.
If $x_2\ne0$, by using the change of variables 
\[
y_3^{a_3}=x_2x_3^{a_3},\quad y_{i-1}y_i^{a_i}=x_{i-1}x_i^{a_i}, i=4,5,\dots, 2m,
\]
we have $\widetilde f_{2m}({\bf x}_{2m})=x_1^{a_1}+x_1x_2^{a_2}+\widetilde f_{2m-2}({\bf y}''_{2m})$. 
\begin{lem}\label{lem : construction of cycle n=2m}
If the assumption in Lemma \ref{lem : construction of cycle n=2m-1} holds, then there exists an element $L^{(2m)}_1\in H_{2m-1}(\widetilde f^{-1}_{2m}(1);\ZZ)$ whose corresponding class $\Gamma^{(2m)}_1\in H_{2m}(\CC^{2m},{\rm Re}(\widetilde f_{2m})\gg0;\ZZ)$ satisfies that  
\begin{equation}\label{eqn : Lefschetz thimble for n=2m}
\begin{split}
& \int_{\Gamma^{(2m)}_1}e^{-\widetilde f_{2m}({\bf x}_{2m})}x_1^{k_1}x_2^{k_2}\phi({\bf x}_{2m}'')d{\bf x}_{2m} \\
= & \int_{\Gamma^{(2)}_1}e^{-(x_1^{a_1}+x_1x_2^{a_2})}x_1^{k_1}x_2^{k_2+\bar{\phi}}dx_1dx_2\cdot
\int_{\Gamma^{(2m-2)}_1}e^{-\widetilde f_{2m-2}({\bf y}_{2m}'')}\phi({\bf y}_{2m}'')d{\bf y}_{2m}'' ,
\end{split}
\end{equation}
for $x_1^{k_1}x_2^{k_2}\cdot\phi({\bf x}_{2m}'')\in B_{\widetilde f_{2m}}$, where $\bar{\phi}$ is a rational number given by $\phi({\bf x}_{2m}'')=x_2^{\bar{\phi}}\cdot\phi({\bf y}_{2m}'')$.
\end{lem}
\begin{pf}
We construct a $(2m-1)$-dimensional sphere $S^{2m-1}$ whose homology class $L_1\in H_{2m-1}(\widetilde f^{-1}_{2m}(1);\ZZ)$ in the statement by gluing two $(2m-1)$-dimensional disks $D_1^{2m-1}$ and $D_2^{2m-1}$ along the boundary. 

First, we construct a disk $D_1^{2m-1}$ as follows.
Consider the projection 
\[
\widetilde f_{2m}^{-1}(1)\longrightarrow\CC,\quad (x_1,\dots,x_{2m}) \mapsto x_1.
\]
Let $\epsilon_1$ be a sufficiently small positive number, 
$l_{0}:=\{t\in\CC\,|\,t\in[\epsilon_1,1]\}$ a line segment. 
Since $x_1\in l_0$ implies $x_1\ne0$, we have $\widetilde f_{2m}({\bf x}_{2m})=x_1^{a_1}+\widetilde f_{2m-1}({\bf y}_{2m}')$.
Therefore, since $S^{2m-2}_{1-x_1^{a_1}}$ vanishes on the fiber over $x_1=1$, by Remark \ref{rem : construction of cycle} and Lemma \ref{lem : construction of cycle n=2m-1} $D_1^{2m-1}:=\bigcup_{x_1\in l_0}S^{2m-2}_{1-x_1^{a_1}}$ is homeomorphic to the $(2m-1)$-dimensional disk.

Next, we construct a disk $D_2^{2m-1}$ as follows.
Consider the projection 
\[
\widetilde f_{2m}^{-1}(1)\longrightarrow\CC^2,\quad (x_1,\dots,x_{2m}) \mapsto (x_1,x_2).
\]
Let $\epsilon_2$ be a sufficiently small positive number.
There exists a counter-clockwise loop $C_{\epsilon_1}=\{x_1(\theta)\in\CC\,|\,0\le\theta\le1\}$ around $x_1=0$ with starting and end point $x_1=\epsilon_1$ such that for each $0\le\theta\le1$, $\lim_{\epsilon_1\to0}x_1(\theta)=0$ and the set $D^{2m-2}_\theta:=\bigcup_{(x_1,x_2)\in l'_\theta}S^{2m-3}_{1-(x_1^{a_1}+x_1x_2^{a_2})}$ is a subset of $\widetilde f_{2m}^{-1}(1)$, where $l'_{\theta}$ is given by
\[
l'_{\theta}:=\left\{\left(x_1(\theta),t{\bf e}\left[-\frac{\theta}{a_2}\right]\right)\in\CC^2\,\middle|\,t\in\left[\epsilon_2,\left(\frac{1-\epsilon_1^{a_1}}{\epsilon_1}\right)^{\frac{1}{a_2}}\right]\right\}.
\]
\begin{figure}[h]
\centering
\includegraphics[pagebox=cropbox]{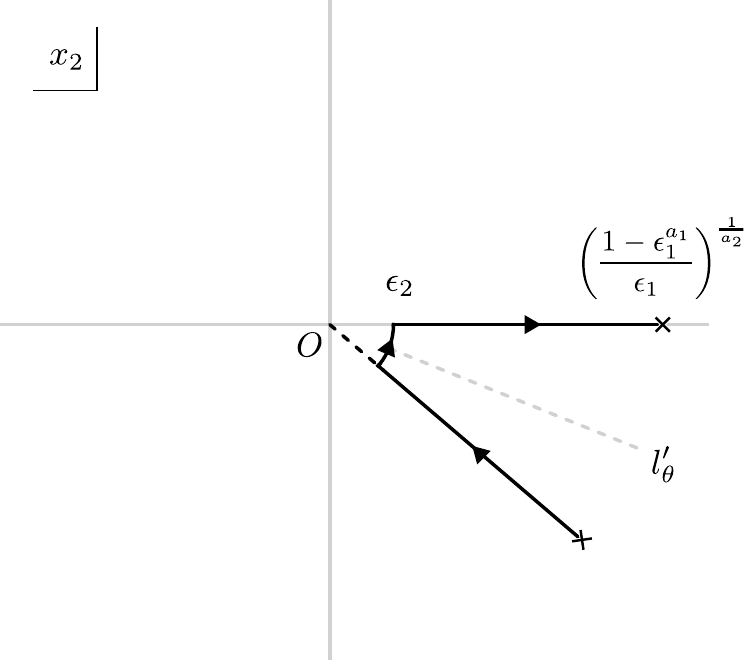}
\caption{}
\label{fig : n=2m}
\end{figure}

Then, for each $0\le\theta\le1$, $D^{2m-2}_\theta$ is homeomorphic to the $(2m-2)$-dimensional disk.
Therefore, $D_2^{2m-1}:=\bigcup_{0\le\theta\le1}D^{2m-2}_\theta$ is homeomorphic to the $(2m-1)$-dimensional disk. 

By taking $\epsilon_2$ as $\epsilon$ in Lemma \ref{lem : construction of cycle n=2m-1} for $w=1-\epsilon_1^{a_1}$, the boundary of $D_2^{2m-1}$ coincides with the one of $D_1^{2m-1}$.
Hence, we obtain a $(2m-1)$-dimensional sphere $S^{2m-1}$ by gluing $D_1^{2m-1}$ and $D_2^{2m-1}$ on their common boundaries.
Define $L^{(2m)}_1\in H_{2m-1}(\widetilde f_{2m}^{-1}(1);\ZZ)$ by the homology class of this $(2m-1)$-dimensional sphere $S^{2m-1}$.
For the corresponding class $\Gamma^{(2m)}_1$, we have 
\begin{align*}
& \int_{\Gamma^{(2m)}_1}e^{-\widetilde f_{2m}({\bf x}_{2m})}x_1^{k_1}x_2^{k_2}\phi({\bf x}_{2m}'')d{\bf x}_{2m} \\
= & \lim_{\substack{\epsilon_1\to0 \\ \epsilon_2\to0}}\left(\int_{\epsilon_1}^\infty\!\!\!\int_{\epsilon_2}^\infty-\int_{\epsilon_1}^\infty\!\!\!\int_{\epsilon_2}^{\ee[-1/a_2]\infty}
+\int_{D_{\epsilon_1,\epsilon_2}}\right)e^{-(x_1^{a_1}+x_1x_2^{a_2})}x_1^{k_1}x_2^{k_2+\bar{\phi}}dx_1dx_2 \\
& \cdot \int_{\Gamma^{(2m-2)}_1}e^{-\widetilde f_{2m-2}({\bf y}_{2m}'')}\phi({\bf y}_{2m}'')d{\bf y}_{2m}'' ,
\end{align*}
where $D_{\epsilon_1,\epsilon_2}:=\{(x_1(\theta), r\ee[-\theta/a_2])\,|\,0\leq \theta \leq 1,\ r\ge\epsilon_2\}$.
It gives the equation \eqref{eqn : Lefschetz thimble for n=2m}.
\qed
\end{pf}
By Lemma \ref{lem : construction of cycle n=2m}, for ${\bf k}=(k_1,\dots,k_{2m})$ satisfying ${\bf x}_{2m}^{\bf k}\in B'_{\widetilde f_{2m}}$ we obtain 
\begin{eqnarray*}
& & \int_{\Gamma^{(2m)}_1}e^{-\widetilde f_{2m}({\bf x}_{2m})}{\bf x}_{2m}^{\bf k}d{\bf x}_{2m} \\
& = & \int_{\Gamma^{(2)}_1}e^{-\widetilde f_{2}(x_1,x_2)}x_1^{k_1}x_2^{a_2\omega^{(2m)}_{{\bf k},2}-1}dx_1dx_2 \cdot
\int_{\Gamma^{(2m-2)}_1}e^{-\widetilde f_{2m-2}({\bf y}_{2m}'')}({\bf y}_{2m}'')^{{\bf k}''}d{\bf y}_{2m}'' \\
& = & \frac{1}{d_2}\Gamma(\omega^{(2m)}_{{\bf k},1})\Gamma(\omega^{(2m)}_{{\bf k},2})(1-{\bf e}[-\omega^{(2m)}_{{\bf k},2}]) 
\int_{\Gamma^{(2m-2)}_1}e^{-\widetilde f_{2m-2}({\bf y}_{2m}'')}({\bf y}_{2m}'')^{{\bf k}''}d{\bf y}_{2m}'' ,
\end{eqnarray*}
where ${\bf k}'':=(k_3,\dots,k_{2m})$ and $({\bf y}_{2m}'')^{{\bf k}''}:=y_3^{k_3}y_4^{k_4}\dots y_{2m}^{k_{2m}}$. 
For each $x_1^{k_1}x_2^{k_2}\cdot\phi({\bf x}_{2m-2}'')\in B_{\widetilde f_{2m-2}}$, it inductively follows that
\begin{eqnarray*}
& & \int_{\Gamma^{(2m)}_1}e^{-\widetilde f_{2m}({\bf x})}x_1^{k_1}x_2^{k_2}\phi({\bf x}_{2m-2}'')x_{2m}^{a_{2m}-1}d{\bf x}_{2m} \\
& = & \int_{\Gamma^{(2)}_1}e^{-\widetilde f_{2}(x_1,x_2)}x_1^{k_1}x_2^{k_2+\bar{\phi}}dx_1dx_2\cdot
\int_{\Gamma^{(2m-2)}_1}e^{-\widetilde f_{2m-2}({\bf y}_{2m}'')}\phi({\bf y}_{2m-2}'')y_{2m-1}^{a_{2m-1}-1}d{\bf y}_{2m}'' \\
& = & \int_{\Gamma^{(2)}_1}e^{-\widetilde f_{2}(x_1,x_2)}x_1^{k_1}x_2^{k_2+\bar{\phi}}dx_1dx_2\cdot
\left(-\frac{2\pi\sqrt{-1}}{a_{2m}}\right)\int_{\Gamma^{(2m-4)}_1}e^{-\widetilde f_{2m-4}({\bf y}_{2m-2}'')}\phi({\bf y}_{2m-2}'')d{\bf y}_{2m-2}'' \\
& = & -\frac{2\pi\sqrt{-1}}{a_{2m}}\int_{\Gamma^{(2m-2)}_1}e^{-\widetilde f_{2m-2}({\bf x}_{2m-2})}x_1^{k_1}x_2^{k_2}\phi({\bf x}_{2m-2}'')d{\bf x}_{2m-2}.
\end{eqnarray*}
Therefore, for $\kappa \in I_{2m}$ we obtain 
\[
{\rm ch}^{(2m)}_{\Gamma,\kappa 1}=\sum_{\lambda\in I_{2m}}\eta^{\lambda\kappa}\int_{\Gamma^{(2m)}_1} e^{-\widetilde f_{2m}} \zeta^{(2m)}_\lambda.
\]
Define a homology class $L_j^{(2m)}\in H_{2m-1}(\widetilde f_{2m}^{-1}(1);\ZZ)$ as the image of $L_1^{(2m)}$ by the action 
$(x_1,\dots,x_k,\dots,x_{2m})\mapsto \left(\ee\left[\frac{(j-1)}{a_1}\right]x_1,\dots,\ee\left[(-1)^{k-1}\frac{(j-1)}{d_{k}}\right]x_k,\dots,\ee\left[-\frac{(j-1)}{d_{2m}}\right]x_{2m}\right)$
on $\widetilde f_{2m}^{-1}(1)$ and $\Gamma^{(2m)}_j\in H_{2m}(\CC^{2m},{\rm Re}(\widetilde f_{2m})\gg0;\ZZ)$ as the corresponding class of $L_j^{(2m)}$.
For $j=1,\dots,\widetilde\mu_{2m}$, we have the equality 
\[
{\rm ch}^{(2m)}_{\Gamma,\kappa j}=\sum_{\lambda\in I_{2m}}\eta^{\lambda\kappa}\int_{\Gamma^{(2m)}_j} e^{-\widetilde f_{2m}} \zeta^{(2m)}_\lambda.
\]

Hence, we have finished the proof of Theorem \ref{thm : Integral structure}. \qed
\section{Stokes and Central connection matrices for $\widetilde f_n$}\label{sec : Stokes matrix}
In this section, we shall consider a Frobenius manifold associated with $\widetilde{f}_n$. 
Denote by $\CC_w$ the complex plane whose coordinate is $w$, which we call $w$-plane for simplicity. 
Define $p:\CC^{n}\times \CC^{\widetilde{\mu}_n} \longrightarrow\CC^{\widetilde{\mu}_n}$ by the natural projection and a function $\widetilde{F}_n : \CC^{n}\times \CC^{\widetilde{\mu}_n} \longrightarrow\CC_w$ by 
\begin{equation}
\widetilde{F}_n({\bf x};{\bf s}):=\widetilde{f}_n({\bf x})+
\sum_{\kappa\in I_n}s_{\kappa}\cdot\phi^{(n)}_\kappa({\bf x}),
\end{equation}
where $\{\phi^{(n)}_\kappa({\bf x})\}_{\kappa\in I_n}$ is the basis of ${\rm Jac}(\widetilde f_n)$ given in Definition \ref{defn : monomial basis on Jac}. 
Let $(\widetilde\omega^{(n)}_1,\dots,\widetilde\omega^{(n)}_n)$ be the rational weights of $\widetilde f_n$. 
By a direct calculation, we have 
\begin{equation}
\widetilde\omega^{(n)}_{i}=\sum_{l=1}^i(-1)^{i-l}\frac{d_{l-1}}{d_i},
\quad i=1,\dots,n.
\end{equation}
Set 
\begin{equation}
q_\kappa:=\deg\,\phi^{(n)}_\kappa({\bf x})=\sum_{i=1}^n k_i\widetilde \omega^{(n)}_i,\quad \kappa\in I'_n, 
\end{equation}
where $(k_1,\dots,k_n)=\psi^{-1}(\kappa)$. 
\begin{prop}
The pair $(\widetilde{F}_n,p)$ is a universal unfolding of $\widetilde{f}_n$, namely it satisfies that 
\begin{enumerate}
\item $\widetilde{F}_n|_{p^{-1}(0)}=\widetilde{f}_n$ in a neighborhood of the origin in $p^{-1}(0)\cong\CC^{\widetilde{\mu}_n}$, 
\item there exists an $\O_{\CC^{\widetilde{\mu}_n},0}$-isomorphism 
\begin{equation}\label{eqn : KS map}
\rho:\T_{\CC^{\widetilde{\mu}_n},0}\longrightarrow p_*\O_{\CC^n\times\CC^{\widetilde{\mu}_n},0}\left/\middle(\frac{\p \widetilde{F}_n}{\p x_1},\dots,\frac{\p \widetilde{F}_n}{\p x_n}\right) ,\quad 
\delta\mapsto\widehat{\delta}\widetilde{F}_n ,
\end{equation}
where $\widehat{\delta}$ is a lifting on $\CC^n\times\CC^{\widetilde{\mu}_n}$ of a vector field $\delta\in\T_{\CC^{\widetilde{\mu}_n},0}$ with respect to the projection $p$.
\end{enumerate}
\end{prop}
Let us fix Euclidean norms $\| \cdot \|$ on $\CC^n\times\CC^{\widetilde{\mu}_n}, \CC^{\widetilde{\mu}_n}$ and $\CC_w$. 
For positive real numbers $\varepsilon_M,\varepsilon_\Delta$ and $\varepsilon_\X$, put 
\begin{subequations}\label{eq: 5.3}
\begin{equation}
M:=\{{\bf s}\in\CC^{\widetilde{\mu}_n}~|~\|{\bf s}\|<\varepsilon_M\}, 
\end{equation}
\begin{equation}
\Delta:=\{w\in\CC_w~|~\| w\|<\varepsilon_\Delta\}
\end{equation}
\begin{equation}
\X:=\{({\bf x},{\bf s})\in\CC^n\times\CC^{\widetilde{\mu}_n}~|~\| ({\bf x},{\bf s}) \|<\varepsilon_\X, \| \widetilde{F}_n({\bf x};{\bf s}) \|<\varepsilon_\Delta\}\cap p^{-1}(M). 
\end{equation}
Denote the fiber of $\X$ with respect to $p$ over ${\bf s}\in M$ by $\X_{\bf s}$ and set
\begin{equation}
Y_{w;{\bf s}}:=\X_{\bf s}\cap\widetilde{F}_n^{-1}(w),\quad w\in\CC_w,~{\bf s}\in M. 
\end{equation}
\end{subequations}
For the choice $1 \gg \varepsilon_\X \gg \varepsilon_\Delta \gg \varepsilon_M > 0$ of radius, the map $(\widetilde{F}_n,p):\X\longrightarrow\Delta\times M$ defines a fibration which the fiber over $0$ is isomorphic to the singularity $Y_{0;{\bf 0}}=\{\widetilde{f}_n=0\}$ and generic fiber is homotopic to a bouquet of $\widetilde{\mu}_n$ copies of $(n-1)$-dimensional sphere due to Milnor \cite{M}. 
In particular, for a regular value $w_0\in\p\Delta$ we have 
\[
H_n(\X_{\bf s},Y_{w_0;{\bf s}};\ZZ)\cong H_{n-1}(Y_{w_0;{\bf s}};\ZZ)\cong\ZZ^{\widetilde{\mu}_n} .
\]
For simplicity, we fix the constants $\varepsilon_M,\varepsilon_\Delta$ and $\varepsilon_\X$ suitably. 
\begin{defn}
We shall denote by $\circ$ the induced product structure on $\T_M$ by the $\O_{\CC^{\widetilde{\mu}_n},0}$-isomorphism $\rho$. 
Namely, we have 
\begin{equation}
\widehat{(\delta\circ\delta')}\widetilde{F}_n=\widehat{\delta}\widetilde{F}_n\cdot\widehat{\delta'}\widetilde{F}_n,\quad\delta,\delta\in\T_M .
\end{equation}
Define two vector fields $e$ and $E$ as follows:
\begin{enumerate}
\item The vector field $e\in\T_M$ corresponding to the unit $1$ by the $\O_{\CC^{\widetilde{\mu}_n},0}$-isomorphism $\rho$ is called the {\it unit vector filed}. 
That is, 
\[
\widehat{e}\widetilde{F}_n=1 .
\]
\item The vector field $E\in\T_M$ corresponding to $\widetilde{F}_n$ by the $\O_{\CC^{\widetilde{\mu}_n},0}$-isomorphism $\rho$ is called the {\it Euler vector filed}. 
That is, 
\[
\widehat{E}\widetilde{F}_n=\widetilde{F}_n .
\]
\end{enumerate}
\end{defn}
Note that the unit vector field and the Euler vector field are given explicitly by 
\begin{equation}
e=\frac{\p}{\p s_{\psi({\bf 0})}},\quad E=\sum_{\kappa\in I_n}(1-q_\kappa)\frac{\p}{\p s_\kappa}. 
\end{equation}
Let $\PP^1_u$ be the complex projective line whose coordinate is $u$. 
In order to define a notion of a primitive form, it is necessary to define a Saito structure associated with the universal unfolding $\widetilde{F}_n$. 
This structure is given as a tuple consisting of the filtered de Rham cohomology group $\H_{\widetilde{F}_n}$ (whose increasing filtration is denoted by $\H_{\widetilde{F}_n}^{(k)}$ $(k\in\ZZ)$), the Gau\ss--Manin connection $\nabla$ on $\H_{\widetilde{F}_n}$ and the higher residue pairing $K_{\widetilde{F}_n}$ on $\H_{\widetilde{F}_n}$. 
In the paper, we omit the details about those objects and refer the interested reader to \cite{ST}. 
Since $\widetilde f_n$ is a weighted homogeneous polynomial, there exists a canonical primitive form defined by exponents of $\widetilde f_n$. 
We will use this primitive form in the paper. 
\begin{prop}[\cite{S-K,S-M}]\label{prop : existence of primitive form}
There exists a unique primitive form $\zeta\in\Gamma(M,\H_{\widetilde{F}_n}^{(0)})$ for the tuple $(\H_{\widetilde{F}_n}^{(0)}, \nabla, K_{\widetilde{F}_n})$ with the minimal exponent 
$r=\sum_{i=1}^n\widetilde\omega^{(n)}_i=\sum_{i=1}^n\omega^{(n)}_i$ and the normalization
\begin{equation}\label{eqn : primitive form}
\zeta|_{{\bf s}=0}=\left[dx_1\wedge\cdots\wedge dx_n\right] .
\end{equation}
\qed
\end{prop}
The higher residue pairing $K_{\widetilde{F}_n}$ and the primitive form $\zeta$ induce a symmetric non-degenerate $\O_M$-bilinear form $\eta_\zeta:\T_M\times\T_M\longrightarrow\O_M$ defined by 
\begin{equation}
\eta_\zeta(\delta,\delta'):=K_{\widetilde{F}_n}\left(r^{(0)}(u\nabla_{\delta}\zeta),r^{(0)}(u\nabla_{\delta'}\zeta)\right),\quad \delta,\delta'\in\T_M .
\end{equation}
Hence, we obtain the following 
\begin{prop}[{cf.~\cite{ST}}]
The tuple $(M,\eta_\zeta,\circ,e,E)$ is a Frobenius manifold of rank $\widetilde{\mu}_n$ and dimension $d_{\widetilde f_n}:=\sum_{i=1}^n(1-2\widetilde\omega^{(n)}_i)=\sum_{i=1}^n(1-2\omega^{(n)}_i)$. 
Namely, it satisfies the following properties;
\begin{enumerate}
\item The product $\circ$ is self-adjoint with respect to $\eta$: that is,
\begin{equation*}
\eta_\zeta(\delta\circ\delta',\delta'')=\eta_\zeta(\delta,\delta'\circ\delta''),\quad
\delta,\delta',\delta''\in\T_M. 
\end{equation*} 
\item The {\rm Levi}--{\rm Civita} connection $\ns:\T_M\otimes_{\O_M}\T_M\longrightarrow\T_M$ with respect to $\eta_\zeta$ is
flat: that is, 
\begin{equation*}
[\ns_\delta,\ns_{\delta'}]=\ns_{[\delta,\delta']},\quad \delta,\delta'\in\T_M.
\end{equation*}
\item The $\O_M$-linear morphism $C_\delta:\T_M\longrightarrow\T_M,~\delta\in\T_M$ defined by $C_\delta(-):=-\circ \delta$ satisfies 
\begin{equation*}
\ns_{\delta}(C_{\delta'}\delta'')-C_{\delta'}(\ns_{\delta}\delta'')-C_{{\nabla}\hspace{-1.0mm}\raisebox{0.3mm}{\text{\fontsize{6pt}{0cm}\selectfont{\bf /}}}_{\delta}\delta'}\delta''=\ns_{\delta'}(C_{\delta}\delta'')-C_{\delta}(\ns_{\delta'}\delta'')-C_{{\nabla}\hspace{-1.0mm}\raisebox{0.3mm}{\text{\fontsize{6pt}{0cm}\selectfont{\bf /}}}_{\delta'}\delta}\delta'' .
\end{equation*}
\item The unit element $e$ of the $\circ $-algebra is a $\ns$-flat holomorphic vector field: that is,
\begin{equation*}
\ns e=0.
\end{equation*} 
\item The bilinear form $\eta_\zeta$ and the product $\circ$ are homogeneous of degree $2-d$, $1$, respectively 
with respect to the Lie derivative ${\rm Lie}_{E}$ of the Euler vector field $E$: that is,
\begin{equation*}
{\rm Lie}_{E}(\eta_\zeta)=(2-d_{\widetilde f_n})\eta_\zeta,\quad {\rm Lie}_{E}(\circ)=\circ.
\end{equation*}
\end{enumerate}
\qed
\end{prop}
For the Frobenius manifold $M$, a connection $\widehat{\nabla}$ on $\T_{\PP^1_u\times M}$, called the {\em first structure connection}, is defined as follows:
\begin{subequations}
\begin{gather}
\widehat{\nabla}_{\delta'} \delta = \ns_{\delta'} \delta - \frac{1}{u} \delta \circ \delta',  \quad\delta, \delta' \in \T_{M},\\
\widehat{\nabla}_{u\frac{d}{du}} \delta  = \left( \frac{1}{u} C_E- \widetilde{\Q} \right) \delta, \quad \delta \in \T_{M}, \\
\widehat{\nabla}_{\delta} \frac{d}{du}   =  \widehat{\nabla}_{\frac{d}{du}} \frac{d}{du} = 0, \quad\delta \in \T_{M}, 
\end{gather}
\end{subequations}
where $\widetilde{\Q}\in{\rm End}_{\O_M}(\T_M)$ is defined by
\begin{equation}
\widetilde{\Q}(\delta)=\frac{2-d_{\widetilde f_n}}{2}\delta-\ns_{\delta} E,\quad\delta \in \T_{M} .
\end{equation}
\begin{rem}
Note that the parameter $z$ in \cite{D1,CDG} is $u^{-1}$ in this paper. 
\end{rem}
\begin{prop}[{\cite[Proposition~3.1]{D1}}]
The connection $\widehat{\nabla}$ is flat. 
\qed
\end{prop}
\begin{rem}
In fact, the primitive form $\zeta$ identifies the first structure connection $\widehat{\nabla}$ with the Gau\ss--Manin connection on the filtered de Rham cohomology (see \cite{ST}), which is flat. 
\end{rem}
Define an $\O_M$-endomorphism $\N\in{\rm End}_{\O_M}(\T_M)$ by 
\begin{equation}
\N(\delta):=\widetilde Q(\delta)+\frac{n}{2}\delta,\quad\delta\in\T_M.
\end{equation}
The eigenvalues of $\N$ are nothing but exponents of $\widetilde f_n$.
There exist local coordinates ${\bf t}=(t_\kappa)_{\kappa\in I_n}$, called {\em flat coordinates}, such that the unit vector field $e$ is given by $\p_{\psi({\bf 0})}$ and ${\rm Ker}\ns$ is spanned by $\p_{\kappa},~\kappa\in I_n$, where we denote $\p/\p t_\kappa$ by $\p_\kappa$. 
Moreover, one can choose flat coordinates ${\bf t}=(t_\kappa)_{\kappa\in I_n}$ satisfying 
\begin{equation}\label{eq : t and s}
t_\kappa({\bf 0})=0,\quad \frac{\p t_\kappa}{\p s_\lambda}({\bf 0})=\delta_{\kappa\lambda},\quad \kappa,\lambda\in I_n .
\end{equation}
In particular, we have 
\[
\eta_\zeta(\p_\kappa,\p_\lambda)=\left.\eta_\zeta\left(\frac{\p}{\p {s_\kappa}},\frac{\p}{\p {s_\lambda}}\right)\right|_{{\bf s}=0}=\eta^{(n)}_{\kappa\lambda},
\]
due to \eqref{eq : t and s}, and 
\[
\widetilde{\Q}(\p_\kappa)=\widetilde{Q}^{(n)}_{\kappa\kappa}\p_\kappa,
\]
since for each $\kappa\in I_n$ we have ${\rm Lie}_E\left(\p_\kappa\right)=(1-q_\kappa)\p_\kappa$ and 
\[
q_\kappa-\frac{d_{\widetilde f_n}}{2}=\sum_{i=1}^nk_i\widetilde\omega^{(n)}_{i}+\sum_{i=1}^n\widetilde\omega^{(n)}_{i}-\dfrac{n}{2}
=\sum_{i=1}^n\omega_{\kappa,i}^{(n)}-\frac{n}{2}=\widetilde Q^{(n)}_{\kappa\kappa},
\]
where $(k_1,\dots,k_n)=\psi^{-1}(\kappa)$. 
To summarize, we obtain the following
\begin{prop}
Matrix representations of $\widetilde{\Q}$ and $\eta_\zeta$ with respect to $\{\p_\kappa\}_{\kappa\in I_n}$ are given by $\widetilde{Q}^{(n)}$ and $\eta^{(n)}$, respectively. 
Moreover, the diagonal matrix $N^{(n)}=(n_\kappa)$ defined by 
\begin{equation}
n_\kappa:=q_\kappa+\frac{n-d_{\widetilde f_n}}{2}=\deg\,\zeta^{(n)}_\kappa, 
\end{equation}
is the matrix representing the $\O_M$-endomorphism $\N$.
\end{prop}
Let $\Z=\Z(u,{\bf t})$ be a function on an open subset in $\PP^1_u\times M$ and $(\eta^{\kappa\lambda})$ the inverse matrix of $\eta^{(n)}$. 
We say that $\Z$ is {\em $\widehat{\nabla}$-flat} if it satisfies 
\[
\widehat{\nabla} \biggl( \sum_{\kappa,\lambda\in I_n} \eta^{\kappa\lambda} u\frac{\partial \Z}{\partial t_{\lambda}} \frac{\p}{\p t_\kappa} \biggr) = 0,
\]
which is equivalent to the following differential system: 
for $\I= (\sum_{\lambda\in I_n} \eta^{\kappa\lambda} u\partial_{\lambda}\Z)$
\begin{subequations}\label{eqn : differential equation}
\begin{eqnarray}
\label{eqn : differential equation 1}& &\frac{\partial}{\partial t_{\kappa}} \I = -\frac{1}{u} C_{\p_\kappa} \I, \quad \kappa\in I_n, \\
\label{eqn : differential equation 2}& &u\frac{d}{du} \I= \left(\frac{1}{u} C_E - \widetilde{Q}^{(n)}\right) \I,
\end{eqnarray}
\end{subequations}
where we denote by the same letters $C_{\p_\kappa}$ and $C_E$ the matrix representations with respect to the basis $\{\p_\kappa\}_{\kappa\in I_n}$. 
The differential system \eqref{eqn : differential equation} can be considered as a family of meromorphic differential equations on $\PP^{1}_u$ parametrized by points on $M$, 
and the flatness of $\widehat\nabla$ means that this family is isomonodromic. 
The equation has a regular singular point at $u = \infty$ and an irregular singular point of Poincar\'e rank one at $u = 0$.
Let us consider a fundamental solution of \eqref{eqn : differential equation} at $u=\infty$. 
\begin{prop}[\cite{D1}]\label{prop : solution at infinity}
There exists a unique matrix fundamental solution $Y_\infty(u,{\bf t})=\Phi(u,{\bf t})u^{-\widetilde Q^{(n)}}$ such that 
$\Phi(u,{\bf t})=1+\Phi_1({\bf t})u^{-1}+\Phi_2({\bf t})u^{-2}+\cdots$ is a matrix-valued convergent power series in $u^{-1}$ satisfying
\begin{equation}
\Phi(u,{\bf 0})=1,\quad \Phi(-u,{\bf t})^T\eta^{(n)}\Phi(u,{\bf t})=\eta^{(n)}.
\end{equation}
In particular, the monodromy around infinity is given by 
\begin{equation}\label{eqn : monodromy at u = infty}
Y_{\infty}(e^{-2\pi\sqrt{-1}}u,{\bf t})=Y_{\infty}(u,{\bf t})\ee[\widetilde{Q}^{(n)}].
\end{equation}
\qed
\end{prop}
Denote by $\B\subset M$ the closed subset, called the {\em bifurcation set}, consisting of points where the values of canonical coordinates $(w_{1},\dots,w_{\widetilde\mu_n})$, namely, the 
critical values of $\widetilde F_n({\bf x};{\bf s})$, coincide.
It is known by \cite[Corollary~3.1]{D3} that the product $\circ$ on the tangent space $T_{\bf s}M$ for ${\bf s}\in M\setminus \B$ is semi-simple
and hence the Frobenius manifold $M$ is semi-simple.
Fix a point ${\bf s}\in M\setminus \B$ and define the matrix $\Theta=(\Theta_{i\lambda})$ of size $\widetilde \mu_{n}$ by 
\begin{equation}
\Theta_{i\lambda} :=\frac{\partial w_i}{\partial t_\lambda}\cdot 
\eta_\zeta\left(\frac{\partial}{\partial w_{i}},\frac{\partial}{\partial w_{i}}
\right)^{\frac{1}{2}},\quad \lambda\in I_n,~i=1,\dots,\widetilde\mu_n.
\end{equation}
We can construct a formal solution of \eqref{eqn : differential equation} near $u=0$ as follows.
\begin{prop}[{\cite[Lemma 4.3]{D3}}]
There exists a unique formal matrix fundamental solution $Y_{\rm formal}(u)$ of the differential equation \eqref{eqn : differential equation} in the form 
\begin{equation}
Y_{\rm formal}(u)=\Theta^{-1}G(u){\rm exp}(-W/{u}),
\end{equation}
where $W={\rm diag}(w_{1},\dots,w_{\widetilde\mu_n})$ and $G(u)=1+G_{1}u+G_{2}u^{2}+\cdots$ is a matrix-valued formal power series satisfying 
\[
G(-u)^T G(u) = 1.
\]
\qed
\end{prop}
Note that such $G(u)$ is a divergent power series in general since $u=0$ is an irregular singular point of Poincar\'e rank one. 
Denote by $\CC_u$ the complex plane whose coordinate is $u$. 
For $0 \le \phi < 1$, a line $\ell = e^{\pi\sqrt{-1}(\phi+1/2)}\RR\setminus\{0\}\in\CC_u$ is called {\it admissible} if 
\[
\frac{w_a-w_b}{u}\notin \RR,\quad u\in e^{\pi\sqrt{-1}\phi}\RR_{>0},\quad a,b=1,\dots,\widetilde\mu_n.
\]
For a small positive number $\varepsilon$, define sectors in $\CC_u$ by 
\begin{gather*}
\Pi_{\rm right}  =  \left\{u \in \CC_u\setminus\{0\}\,\left|\, 
\phi - \frac{1}{2} - \varepsilon < \frac{\arg{u}}{\pi} < \phi +\frac{1}{2}+ \varepsilon \right.\right\}, \\
\Pi_{\rm left}  =  \left\{u \in \CC_u\setminus\{0\}\,\left|\,
\phi +\frac{1}{2}- \varepsilon < \frac{\arg{u}}{\pi} < \phi + \frac{3}{2} + \varepsilon \right.\right\}. 
\end{gather*} 
\begin{figure}[h]
\centering
\includegraphics[pagebox=cropbox]{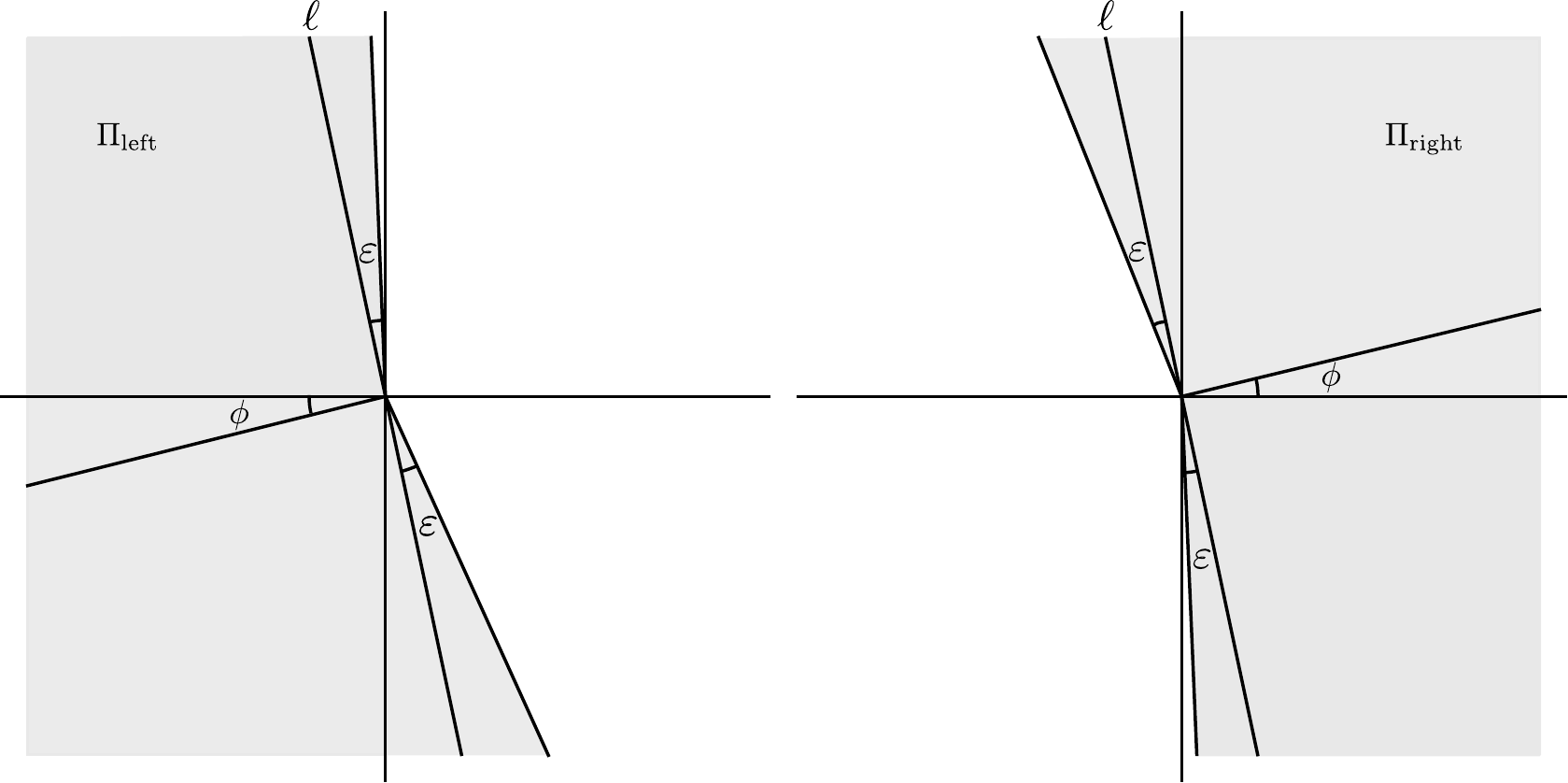}
\caption{Left sector $\Pi_{\rm left}$ and right sector $\Pi_{\rm right}$}
\end{figure}

It is known by the general theory for ordinary differential equations (cf. \cite[Theorem~A]{BJL}) that 
there exist unique solutions $Y_{\rm right/left}$ of \eqref{eqn : differential equation 2} analytic in $u$ in the sectors $\Pi_{\rm right/left}$ 
satisfying the following asymptotic properties:
\begin{gather*}
Y_{\rm right}(u)  \sim  Y_{\rm formal}(u)\quad \text{as} \hspace{+.5em} u \rightarrow 0,  \hspace{+.3em} u \in \Pi_{\rm right}, \\
Y_{\rm left}(u)  \sim  Y_{\rm formal}(u)\quad
\text{as} \hspace{+.5em} u \rightarrow 0, \hspace{+.3em} u \in \Pi_{\rm left}.
\end{gather*}
\begin{figure}[h]
\centering
\includegraphics[pagebox=cropbox]{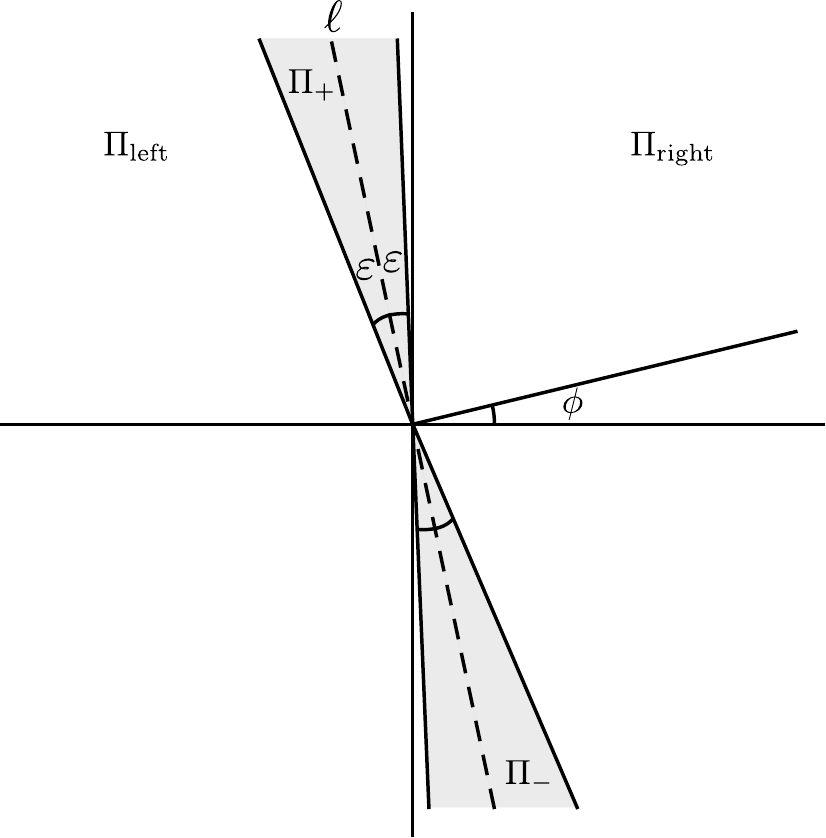}
\caption{$\Pi_{+}$ and $\Pi_{-}$}
\end{figure}
In the sectors $\Pi_+$ and $\Pi_-$ in $\Pi_{\rm right} \cap \Pi_{\rm left}$ defined by
\begin{gather*}
\Pi_+ := \left\{u \in \CC_u\setminus\{0\}\,\left|\, 
\phi + \frac{1}{2} - \varepsilon < \frac{\arg{u}}{\pi} < \phi +\frac{1}{2}+ \varepsilon \right.\right\}, \\
\Pi_- := \left\{u \in \CC_u\setminus\{0\}\,\left|\, 
\phi + \frac{3}{2} - \varepsilon < \frac{\arg{u}}{\pi} < \phi +\frac{3}{2}+ \varepsilon \right.\right\}, 
\end{gather*}
we have two analytic solutions $Y_{\rm right}(u)$ and $Y_{\rm left}(u)$.
They are related by
\begin{gather}
Y_{\rm right}(u) = Y_{\rm left}(u) S^\phi, \quad u \in \Pi_{+},\\
Y_{\rm right}(e^{-2\pi\sqrt{-1}}u) = Y_{\rm left}(u) (S^\phi)^T, \quad u \in \Pi_{-},
\end{gather}
with a matrix $S^\phi$ independent of $u$. 
\begin{defn}
The matrix $S^\phi$ is called the {\it Stokes matrix} of the 
first structure connection of $M$ 
(for the admissible line $\ell$ and the chosen point on $M\setminus\B$).
\end{defn}
The fundamental solution $Y_\infty(u):=Y_\infty(u,{\bf s})$ around $u=\infty$ is related with the analytic solution $Y_{\rm right}(u)$ by 
\begin{equation}\label{eqn : central connection matrix}
Y_{\rm right}(u)=Y_{\infty}(u)C^\phi,\quad u\in\Pi_{\rm right},~|u|\gg0 ,
\end{equation}
with a matrix $C^\phi$ independent of $u$. 
The matrix $C^\phi$ is called the {\em central connection matrix} (for the admissible line $\ell$ and the chosen point ${\bf s}\in M\setminus\B$). 
The following proposition is a consequence of the isomonodromic property 
of the differential equation \eqref{eqn : differential equation 2}.
\begin{prop}[{\cite[Isomonodromicity Theorem]{D3}}]
The Stokes matrix $S^\phi$ and the central connection matrix $C^\phi$ are locally constant as a function on $M\setminus\B$. 
\qed
\end{prop}
\begin{prop}[{\cite{D3}}]\label{prop : D3}
We have the following relations; 
\begin{subequations}
\begin{equation}
(C^\phi)^{-1}\ee[\widetilde{Q}^{(n)}]C^\phi=(S^\phi)^{-1} (S^\phi)^T, 
\end{equation}
\begin{equation}
(C^\phi)^T\ee\left[\frac{1}{2}\widetilde{Q}^{(n)}\right]\eta^{(n)} C^\phi=S^\phi.
\end{equation}
\end{subequations}
\qed
\end{prop}
Consider an unfolding 
\begin{equation}
\widetilde f_{n;s}({\bf x}):=\widetilde f_n({\bf x})+s\cdot x_n ,
\end{equation}
with one parameter $s\in\CC$ with $|s| \ll1$. 
Since we have $0<\deg\,s<1$ where $[\frac{\p}{\p s},E]=(\deg\,s)\frac{\p}{\p s}$, 
the restriction of the domain of $\widetilde f_{n;s}$ to a small neighborhood of the origin in $\CC^n$ 
(such as \eqref{eq: 5.3}) is not necessary, that is, the unfolding $\widetilde f_{n;s}$ is defined globally on $\CC^n$ for each $|s| \ll1$. 
Moreover, we have the following
\begin{prop}[{cf.~\cite[Appendix A]{V}}]\label{prop : small deformation}
For each $s$, the point on the base space $M$ of the universal unfolding $\widetilde F_n$ corresponding to $\widetilde f_{n;s}$ belongs to $M\setminus \B$. 
\qed
\end{prop}
In what follows, we determine the monodromy data at the point in Proposition \ref{prop : small deformation}.
Let $\phi=0$ if the line $\sqrt{-1}\RR\setminus\{0\}$ is admissible, if not, be a sufficiently small positive number so that all lines of the form $e^{\pi\sqrt{-1}(\phi'+1/2)}\RR\setminus\{0\}$, $0<\phi'\le2\phi$ are admissible.
Let ${\bf p}_{1}, \dots, {\bf p}_{\widetilde\mu_{n}}$ be the critical points of the holomorphic function $\widetilde f_{n;s}({\bf x}): \CC^n\longrightarrow\CC_w$ 
and $w_{j}:=\widetilde f_{n;s}({\bf p}_j)$. 
For $u\in e^{\pi\sqrt{-1}\phi}\RR_{>0}$ and each critical point ${\bf p}_{j}$, we can define a relative $n$-cycle 
\begin{equation}\label{eqn : Lefschetz thimble}
\Gamma_j(u) \in H_n(\CC^n,{\rm Re}(\widetilde f_{n;s}({\bf x})/u)\gg 0;\ZZ),
\end{equation} 
called the {\em Lefschetz thimble} for ${\bf p}_{j}$ as follows. 
The image of $\Gamma_{j}(u)$ by $\widetilde f_{n;s}({\bf x})$ is a half-line $w_{j}+e^{\pi\sqrt{-1}\phi}\RR_{\ge 0}\subset \CC_w$ 
and the fiber over a point $w$ on this half-line is the $(n-1)$-cycle in $Y_{w;s}:=\{\widetilde f_{n;s}({\bf x})=w\}\subset\CC^n$ which shrinks to the critical point ${\bf p}_{j}$ 
by the parallel transport along this half-line. 
Choose $w_0\in\CC_w$ so that ${\rm Re}(w_0/u)\gg 0$ and we have 
\begin{equation}\label{eq: 5.22}
H_n(\CC^n,{\rm Re}(\widetilde f_{n;s}({\bf x})/u)\gg 0;\ZZ)\cong H_n(\CC^n,Y_{w_0;s};\ZZ).
\end{equation}
Take paths $\gamma_{j}:[0,1]\longrightarrow \CC_w$ with $\gamma_{j}(0)=w_{j}$ and $\gamma_{j}(1)=w_0$ for $j = 1,\dots,\widetilde \mu_{n}$ satisfying the following conditions:
\begin{itemize}
\item
For each $j=1,\dots, \widetilde\mu_n$, the path $\gamma_j$ has no self-intersection.
\item 
For $i \ne j$, $\gamma_{i}$ and $\gamma_{j}$ have a unique common point $w_0$.
\item 
The paths $\gamma_{1},\dots,\gamma_{\widetilde \mu_{n}}$ are ordered counter-clockwise, namely, $\arg \gamma_i'(1)<\arg \gamma_j'(1)$ if $i<j$.
\end{itemize}
We obtain an $(n-1)$-sphere in $Y_{w_0;s}$ which shrinks to a point at $w_{j}$ by the parallel transport along the path $\gamma_{j}$, whose 
homology class $L_j\in H_{n-1}(Y_{w_0;s};\ZZ)$ is called a {\em vanishing cycle} along $\gamma_{j}$. 
The ordered set $(L_{1},\dots,L_{\widetilde \mu_{n}})$ of vanishing cycles is called a {\em (counter-clockwise) distinguished basis of vanishing cycles}.
\begin{rem}
Note that our ordering of vanishing cycles is inverse to the usual one. 
Since we only use this counter-clockwise ordering of vanishing cycles in this paper, we just call it a distinguished basis for simplicity.
\end{rem}
\begin{figure}[h]
\centering
\includegraphics[pagebox=cropbox]{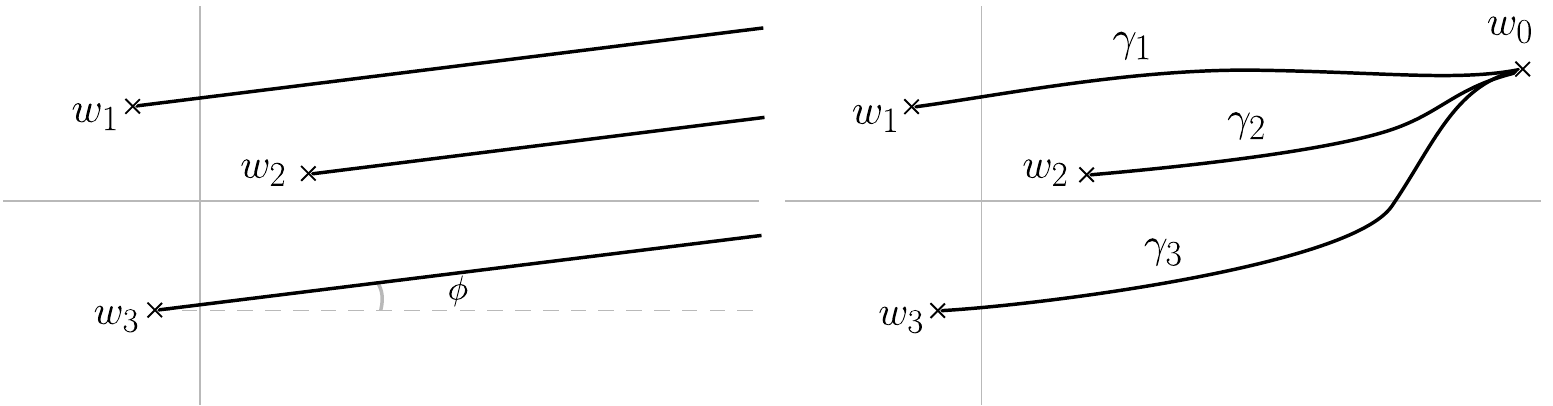}
\caption{Left: The image of $\Gamma_j(u=e^{\pi\sqrt{-1}\phi})$ by $\widetilde f_{n;s}$,\quad Right: Paths $\gamma_j$.}
\end{figure}
Since we have the natural isomorphism 
\begin{equation}
H_n(\CC^n,Y_{w_0;s};\ZZ)\stackrel{\cong}{\longrightarrow} H_{n-1}(Y_{w_0;s};\ZZ),
\end{equation}
it turns out that homology class represented by the Lefschetz thimble $\Gamma_{j}(u)$ for ${\bf p}_{j}$ are uniquely determined by the vanishing cycle $L_{j}$. 
The primitive form $\zeta$ is given by $[d{\bf x}]$ over the subset of $M$ 
where $\widetilde F_n=\widetilde f_{n;s}$ since $0<\deg\,s<1$ and the normalization \eqref{eqn : primitive form}. Therefore, it follows from the definition of the primitive form that 
\[
\I_{j}:=(\I_{\kappa j}),\quad 
\I_{\kappa j}:=\frac{1}{(2\pi u)^{\frac{n}{2}}}\sum_{\lambda\in I_n} \eta^{\lambda\kappa}
\int_{\Gamma_{j}(u)}e^{-\frac{\widetilde f_{n;s}({\bf x})}{u}}\zeta^{(n)}_\lambda,
\]
satisfies
\begin{subequations}
\begin{eqnarray}
& &\frac{\partial}{\partial s} \I_{j} = -\frac{1}{u} C_{\frac{\p}{\p s}} \I_{j}, \\
& &u\frac{d}{du} \I_{ j}= \left(\frac{1}{u} C_{(\deg\,s)s\frac{\p}{\p s}} - \widetilde Q^{(n)}\right) \I_{j}.
\end{eqnarray}
\end{subequations} 
where $\zeta^{(n)}_\lambda$ is the class in $\HH^n(\Omega^\bullet_{\CC^n}, d-d\widetilde f_n\wedge)$ corresponding to the class $\zeta^{(n)}_\lambda\in \Omega_{\widetilde f_n}$ under the isomorphism $\HH^n(\Omega^\bullet_{\CC^n}, d-d\widetilde f_n\wedge)\cong\Omega_{\widetilde f_n}$.
By the saddle-point method, the oscillatory integral
\begin{equation*}
\Z_{j}(u) :=\frac{1}{(2\pi u)^{\frac{n}{2}}}
\int_{\Gamma_{j}(u)}e^{-\frac{\widetilde f_{n;s}({\bf x})}{u}}d{\bf x},\quad j = 1,\dots,\widetilde \mu_{n}
\end{equation*}
has an asymptotic expansion
\begin{equation*}
\Z_{j}(u) = 
\frac{e^{-\frac{w_{j}}{u}}}
{\sqrt{\Delta_{j}}}(1 + O(u)),\quad u \rightarrow 0,\quad j=1,\dots,\widetilde\mu_n,
\end{equation*}
where $\Delta_{j}$ is the Hessian at ${\bf p}_{j}$:
\begin{equation*}
\Delta_{j} := {\rm det}\left( \frac{\partial^{2} \widetilde f_{n;s}({\bf x})}
{\partial x_{k} \partial x_{l}}({\bf p}_{j})\right)_{k,l=1,\dots,n},\quad j=1,\dots,\widetilde\mu_{n}.
\end{equation*}
Since $(w_1,\dots, w_n)$ are canonical coordinates, we have
\[
\frac{\p \widetilde F_n}{\p w_i}({\bf p}_j)=\delta_{ij},\quad 
\eta_\zeta\left(\frac{\p }{\p w_{i}},\frac{\p }{\p w_{j}}\right)=
\frac{\delta_{ij}}{\Delta_{i}},\quad i,j=1,\dots,\widetilde\mu_{n}.
\]
Hence, it turns out that the asymptotic expansion of  $\I(u):=(\I_{\kappa j}(u))$ coincides with 
$Y_{{\rm formal}}(u)$ if a suitable branch of the square root is chosen. 
The Lefschetz thimble $\Gamma_{j}(u)$ discontinuously changes when $u$ cross a half-line of the form $e^{\sqrt{-1}\arg (w_j-w_i)}\RR_{>0}\subset \CC_u$.
The discontinuity causes the Stokes phenomena for the oscillatory integrals. 
In order to determine the Stokes matrix $S^\phi$ of the Frobenius manifold $M$ for the admissible line $\ell$, 
we establish the correspondence between the analytic solutions $Y_{\rm right/left}$ and $\I$.
Order critical values $\{w_{j}\}_{j=1}^{\widetilde \mu_{n}}$ so that 
\[
|e^{-{w_{1}}/{u}}| \ll |e^{-{w_{2}}/{u}}| \ll \cdots \ll |e^{-{w_{\widetilde \mu_{n}}}/{u}}|
\]
holds as $u \rightarrow 0$ along the line $e^{\pi\sqrt{-1}(\phi+1/2)}\RR_{>0}$.
Consider the local system on $\CC_u\setminus\{0\}$ whose fiber over $u \in \CC_u\setminus\{0\}$ is the relative homology group 
$H_n(\CC^n,{\rm Re}(\widetilde f_{n;s}({\bf x})/u) \gg 0;\ZZ)$, and let $\Gamma_{{\rm right/left}, j}(u)$ be a section of the local system on 
$\Pi_{\rm right/left}$ satisfying the following condition; 
\begin{itemize}
\item[(A)] for $u \in \{ u \in \CC_u\setminus\{0\}~|~0<\dfrac{\arg{u}}{\pi}<2\phi \} 
\subset \Pi_{\rm right}$ (resp., $u \in \{ u \in \CC_u\setminus\{0\}~|~
1<\dfrac{\arg{u}}{\pi}<2\phi+1\} \subset \Pi_{\rm left}$),
$\Gamma_{{\rm right},j}(u)$ (resp., $\Gamma_{{\rm left},j}(u)$) 
coincides with the relative homology class represented by 
the Lefschetz thimble $\Gamma_j(u)$ for the $j$-th critical point ${\bf p}_{j}$. 
\end{itemize}
Since $\Pi_{\rm right/left}$ is simply-connected, the condition (A) determines the section $\Gamma_{{\rm right/left},j}(u)$ uniquely. 
Let 
\begin{subequations}
\begin{eqnarray}
\label{eqn : definition of Y right}& & Y_{\rm right}(u)_{\kappa j} :=
\frac{1}{(2\pi u)^{\frac{n}{2}}}\sum_{\lambda\in I_n}\eta^{\lambda\kappa}
\int_{\Gamma_{{\rm right}, j}(u)}e^{-\frac{\widetilde f_{n;s}({\bf x})}{u}}\zeta^{(n)}_\lambda,\\ 
& & Y_{\rm left}(u)_{\kappa j} :=
\frac{1}{(2\pi u)^{\frac{n}{2}}}\sum_{\lambda\in I_n}\eta^{\lambda\kappa}
\int_{\Gamma_{{\rm left}, j}(u)}e^{-\frac{\widetilde f_{n;s}({\bf x})}{u}}\zeta^{(n)}_\lambda, 
\end{eqnarray}
\end{subequations}
for $\kappa\in I_n,~j=1,\dots,\widetilde\mu_n$. 
The matrix valued function $Y_{\rm right/left}(u)=(Y_{\rm right/left}(u)_{\kappa j} )$ in $u$ is a fundamental solution of \eqref{eqn : differential equation} which is asymptotic to 
$Y_{\rm formal}(u)$ as $u \rightarrow 0$ in the whole sector $\Pi_{\rm right/left}$ by the condition (A). 
Therefore, the Stokes matrix $S^\phi$ can be read off 
from the monodromy of integration cycles since the integrand is single-valued. 
That is, for $u \in \Pi_{+}$, $S^\phi = (S^\phi_{ij})_{i,j=1,\cdots,\widetilde{\mu}_{n}}$ satisfies 
\begin{equation}
\Gamma_{{\rm right},j}(u) = \sum_{i = 1}^{\widetilde\mu_{n}} S^\phi_{ij} \Gamma_{{\rm left},i}(u), \quad j = 1,\cdots,\widetilde \mu_n.
\end{equation} 
By the Picard-Lefschetz formula, it turns out that the Stokes matrix coincides with the Seifert matrix.
\begin{prop}\label{prop : Stokes and intersection form}
Let the notations be as above.
Let ${\mathbb S}$ be the Seifert form on $H_{n-1}(Y_{w_0;s};\ZZ)$ whose matrix representation is given by 
\begin{equation}\label{eq: 5.27}
{\mathbb S}(L_{i},L_{j}):=
\begin{cases}
0 \quad &\text{if}\quad  i>j,\\
1 \quad &\text{if}\quad  i=j,\\
(-1)^{\frac{(n-1)(n-2)}{2}}I(L_{i},L_{j}) \quad &\text{if}\quad  i<j,
\end{cases}
\end{equation}
where $I$ is the intersection form on $H_{n-1}(Y_{w_0;s};\ZZ)$. We have 
\[
S^\phi_{ij}={\mathbb S}(L_{i},L_{j}),\quad i,j=1,\cdots,\widetilde{\mu}_{n}.
\]
\end{prop}
\begin{pf}
For $u\in\Pi_+$, choose $w_0\in\CC_w$ so that ${\rm Re}(w_0/u)\gg 0$ and the isomorphism~\eqref{eq: 5.22} holds. 
For each $j$, let $L_{{\rm right/left},j}\in  H_{n-1}(Y_{w_0;s};\ZZ)$ be the vanishing cycle corresponding to the Lefschetz thimble $\Gamma_{{\rm right/left},j}(u)$. 
Note that $L_j=L_{{\rm right},j}$. 
By the ordering of $\gamma_1,\dots,\gamma_{\widetilde{\mu}_n}$, we have $L_{{\rm left},j}=h^{-1}_{L_{{\rm right},1}}\circ\cdots\circ h^{-1}_{L_{{\rm right},j-1}}(L_{{\rm right},j})$, where $h_{L_{{\rm right},k}}$ is the Picard-Lefschetz transformation associated to the vanishing cycle $L_{{\rm right},k}$ 
(see Figure \ref{fig : PL transform}). 
By the Picard-Lefschetz formula
\[
h^{-1}_{L_{{\rm right},i}}(L)=L-(-1)^{\frac{(n-1)(n-2)}{2}}I(L_{{\rm right},i},L)L_{{\rm right},i},
\]
we obtain the statement.
\begin{figure}[h]
\begin{tabular}{cc}
\begin{minipage}{0.5\hsize}
\centering
\includegraphics[pagebox=cropbox]{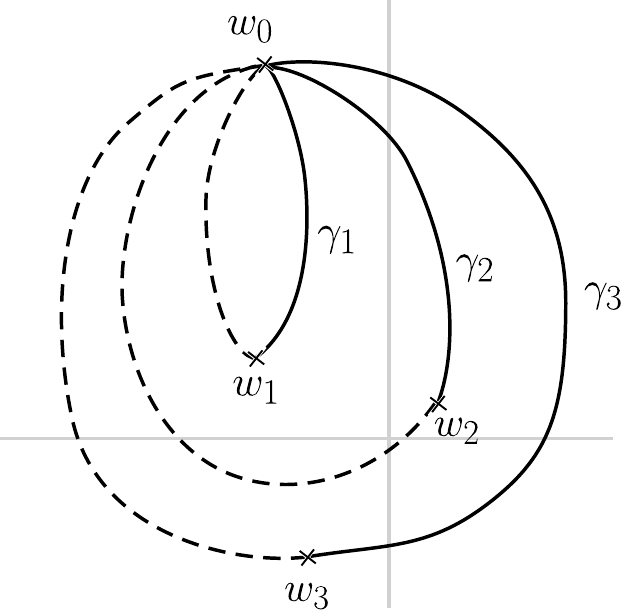}
\caption{}
\end{minipage}
\begin{minipage}{0.5\hsize}
\centering
\includegraphics[pagebox=cropbox]{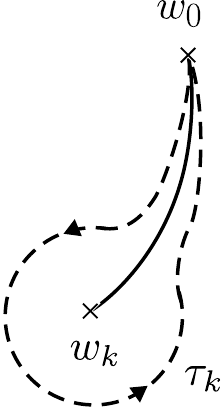}
\caption{$h_{L_{{\rm right},k}}$ is the monodromy operator along $\tau_k$.}
\label{fig : PL transform}
\end{minipage}
\end{tabular}
\end{figure}
\qed
\end{pf}
\begin{prop}\label{prop : central connection and oscillatory integral}
For each $j=1,\dots,\widetilde\mu_n$, let $\Gamma^0_j(u):=\lim_{s\to 0}\Gamma_{{\rm right},j}(u)$ be the element in $H_n(\CC^n, {\rm Re}(\widetilde f_n/u)\gg 0;\ZZ)$. 
The central connection matrix $C^\phi$ is given by 
\begin{equation}
C^\phi=\left(C^\phi_{\kappa j}\right),\quad C^\phi_{\kappa j}=\frac{1}{(2\pi)^{\frac{n}{2}}}\sum_{\lambda\in I_n}\eta^{\lambda\kappa}
\int_{\Gamma^0_j(1)}e^{-\widetilde f_n}\zeta^{(n)}_\lambda .
\end{equation}
\end{prop}
\begin{pf}
Obviously, we have $Y_\infty(u)|_{{\bf s}=0}=u^{-\widetilde{Q}^{(n)}}$. Therefore, we have 
\[
u^{-\widetilde{Q}^{(n)}}\cdot C^\phi=Y_\infty(u)|_{{\bf s}=0}\cdot C^\phi=Y_{\rm right}(u)|_{{\bf s}=0}
=\left(\frac{1}{(2\pi u)^{\frac{n}{2}}}\sum_{\lambda\in I_n}\eta^{\lambda\kappa}
\int_{\Gamma^0_{j}(u)}e^{-\frac{\widetilde f_{n}({\bf x})}{u}}\zeta^{(n)}_\lambda\right).
\]
Since $\widetilde f_{n}$ is weighted homogeneous, we have
\[
u^{-\frac{n}{2}}\int_{\Gamma^0_{j}(u)}e^{-\frac{\widetilde f_{n}({\bf x})}{u}}\zeta^{(n)}_\lambda=
u^{\widetilde Q^{(n)}_{\lambda\lambda}}\int_{\Gamma^0_{j}(1)}e^{-\widetilde f_{n}({\bf x})}\zeta^{(n)}_\lambda,
\]
and hence the statement follows.
\qed 
\end{pf}
\begin{cor}
Let $\{L^{(n)}_i\}_{i=1}^{\widetilde\mu_n}$ be the cycles constructed in Section \ref{sec : proof of integral structure}.
We have 
\[
{\mathbb S}(L^{(n)}_i,L^{(n)}_j)=\chi^{(n)}_{ij},\quad i,j=1,\dots,\widetilde\mu_n.
\]
\end{cor}
\begin{pf}
Proposition \ref{prop : Stokes and intersection form} implies that ${\mathbb S}(L^{(n)}_i,L^{(n)}_j)=(U^T S^0 U)_{ij}$ at the limit $\phi\to0$, where $U=(U_{ij})$ is the change of basis matrix with respect to $\{L_1^{(n)},\dots,L_{\widetilde\mu_n}^{(n)}\}$ and $\{L_1,\dots,L_{\widetilde\mu_n}\}$.
By Theorem \ref{thm : Integral structure} and Proposition \ref{prop : central connection and oscillatory integral}, we have $C^0U=(2\pi)^{-\frac{n}{2}}{\rm ch}^{(n)}_\Gamma$.
Hence, it follows from Proposition \ref{prop : D3} and Theorem \ref{thm : main theorem} that 
\begin{align*}
U^T S^0 U & = U^T (C^0)^T \ee\left[\frac{1}{2}\widetilde{Q}^{(n)}\right]\eta^{(n)} C^0 U \\
& = \left(\dfrac{1}{(2\pi)^{\frac{n}{2}}}\ch_{\Gamma}^{(n)}\right)^T\ee\left[\dfrac{1}{2}\widetilde Q^{(n)}\right]\eta^{(n)}\left(\dfrac{1}{(2\pi)^{\frac{n}{2}}}\ch_{\Gamma}^{(n)}\right) \\
& = \chi^{(n)}.
\end{align*}
\qed
\end{pf}
Therefore, we obtain the commutativity of the following diagram:
\[
\xymatrix{
(K_0({\rm HMF}^{L_{f_n}}_{S_n}(f_n)),\chi) \ar[rr]^{\cong}  \ar[d]_{(2\pi\sqrt{-1})^{-n}{\rm ch}^{(n)}_\Gamma} & & (H_{n-1}(\widetilde f_n^{-1}(1);\ZZ),{\mathbb S}) \ar[d]^{{\mathbb D}} \\
(\Omega_{f_n,G_{f_n}}, (2\pi\sqrt{-1})^n {\mathbb S}_{f_n,G_{f_n}}) \ar[rr]^{\bf mir} & & (\Omega_{\widetilde f_n},(2\pi\sqrt{-1})^n {\mathbb S}_{\widetilde f_n})
},
\]
where ${\mathbb S}$ denotes the Seifert form, 
and the $\CC$-bilinear forms ${\mathbb S}_{f_n,G_{f_n}}$ and ${\mathbb S}_{\widetilde f_n}$ are defined by
\[
{\mathbb S}_{f_n,G_{f_n}}(\xi_1,\xi_2)=J_{f_n,G_{f_n}}\left({\bf e}\left[\frac{1}{2}N^{(n)}\right]\xi_1,\xi_2\right),\quad \xi_1,\xi_2\in\Omega_{f_n,G_{f_n}} ,
\]
\[
{\mathbb S}_{\widetilde f_n}(\zeta_1,\zeta_2)=J_{\widetilde f_n}\left({\bf e}\left[\frac{1}{2}\N\right]\zeta_1,\zeta_2\right),\quad \zeta_1,\zeta_2\in\Omega_{\widetilde f_n} ,
\]
$N^{(n)}$ and $\N$ are regarded as elements of ${\rm End}_\CC(\Omega_{f_n,G_{f_n}})$ and ${\rm End}_\CC(\Omega_{\widetilde f_n})$, respectively.
Since the intersection form is the symmetrized Seifert form, we have $I=S^\phi+(S^\phi)^T$ by Proposition \ref{prop : Stokes and intersection form}.
Hence, we obtain the formula between the intersection form and periods introduced by K. Saito. 
\begin{prop}
For each $w\in \RR_{>0}$, we have the Saito's formula of the intersection form (see \cite[Section~3.4, Theorem]{S-K}) :
\[
I(L_i,L_j)=\frac{(-1)^{\frac{n(n-1)}{2}}}{(2\pi\sqrt{-1})^{n-1}}\sum_{\lambda,\kappa\in I_n}\oint_{L_i(w)}\frac{\zeta^{(n)}_\lambda}{d\widetilde f_n}
\cdot \eta^{\lambda\kappa}\cdot w\p_w^{n-1}\oint_{L_j(w)}\frac{\zeta^{(n)}_\kappa}{d\widetilde f_n},\quad i,j=1,\dots ,\widetilde \mu_n,
\]
where $L_i(w)$, $i=1,\dots, \widetilde \mu_n$ is the image of $L_i$ by the isomorphism $H_{n-1}(\widetilde{f}_n^{-1}(1);\ZZ)\cong H_{n-1}(\widetilde{f}_n^{-1}(w);\ZZ)$ induced by the parallel transformation and $\frac{\zeta^{(n)}_\lambda}{d\widetilde f_n}$ denotes the Gelfand--Leray form defined by
\[
\zeta^{(n)}_\lambda=d\widetilde f_n\wedge \frac{\zeta^{(n)}_\lambda}{d\widetilde f_n}.
\]
\end{prop}
\begin{pf}
By the definition, the intersection matrix $I$ is given by
\begin{eqnarray*}
I&=&(-1)^{\frac{(n-1)(n-2)}{2}}\left(S^\phi+(-1)^{n-1}(S^\phi)^T\right)\\
&=&(-1)^{\frac{(n-1)(n-2)}{2}}\left((C^\phi)^T\ee\left[\frac{1}{2}\widetilde Q^{(n)}\right]\eta^{(n)} C^\phi+(-1)^{n-1}(C^\phi)^T\eta^{(n)}\ee\left[\frac{1}{2}\widetilde Q^{(n)}\right]C^\phi\right)\\
&=&(-1)^{\frac{n(n+1)}{2}}(\sqrt{-1})^nC^T\eta^{(n)}\left(\ee\left[\frac{1}{2} N^{(n)}\right]-\ee\left[-\frac{1}{2} N^{(n)}\right]\right)C^\phi\\
&=&(-1)^{\frac{(n-1)(n-2)}{2}}(\sqrt{-1})^{n-1}\cdot 2\pi (C^\phi)^T\eta^{(n)}\frac{\sin\pi N^{(n)}}{\pi}C^\phi.
\end{eqnarray*}
Therefore, by the Euler's reflection formula and the inverse Laplace transforms, we have
\begin{eqnarray*}
I(L_i,L_j)&=&\frac{(-1)^{\frac{n(n-1)}{2}}}{(2\pi\sqrt{-1})^{n-1}}\sum_{\lambda,\kappa\in I_n}\int_{\Gamma^0_i(1)}e^{-\widetilde f_n}\zeta^{(n)}_\lambda\cdot
\frac{\eta^{\lambda\kappa}}{\Gamma(n_\lambda)\Gamma(1-n_\lambda)}\cdot \int_{\Gamma^0_j(1)}e^{-\widetilde f_n}\zeta^{(n)}_\kappa\\
&=&\frac{(-1)^{\frac{n(n-1)}{2}}}{(2\pi\sqrt{-1})^{n-1}}\sum_{\lambda,\kappa\in I_n}\frac{w^{n_\lambda-1}}{\Gamma(n_\lambda)}\int_{\Gamma^0_i(1)}e^{-\widetilde f_n}\zeta^{(n)}_\lambda\cdot
\eta^{\lambda\kappa}\cdot \frac{w^{n_\kappa-n+1}}{\Gamma(n_\kappa-n+1)}\int_{\Gamma^0_j(1)}e^{-\widetilde f_n}\zeta^{(n)}_\kappa\\
&=&\frac{(-1)^{\frac{n(n-1)}{2}}}{(2\pi\sqrt{-1})^{n-1}}
\sum_{\lambda,\kappa\in I_n} \frac{1}{2\pi\sqrt{-1}}\int_{c-\sqrt{-1}\infty}^{c+\sqrt{-1}\infty}e^{\frac{w}{u}}u^{n_\lambda}d\left(\frac{1}{u}\right)\cdot \int_{\Gamma^0_i(1)}e^{-\widetilde f_n}\zeta^{(n)}_\lambda\\
& &\cdot \eta^{\lambda\kappa}\cdot\frac{1}{2\pi\sqrt{-1}}w\p_w^{n-1}\int_{c-\sqrt{-1}\infty}^{c+\sqrt{-1}\infty}e^{\frac{w}{u}}u^{n_\kappa}d\left(\frac{1}{u}\right)\cdot \int_{\Gamma^0_j(1)}e^{-\widetilde f_n}\zeta^{(n)}_\kappa\\
&=&\frac{(-1)^{\frac{n(n-1)}{2}}}{(2\pi\sqrt{-1})^{n-1}}\sum_{\lambda,\kappa\in I_n}
\frac{1}{2\pi\sqrt{-1}}\int_{c-\sqrt{-1}\infty}^{c+\sqrt{-1}\infty}e^{\frac{w}{u}}\left(\int_{\Gamma^0_i(u)}e^{-\frac{\widetilde f_n}{u}}\zeta^{(n)}_\lambda\right) d\left(\frac{1}{u}\right) \\
& &\cdot \eta^{\lambda\kappa}\cdot\frac{1}{2\pi\sqrt{-1}}w\p_w^{n-1}\int_{c-\sqrt{-1}\infty}^{c+\sqrt{-1}\infty}e^{\frac{w}{u}}\left(\int_{\Gamma^0_j(u)}e^{-\frac{\widetilde f_n}{u}}\zeta^{(n)}_\kappa\right) d\left(\frac{1}{u}\right)\\
&=&\frac{(-1)^{\frac{n(n-1)}{2}}}{(2\pi\sqrt{-1})^{n-1}}\sum_{\lambda,\kappa\in I_n}\oint_{L_i(w)}\frac{\zeta^{(n)}_\lambda}{d\widetilde f_n}
\cdot \eta^{\lambda\kappa}\cdot w\p_w^{n-1}\oint_{L_j(w)}\frac{\zeta^{(n)}_\kappa}{d\widetilde f_n}.\\
\end{eqnarray*}
Here, since $\widetilde f_n$ is weighted homogeneous, note also that
\[
\int_{\Gamma^0_i(u)}e^{-\frac{\widetilde f_n}{u}}\zeta^{(n)}_\lambda=u^{n_\lambda}\int_{\Gamma^0_i(1)}e^{-\widetilde f_n}\zeta^{(n)}_\lambda
=\int_0^\infty e^{-\frac{w}{u}}\left(\oint_{L_i(w)}\frac{\zeta^{(n)}_\lambda}{d\widetilde f_n}\right)dw
\]
for ${\rm Re}(u)>0$.
\qed
\end{pf}
A distinguished basis can be chosen in other ways.
The Artin's {\it braid group} $B_{\widetilde\mu_n}$ on $\widetilde\mu_n$-stands is a group presented by the following generators and relations: 
\begin{description}
\item[{\bf Generators}] $\{b_i~|~i=1,\dots, \widetilde\mu_n-1\}$
\item[{\bf Relations}] $b_{i}b_{j}=b_{j}b_{i}$ for $|i-j|\ge 2$, $b_{i}b_{i+1}b_{i}=b_{i+1}b_{i}b_{i+1}$ for $i=1,\dots, \widetilde\mu_n-2$.
\end{description}
Consider the group $B_{\widetilde\mu_n}\ltimes (\ZZ/2\ZZ)^{\widetilde\mu_n}$, the semi-direct product of the braid group $B_{\widetilde\mu_n}$ and the 
abelian group $(\ZZ/2\ZZ)^{\widetilde\mu_n}$, defined by the group homomorphism $B_{\widetilde\mu_n}\rightarrow {\mathfrak S}_{\widetilde\mu_n} \rightarrow {\rm Aut}_{\ZZ} (\ZZ/2\ZZ)^{\widetilde\mu_n}$,
where the first homomorphism is $b_{i}\mapsto (i, i+1)$ and the second one is
induced by the natural actions of the symmetric group ${\mathfrak S}_{\widetilde\mu_n}$ on $(\ZZ/2\ZZ)^{\widetilde\mu_n}$.  
The group $B_{\widetilde\mu_n}\ltimes (\ZZ/2\ZZ)^{\widetilde\mu_n}$ acts on the set of distinguished basis of vanishing cycles by 
\begin{gather*}
b_{i}(L_{1},\dots, L_{\widetilde\mu_n}):=(L_{1}, \dots, L_{i-1}, L_{i+1}, h_{L_{i+1}}(L_i),L_{i+2},\dots, L_{\widetilde\mu_n}),\\
b^{-1}_{i}(L_{1},\dots, L_{\widetilde\mu_n}):=(L_{1}, \dots, L_{i-1}, h^{-1}_{L_{i}}(L_{i+1}), L_{i}, L_{i+2},\dots, L_{\widetilde\mu_n}),\\
p_{i}(L_{1},\dots, L_{\widetilde\mu_n}):=(L_{1}, \dots, L_{i-1}, -L_{i}, L_{i+1},\dots, L_{\widetilde\mu_n}),
\end{gather*}
where 
\begin{gather*}
h_{L_{i+1}}(L_i)=L_{i}-(-1)^{\frac{(n-1)(n-2)}{2}}I(L_{i},L_{i+1})L_{i+1},\\ h^{-1}_{L_{i}}(L_{i+1})=L_{i+1}-(-1)^{\frac{(n-1)(n-2)}{2}}I(L_{i},L_{i+1})L_{i}.
\end{gather*}
On the other hand, the group $B_{\widetilde\mu_n}\ltimes (\ZZ/2\ZZ)^{\widetilde\mu_n}$ acts on the set of sequences in $K_0({\rm HMF}_{S_n}^{L_{f_n}}(f_n))$
\[
\{([\E_{1}],\dots, [\E_{\widetilde\mu_n}])\,|\, (\E_{1},\dots, \E_{\widetilde\mu_n})\text{ is a full exceptional collection in }{\rm HMF}_{S_n}^{L_{f_n}}(f_n) \}
\]
by mutations and parity transformations ({cf. \cite[Proposition~2.1]{BP}}):
\begin{gather*}
b_{i}([\E_{1}],\dots, [\E_{\widetilde\mu_n}]):=([\E_{1}], \dots, [\E_{i-1}], [\E_{i+1}], [\RR_{\E_{i+1}}\E_i],[\E_{i+2}],\dots, [\E_{\widetilde\mu_n}]),\\
b^{-1}_{i}([\E_{1}],\dots, [\E_{\widetilde\mu_n}]):=([\E_{1}], \dots, [\E_{i-1}], [\LL_{\E_{i}}\E_{i+1}], [\E_{i}], [\E_{i+2}],\dots, [\E_{\widetilde\mu_n}]),\\
p_{i}([\E_{1}],\dots, [\E_{\widetilde\mu_n}]):=([\E_{1}], \dots, [\E_{i-1}], -[\E_{i}], [\E_{i+1}],\dots, [\E_{\widetilde\mu_n}]),
\end{gather*}
where 
\[
[\RR_{\E_{i+1}}\E_i]=[\E_{i}]-\chi(\E_{i},\E_{i+1})[\E_{i+1}],\quad [\LL_{\E_{i}}\E_{i+1}]=[\E_{i+1}]-\chi(\E_{i},\E_{i+1})[\E_{i}],
\]
and $p_{i}$ is the $i$-th generator of $(\ZZ/2\ZZ)^{\widetilde\mu_n}$.
Therefore, if a distinguished basis of vanishing cycles $(L_{1},\dots, L_{\widetilde\mu_n})$ is homological mirror dual to a full exceptional collection $([\E_1], \dots, [\E_{\widetilde\mu_n}])$, 
then we should have ${\mathbb S}(L_{i},L_{j})=\chi(\E_i,\E_j)$ in order for $B_{\widetilde\mu_n}\ltimes (\ZZ/2\ZZ)^{\widetilde\mu_n}$-actions to be compatible.
\begin{prop}[\cite{G}]\label{prop : transitivity of Braid action}
The group $B_{\widetilde\mu_n}\ltimes (\ZZ/2\ZZ)^{\widetilde \mu_n}$ acts transitively on the set of distinguished bases of vanishing cycles. 
\qed
\end{prop}
Varolgunes shows that there exist cycles $L_1,\dots,L_{\widetilde \mu_n}\in H_{n-1}(\widetilde{f}_n^{-1}(1);\ZZ)=H_{n-1}(Y_{1;0};\ZZ)$, which is isomorphic to $H_{n-1}(Y_{w_0;s};\ZZ)$ by the parallel transport, 
such that 
\[
(-1)^{\frac{n(n-1)}{2}}I(L_{j},L_{i})=\chi^{(n)}_{ij},\quad i<j,
\]
and $(L_{\widetilde \mu_n},\dots, L_1)$ forms a usual (clockwise) distinguished basis (see {\cite[Theorem 1.3]{V}}).
\begin{cor}
There exists an element $B\in B_{\widetilde\mu_n}\ltimes (\ZZ/2\ZZ)^{\widetilde\mu_n}$ such that 
\begin{equation}
B^T S^\phi B=\chi^{(n)},\quad C^\phi B=\left(\frac{1}{(2\pi)^{\frac{n}{2}}}\sum_{\lambda\in I_n}\eta^{\lambda\kappa}\int_{\Gamma^0_j(1)}e^{-\widetilde f_n}\zeta^{(n)}_\lambda \right),
\end{equation}
where $\Gamma^0_j(1)$ is the image of $L_j$ under the isomorphism 
$H_{n-1}(\widetilde{f}_n^{-1}(1);\ZZ)\cong H_n(\CC^n,{\rm Re}(\widetilde f_n)\gg0;\ZZ)$.
\qed
\end{cor}
\begin{rem}
$S:=B^T S^0 B$ and $C:=C^0 B$ are also the Stokes matrix and the central connection matrix of the Frobenius manifold with respect to some chosen data.
\end{rem}
It is natural to expect that the cycles $\{L^{(n)}_1,\dots,L^{(n)}_{\widetilde \mu_n}\}$ on $H_{n-1}(\widetilde f^{-1}_n(1);\ZZ)$ constructed in Section \ref{sec : proof of integral structure} can be obtained from those given by Varolgunes by the action of an automorphism of $\widetilde f_n^{-1}(1)$ and $B_{\widetilde\mu_n}\ltimes (\ZZ/2\ZZ)^{\widetilde \mu_n}$.
In particular, we expect that a set of cycles $\{L^{(n)}_1,\dots,L^{(n)}_{\widetilde \mu_n}\}$ forms a distinguished basis of vanishing cycles.
In fact, it is obvious for the case of $n=1$ (see Lemma \ref{lem : integral structure for n=1}).

\end{document}